\documentclass[12pt]{article}

\usepackage{arydshln}
\usepackage{mathrsfs}
\usepackage{amsmath}
\usepackage{amsfonts}
\usepackage{amssymb}
\usepackage[all,cmtip]{xy}
\usepackage{hyperref}

\usepackage{ulem}
\usepackage{eufrak}

\usepackage{amsmath, amsthm, amsfonts, amssymb}
\def\<{\langle}
\def\>{\rangle}

\def\p{\partial}
\def\a{\alpha}
\def\wh{\widehat}
\def\wt{\widetilde}

\def\-{\overline}

\def\ov{\overline}

\def\d{\delta}
\def\endpf{\hbox{\vrule height1.5ex width.5em}}

\def\d{\delta}
\def\d{\delta}
\def\b{\beta}
\def\a{\alpha}
\def\endpf{\hbox{\vrule height1.5ex width.5em}}

\def\CC{{\mathbb C}}

\def\-{\overline}

\def\ld{\lambda}

\def\t{\tau}

\def\wt{\widetilde}

\def\p{\partial}

\def\d{\delta}

\def\endpf{\hbox{\vrule height1.5ex width.5em}}

\def\a{\alpha}
\def\d{\delta}

\def\d{\delta}
\def\b{\beta}

\def\a{\alpha}

\def\endpf{\hbox{\vrule height1.5ex width.5em}}
\def\a{\alpha}

\def\CC{\mathbb C}

\def\p{\partial}
\def\wt{\widetilde}
\def\-{\overline}

\newtheorem{theorem}{Theorem}[section]
\newtheorem{lemma}[theorem]{Lemma}
\newtheorem{corollary}[theorem]{Corollary}
\newtheorem{proposition}[theorem]{Proposition}

\newtheorem{definition}[theorem]{Definition}

\newtheorem{remark}[theorem]{Remark}

\date{}

\begin{document}

\medskip
\title{\bf  Flattening  a non-degenerate   CR singular
point of real codimension two  }
%\title{\bf Flattening of CR singular
%points of real codimension two}

\medskip
\author {Hanlong Fang and  Xiaojun Huang \footnote{Supported in part by  NSF-1363418.}}

%\centerline{ Dedicated to the memory of  Professor Qi-keng Lu }%
\medskip
%\centerline{ Dedicated to the memory of Professor Qi-Keng Lu}

\vspace{3cm} \maketitle

\centerline{\it Dedicated to Ngaiming Mok on the occasion of his
60th birthday }
%\medskip
%\centerline{ Dedicated to the memory of Professor Qi-Keng Lu}

%\vspace{3cm} \maketitle

\bigskip\bigskip
{\bf Abstract}: This paper continues the previous studies in two
papers of Huang-Yin [HY3-4] on the flattening problem of a CR
singular point of real codimension two sitting in a submanifold in
${\mathbb C}^{n+1}$ with $n+1\ge 3$, whose CR points are
non-minimal. Partially based on the geometric approach initiated in
[HY3] and a formal theory approach used in [HY4], we are able to
provide  a very general flattening theorem for a non-degenerate CR
singular point. As an application, we provide a solution to the
local  complex Plateau problem and  obtain the analyticity of the
local hull of holomorphy near a real analytic definite CR singular
point in a  general setting.

\tableofcontents
\section{Introduction} Let $M\subset {\mathbb C}^{n+1}$ be  a smooth real submanifold. For
a point $q\in M$, there is an immediate first order holomorphic
invariant $Rk(q):=\hbox{dim}_CT^{(1,0)}_qM$ for the germ of $M$ at
$q$. $Rk(q)$ is an upper semi-continuous function. When $Rk(q)$ is
constant for $q\approx p\in M$, we call $p$ a CR point of $M$.
Otherwise, $p$ is called a CR singular point. The study of the
geometric and analytic properties for $M$ near a CR singular point
has attracted considerable attentions since the celebrated paper of
Bishop in 1965 [Bis]. Bishop considered the case when $M$ is a real
$(n+1)$-manifold in ${\mathbb C}^{n+1}$ with a CR singular point at
$p$ and with $Rk(p)=1$. Bishop discovered that under a certain
natural non-degeneracy assumption and a certain holomorphically
invariant convexity of $M$ near $p$, $M$ has a non-trivial local
hull of holomorphy $\widehat{M}$ and has a very rich holomorphically
invariant geometric structure. Bishop conjectured that $\widehat{M}$
is a Levi-flat submanifold which has more or less the same
regularity as $M$ does even up to $M$ near $p$. Bishop's problem was
confirmed in a sequence of papers by Kenig-Webster [KW1-2],
Moser-Webster [MW], Huang-Krantz [HK], and finally in Huang [Hu1].
The global version of the Bishop problem was investigated in the
work of Bedford-Gaveau [BG] and Bedford-Klingenberg [BK]. Other
closely related work at least includes the papers by Gong [Gong1-2],
Gong-Lebl [GL], Gong-Stolotvich [GS1-2], Lebl [Leb],  Burcea
[Bur1-2], Coffman [Co], Huang-Yin [HY1-2], Lebl-Noell-Ravisankar
[LNR1-2], etc.

In 2010,  Dolbeault, Tomassini and Zaitsev [DTZ1-2] initiated the
study of  the generalized Bishop problem for a real codimension two
submanifold $M\subset {\mathbb C}^{n+1}$ with $n+1\ge 3$. In this
setting,  the CR singularity must have complex codimension one and
CR points have CR dimension $n-1\ge 1$. For $M$ to bound a Levi-flat
submanifold $\widehat{M}$, a CR point must be a CR non-minimal
point, namely, for each CR point $q\in M$, there is a proper CR
submanifold  in ${\mathbb C}^{n+1}$ passing through $q$, that is
contained in $M$ and has  CR dimension $n-1$. (This CR submanifold
is in fact the intersection of a leaf of the Levi-foliation in
$\widehat{M}$ with $M$). Moreover, more restricted geometric
assumptions have to be imposed at the  CR singular points. A
solution in certain cases for this generalized Bishop problem was
obtained in Huang-Yin [HY3].
% though the question in the smooth category still remains open.
In a subsequent paper of
Huang-Yin [HY4], a formal version of problems similar to that
considered in [HY3] has been studied.

To be more detailed, we assume that $p\in M$ is a CR singular point
of $M$. We write $(z_1,\cdots,z_n,w)$ for the coordinates of
${\mathbb C}^{n+1}$. After a holomorphic change of coordinates, we
assume that $p=0,\ T^{(1,0)}_pM=\{w=0\}$. Then $M$ near $p=0$ is the
graph of a function of the form:
\begin{equation}\label{001}
 w=F(z,\-{z})=q^{(2)}(z,\ov{z})+o(|z|^2),
 \end{equation}
where $q^{(2)}(z,\ov{z})$ is a polynomial of degree two in
$(z,\ov{z})$. In the classical  case, namely,  $n+1=2$, after a
holomorphic change of variables, we can always make
$q^{(2)}(z,\ov{z})$-real-valued. However, this is no longer the case
for $n+1\ge 3$. Indeed, after a simple holomorphic change of
coordinates, if needed, we can  write
\begin{equation} \label {0011}
q^{(2)}(z,\ov{z})=2\Re{(z\cdot \mathcal A\cdot z^t)}+z\cdot \mathcal B\cdot \ov{z}^t
\end{equation}
with $\mathcal A$ and $ \mathcal B$ being two $(n\times
n)$-matrices. Suppose that $z=\wt{z}\cdot
P+\vec{a}\wt{w}+O(|(z,w)|^2);\ w=\mu \wt{w}+z\cdot b^t+
O(|(z,w)|^2)$  is a holomorphic transformation preserving the form
as in (\ref{001}) and (\ref{0011}). Then $b=0,\mu\not =0$, and $P$
is an $(n\times n)$-invertible matrix. Moreover, if $(M,0)$ is
defined in the new coordinates by
\begin{equation}\label{002}
\wt{w}=\wt{q}^{2}(\wt{z},\ov{\wt{z}})+o(|\wt{z}|^2),
 \end{equation}
 with $\wt{q}^{(2)}(\wt{z},\ov{\wt{z}})=2\Re{\wt{z}\cdot \wt{\mathcal A}\cdot
\wt{z}^t}+\wt{z}\cdot \wt{\mathcal B}\cdot \ov{\wt{z}}^t.$ Then

\begin{equation}\label{00111}
\wt{\mathcal B}=\frac{1}{\mu} P\cdot \mathcal B\cdot \ov{P}^t,\ \
\wt{\mathcal A}=\frac{1}{\ov\mu} P\cdot\mathcal A\cdot {P}^t.
\end{equation}

 When there do not exist
a $\mu\not =0$ and an invertible $P$ such that $\frac{1}{\mu} P\cdot
\mathcal B\cdot \ov{P}^t$ is Hermitian, one can never make
$\wt{q}^{(2)}(z,\ov{z})$ real-valued. Also notice that the
non-degeneracy of the matrix ${\mathcal B}$ is a holomorphic
invariant property.
%We say $M$ has a non-degenerate CR singularity
%at $p=0$ if $\det{\mathcal B}\not =0$.
More general, we make   the following definition:

\begin{definition} (A).  $M$ is said to have a non-degenerate CR singularity at $p$ if  there is a
holomorphic change of variables such that in the new coordinates,
$p=0$, $M$ is defined by an equation of the form as in (\ref{001})
and (\ref{0011}) with $\det{\mathcal B}\not =0$. If there is a
holomorphic change of variables such that ${\mathcal B}$ is a
definite Hermitian matrix, we call $p$ a definite  CR singular point
of $M$. (B). Let $M$ be a real-codimension two real submanifold with
$p\in M$ a CR singular point. We say that $M$ is quadratically
flattenable if there is a change of coordinates such that in the new
coordinates, $p=0$, $M$ near $p=0$ is defined by an equation of the
form as in (\ref{001}) with $q^{(2)}(z,\ov{z})$ real-valued. One
says that $M$ can be holomorphically flattened at $p$ if there is a
holomorphic change of variables such that in the new coordinates,
$p=0$, $M$ is defined by an equation of the form as in (\ref{001})
with $\Im\left(F(z,\-{z})\right)\equiv 0.$
\end{definition}

In the above mentioned work of Dolbeault-Tomassini-Zaitsev [DTZ1-2]
and Huang-Yin [HY3-4], the starting point is to define a generalized
notion of the  Bishop non-degeneracy and generalized Bishop
invariants at a CR singular point $p$. For that purpose, one needs
to assume that $M$ near $p$ is quadratically flattenable. However,
in the setting considered in [DTZ1-2] and [HY3-4], $M$ is always CR
non-minimal at its CR points. This raises a natural question ([Zat])
to understand the implication of CR non-minimality to the quadratic
flattenabilty of $M$  near a CR singular point. The first part of
the  present paper makes an effort along these lines and  we obtain
the following, which, together with  Theorem \ref{04} and the
examples there,  more or less gives an answer to a question of
Zaitsev [Zat]:

\begin{theorem}\label{003} Let $M$ be a real codimension two smooth submanifold in
${\mathbb C}^{n+1}$ with $p\in M$ a  non-degenerate  CR singular
point.
% Suppose that $M$ near $p=0$ is  defined by an equation as in
%(\ref{001}) and (\ref{0011}).
Assume that $M$ is CR non-minimal at its CR points near $p$. Then
$M$ is quadratically flattenable.
% if $\mathcal B$ is non-degenerate,
%namely, $\hbox{det}(\mathcal B)\not =0$.
\end{theorem}

%\begin{remark} Notice that by (\ref{00111}), the non-degeneracy of $\mathcal B$ is an
%invariant property. Also, by an example presented in [Hu?], Theorem
%\ref{003} fails when the non-degeneracy assumption for $\mathcal B$
%is dropped.
%\end{remark}
In the case of $n+1=3$, we will give a much more detailed result
(see Theorem \ref{04}) in terms of the normal form for the pair
$\{{\mathcal A},{\mathcal B}\}$ given by Coffman in [Co], even if
$\mathcal B$ is degenerate. This result, due to its technical
nature,  will be stated as Theorem \ref{04} in $\S 3$.

\medskip

Let $M\subset {\mathbb C}^{n+1}$ be a codimension two real
submanifold with $p\in M$ a non-degenerate CR singular point. Also
assume that after a holomorphic change of variables, $p=0$ and $M$
is defined as in (\ref{001}) (\ref{0011}) with ${\mathcal B}$ a
Hermitian matrix. When ${\mathcal B}$ is  definite, then by the
classical Takagi theorem [HY2] [HK], we can further make
$q^{(2)}(z,\-{z})=\sum_{j=1}^{n}\left(|z_j|^2+\lambda_j(z_j^2+\-{z_j^2})\right),$
where $0\le \lambda_1\le\cdots, \lambda_n<\infty.$ The set
$\{\lambda_1,\cdots,\lambda_n\}$ is called the set of generalized
Bishop invariants of $M$ at the CR singular point. $\lambda_j$ is
called an elliptic, parabolic or hyperbolic Bishop invariant, if
$0\le \lambda_j<\frac{1}{2}$,  $ \lambda_j=\frac{1}{2}$ or
$\lambda_j>\frac{1}{2}$. This terminology coincides with the
classical definition of Bishop when $n+1=2$ ([Bis] [KW1] [MW] [HK]).
However, when ${\mathcal B}$ is not a definite matrix, we can not,
in general, simultaneously diagonalize ${\mathcal A}$ and ${\mathcal
B}$. In the case of $n+1=3$, Coffman gave a list of the forms that
the pair $\{{\mathcal A},{\mathcal B}\}$ can be transformed to. (See
the list given at the beginning of $\S 3$). Two cases in his list
are geometrically quite special, in which the corresponding
quadratic term takes one of the following forms after a holomorphic
change of coordinates:
\begin{equation}
\label{P1}
q^{(2)}=|z_1|^2+|z_2|^2+\frac{1}{2}(z_1^2+\-{z_1^2})+\frac{1}{2}(z_2^2+\-{z_2^2})\
\hbox{or}
\end{equation}
\begin{equation} \label{M1}
q^{(2)}=|z_1|^2-|z_2|^2+\lambda(z_1^2+\-{z_1^2})+\lambda(z_2^2+\-{z_2^2}),\
\ \lambda\ge \frac{1}{2}.
\end{equation}
In the case of (\ref{P1}), the two generalized Bishop invariants of
the CR singular point at the origin are both  parabolic.  A
consequence of this is that  the set of CR singular points may have
real dimension $n=2$, which does create a lot of problems for the
geometric studies of $M$ near $0$.

To explain the speciality of (\ref{M1}), we recall a definition from
[HY3]: Let $(M,p)$ be a codimension two real submanifold in
${\mathbb C}^{n+1}$ with $p\in M$ a CR singular point. We say
$(M,p)$ possesses an elliptic complex tangent direction if there is
an affine complex plane ${\mathcal H}$ that passes through $p$ and
is transversal to the complex tangent space of $M$ at $p$ such that
$M\cap {\mathcal H}$ is an elliptic Bishop surface inside ${\mathcal
H}$ in the classical sense ([Bis]). (See Definition \ref{0.0}). Now,
a simple  algebraic computation  shows that a codimension two real
submanifold $M\subset {\mathbb C}^{n+1}$ with a non-degenerate
quadratically flattenable CR singular point at $p$ has no elliptic
directions at $p$ if and only if $n+1=3$ and after a holomorphic
change of variables sending $p$ to $0$,  $M$ near $p=0$ is defined
by an equation of the form as in (\ref{001}) with $q^{(2)}$ being
given by (\ref{M1}). (See the paper by Lebl-Noell-Ravisankar[LNR]).

% In (\ref{P1}), two generalized Bishop invariants are
%parabolic and the CR singular set ${\mathcal S}$ of $M$  near $p=0$
%may have the largest possible dimension: real dimension two. In
%(\ref{M1}), $(M,0)$ does not have any elliptic direction. Here, we
%say $M\subset {\mathbb C}^{n+1}$ defined by an equation of the form
%as in (\ref{001}) has an elliptic direction $c=(c_1,\cdots,c_n)\not
%= 0$ at its CR singular point $0$, if the the real surface $M\cap
%\{(z_1,\cdots,z_n)=\xi c,\ \xi\in {\mathbb C}\}$ contained in the
%complex 2-subspace defined by $\{(z_1,\cdots, z_n,w)\in {\mathbb
%C}^{n+1}: \ (z_1,\cdots, z_n)=\xi c\}$ is a Bishop surface with an
%elliptic singularity at $0$ in the classical case [Bis] [HK]. In all
%other cases with a quadratically flattenable non-degenerate CR
%singular point at $0$, $(M,0)$ has an elliptic direction and the
%dimension of its CR singular set near $0$ has real dimension bounded
%by  $n$.

A major part of the paper continues the study in [HY3] and [HY4],
which is devoted to the understanding of the holomorphically
flattening problem near a CR singular point when  $M$ is real
analytic. This problem is equivalent to  finding a real analytic
Levi-flat hypersurface with $M$ as part of its real analytic
boundary and with leaves moving along a transversal direction to the
complex tangent space of the CR singular point. This   has an
immediate application to the  study of the precise description of
the local hull of holomorphy of $M$.
%by making use Theorem \ref{003} and
%partially applying the work   from Huang-Yin [HY3-4],
Our purpose is to provide the following  general holomorphic
flattening theorem:

\begin{theorem}\label{005} Let $M$ be a real analytic real-codimension two submanifold in
${\mathbb C}^{n+1}$ with $n\ge 2$ and with  $p\in M$ a
non-degenerate CR singular point.
% Suppose that $M$ near $p=0$ is  defined by an equation as in
%(\ref{001}) and (\ref{0011}).
Assume that $M$ is CR non-minimal at its CR points near $p$. Then
$(M,p)$ can be holomorphically  flattened if $M$ has an elliptic
direction at $p$. More precisely, $(M,p)$ can be holomorphically
flattened if either $n+1\ge 4$ or  $n+1=3$ but $(M,p)$ is not
holomorphically equivalent to a submanifold $(M',0)$ whose quadratic
term takes the form  in (\ref{M1}).
% if $\mathcal B$ is non-degenerate,
%namely, $\hbox{det}(\mathcal B)\not =0$.
\end{theorem}

\begin{corollary} \label{44.44}
Let $M$ be a real analytic real-codimension two submanifold in
${\mathbb C}^{n+1}$ with $n\ge 2$ and with  $p\in M$ a
non-degenerate CR singular point.
% Suppose that $M$ near $p=0$ is  defined by an equation as in
%(\ref{001}) and (\ref{0011}).
Assume that $M$ is CR non-minimal at its CR points near $p$. Assume
that either $n+1\ge 4$ or  $n+1=3$ but $(M,p)$ is not
holomorphically equivalent to a submanifold $(M',0)$ whose quadratic
term takes the form  in (\ref{M1}). Then there is  a real analytic
Levi-flat hypersurface $\wh{M}$, which has $M$ near $p$ as part of
its real analytic boundary and is foliated by complex hypersurfaces
shrinking down to $p$ along the normal direction of $M$ in $\wh{M}$
at $p$. Moreover, when $p$ is a definite CR singular point, then
there is a small $\epsilon_0>0$ such that for any
$0<\epsilon<<\epsilon_0$, $\wh{M}\cap B_p(\epsilon)$ is a connected
open
 piece containing the origin  of  the hull of holomorphy of $M\cap B_p(\epsilon_0)$. Here
$B_p(\epsilon)$ denotes the ball  centered at $p$ of radius
$\epsilon$.

\end{corollary}

% Let $M\subset
%{\mathbb C}^3$ be a codimension two real submanifold with $p\in M$ a
%non-degenerate CR singular point. We say that $(M,p)$ is of an
%exceptional class, if after a holomorphic change of variables, $M$
%is defined as in (\ref{001}) with
Theorem \ref{005} is contained in Huang-Yin [HY3] when $p$ is a
definite CR singular point with one of the generalized Bishop
invariants elliptic.

 We next say a few words about the organization
of the paper. In $\S 2$, we review a fundamental identity for
$n+1=3$ from the non-minimality condition first obtained in [HY3].
In $\S 3$, we first give a list of the normal form for the quadratic
terms when $n+1=3$. We then state Theorem \ref{04} which gives an
understanding about the quadratic flattening problem when $n+1=3$.
We also show that the result in Theorem \ref{04} is optimal by
presenting several examples. In $\S 4$, we give a proof of Theorem
\ref{04} by making an extensive use of the  identity discussed in
$\S 2$. In $\S 5$, we give a proof of Theorem \ref{003}. $\S 6-\S7$
will be devoted to the proof of Theorem \ref{005}. In $\S 6$, we use
a geometric argument initiated in [HY3] to approach Theorem
\ref{005}. The nice feature for this approach is that we do not need
to know much about the quadratic normal form, which is almost
impossible to obtain when $n+1>3$. However this argument, though
very general, needs the real dimension of the real analytic set of
CR singular points of $M$  is no  more than $2n-2$. (And there is an
example in Remark \ref{22.221} showing that this approach fails when
the singular set is too large). This excludes the case when $n+1=3$
and the quadratic term of the defining equation of the manifold
takes the normal form as in (\ref{P1}). In $\S 7$, we give a proof
of Theorem \ref{005} in this exceptional case. A good thing about
this exceptional  case is that the associated quadratic term has the
simplest possible symmetric form. This makes the formal argument
developed in [HY4] to be very much adaptable to this setting.
Indeed, we can first prove that in this case $M$ can always be
formally flattened by a special form of holomorphic transformations,
which together with the Huang-Krantz construction of holomorphic
disks give a convergent flattening as in the  other cases considered
in [HY4].

\bigskip
\section{Implication of integrability conditions}

 Let $(M,0)$ be a smooth submanifold of codimension two in
${\mathbb{C} }^{3}$ with $0\in M$ as a CR singular point. Assume
that all CR points are non-minimal with CR dimension one. Use
$(z,w)=(z_1,z_2,w)$ for the coordinates of ${\mathbb C}^3$. Assume
that, after a holomorphic change of coordinates, $M$ near $0$ is
defined by an equation of the form:
%We also assume that $T^{(1,0)}_0M=\{w=0\}$ and $0$ is not completely
%degenerate flat CR singular point(see [Sto],[DTZ] and [HY]). Then
%$M$ can be defined by a real analytic equation of the form (see [HY]):
\begin{equation}
w=q^{(2)}(z,\ov{z})+p(z,\ov{z})+iE(z,\ov{z}), \label{429eq1}
\end{equation}
where $q^{(2)}(z,\ov{z})=2\Re{({z}\cdot {\mathcal A}\cdot
{z}^t)}+{z}\cdot {\mathcal B}\cdot \ov{{z}}^t$,
$p(z,\ov{z}),E(z,\ov{z})=O(|z|^3)$ and both $p(z,\ov{z})$ and
$E(z,\ov{z})$ are real-valued smooth functions. For convenience of
notation, we also write
\begin{equation}
F(z,\ov{z})=p(z,\ov{z})+iE(z,\ov{z})\ \text{and}\ G(z,\ov{z})=
q^{(2)}(z,\ov{z})+p(z,\ov{z}). \label{82eq1}
\end{equation}
Then we have
$$
w=q^{(2)}(z,\ov{z})+F(z,\ov{z})=G(z,\ov{z})+iE(z,\ov{z}).
$$

%First, a point $(z,w)\in M$ is a CR singular point it is easy to
%verify that the set ${\cal S}$ of CR singular points of $M$ near $0$
%is defined by the foll
% a thin
%set in $M$.
In what follows, we write $ \chi_{\a}=\frac{\p \chi}{\p z_\a},\
\chi_{\ov{\a}}=\frac{\p \chi}{\p \ov{z}_{\a}}$ with $ \a=1,2$ for a
smooth function $\chi(z,\ov{z})$ in $z$.  We define
\begin{equation}\begin{split}\label{325eq2}
L:=&(G_2-iE_2)\frac{\p}{\p z_1}-(G_1-iE_1)\frac{\p}{\p
z_2}+2i(G_2E_1-G_1E_2)\frac{\p}{\p w}\\
=&A\frac{\p}{\p z_1}-B\frac{\p}{\p z_2}+C\frac{\p}{\p w}.
\end{split}\end{equation}
%Then we have
%\begin{equation*}\begin{split}
%L(-w+G+iE)&=(G_2-iE_2)(G_1+iE_1)-(G_1-iE_1)(G_2+iE_2)\\
%&\ \ -2i(G_2E_1-G_1E_2)=0,\\
%L(\ov{-w+G+iE})&=(G_2-iE_2)(G_1-iE_1)-(G_1-iE_1)(G_2-iE_2) =0.
%\end{split}\end{equation*}
Then $L$ is a complex tangent vector field of type $(1,0)$ along $M$
near $0$. (See [$\S 2$, HY4]).  Moreover, a straightforward
computation shows that
\begin{equation*}\begin{split}
T:=[L,\ov{L}]=&\left[A\frac{\p}{\p z_1}-B\frac{\p}{\p
z_2}+C\frac{\p}{\p w},\ov{A}\frac{\p}{\p
\ov{z_1}}-\ov{B}\frac{\p}{\p
\ov{z_2}}+\ov{C}\frac{\p}{\p \ov{w}}\right]\\
 =&\lambda_{(1)}\frac{\p}{\p
\ov{z_1}}+\lambda_{(2)}\frac{\p}{\p
\ov{z_2}}+\lambda_{(3)}\frac{\p}{\p
\ov{w}}+\lambda_{(4)}\frac{\p}{\p z_1}+\lambda_{(5)}\frac{\p}{\p
z_2}+\lambda_{(6)}\frac{\p}{\p w},
\end{split}\end{equation*}
where
\begin{equation}\begin{array}{lll}   \label{325eq3}
&\lambda_{(1)}=A\cdot (\ov{A})_1-B\cdot(\ov{A})_2,& \ \lambda_{(4)}
=-\ov{A}\cdot A_{\ov{1}}+\ov{B}\cdot{A}_{\ov{2}},\\
&\lambda_{(2)}=-A\cdot (\ov{B})_1+B\cdot(\ov{B})_2,&\
\lambda_{(5)}=\ov{A}\cdot {B}
_{\ov{1}}-\ov{B}\cdot{B}_{\ov{2}},\\
&\lambda_{(3)}=A\cdot (\ov{C})_1-B\cdot(\ov{C})_2,&\
\lambda_{(6)}=-\ov{A}\cdot {C}_{\ov{1}}+\ov{B}\cdot{C}_{\ov{2}}.
\end{array}\end{equation}
Notice that
\begin{equation*}\label{57eq8}
\lambda_{(1)}=-\ov{\lambda_{(4)}},\
\lambda_{(2)}=-\ov{\lambda_{(5)}},\
\lambda_{(3)}=-\ov{\lambda_{(6)}}.
\end{equation*}

Write $[L,T]=:\Gamma_{(1)}\frac{\p}{\p \ov
z_1}+\Gamma_{(2)}\frac{\p}{\p \ov z_2}+\Gamma_{(3)}\frac{\p}{\p \ov
\omega}+\Gamma_{(4)}\frac{\p}{\p z_1}+\Gamma_{(5)}\frac{\p}{\p
z_2}+\Gamma_{(6)}\frac{\p}{\p \omega}.$  By a direct computation as
in  [HY4],  we have the following explicit expressions for
$\Gamma_{(1)},\cdots,\Gamma_{(6)}$:
\begin{equation}\begin{array}{lll}\label{325eq4}
&\Gamma_{(1)}=A\cdot{(\lambda_{(1)})}{_1}-B\cdot{(\lambda_{(1)})}{_2},\,\,\,\Gamma_{(2)}=A\cdot{(\lambda_{(2)})}{_1}-B\cdot{(\lambda_{(2)})}{_2},\,\,\,\Gamma_{(3)}=A\cdot{(\lambda_{(3)})}{_1}-B\cdot{(\lambda_{(3)})}{_2},\\
&\Gamma_{(4)}=A\cdot{(\lambda_{(4)})}{_1}-B\cdot{(\lambda_{(4)})}{_2}-{\lambda_{(1)}}\cdot{ A}_{\ov{ 1}}-{\lambda_{(2)}}\cdot{ A}_ {\ov {2}}-{\lambda_{(4)}}\cdot{ A}{_1}-{\lambda_{(5)}}\cdot{ A}{_2},\\
&\Gamma_{(5)}=A\cdot{(\lambda_{(5)})}{_1}-B\cdot{(\lambda_{(5)})}{_2}+{\lambda_{(1)}}\cdot{ B}_ {\ov {1}}+{\lambda_{(2)}}\cdot{ B}_ {\ov {2}}+{\lambda_{(4)}}\cdot{ B}{  _1}+{\lambda_{(5)}}\cdot{ B}{_2},\\
&\Gamma_{(6)}=A\cdot{(\lambda_{(6)})}{_1}-B\cdot{(\lambda_{(6)})}{_2}-{\lambda_{(1)}}\cdot{ C}_ {\ov {1}}-{\lambda_{(2)}}\cdot{ C}_ {\ov {2}}-{\lambda_{(4)}}\cdot{C}{_1}-{\lambda_{(5)}}\cdot{C}{_2}.
\end{array}\end{equation}
Since we assumed that all CR points are  non-minimal with CR
dimension one,  $[L,T]$ is spanned by $\{L,\ov L,T\}$ over a dense
subset of $M$ near $0$. Hence, at these  points, there are complex
numbers $k,\sigma,\tau$, such that
$$[L,T]=kL+\sigma\ov L+\tau T.$$

Comparing the coefficients of both sides, we have as in [HY4]:
\begin{equation}\label{1} \Gamma_{(1)}=\sigma \ov A+\tau\lambda_{(1)},
\end{equation}
\begin{equation} \label{2}\Gamma_{(2)}=-\sigma \ov B+\tau\lambda_{(2)},
\end{equation}
\begin{equation}\label{3} \Gamma_{(4)}=kA+\tau\lambda_{(4)},
\end{equation}
\begin{equation}\label{4} \Gamma_{(5)}=-kB+\tau\lambda_{(5)}.
\end{equation}

 An elementary algebraic computation gives the following equality by eliminating $\sigma$ from $(\ref{1})$ and  $(\ref{2})$:
 \begin{equation}\label{5*}
    \ov B\Gamma_{(1)}+\ov A\Gamma_{(2)}=\ov B(\sigma \ov A+\tau\lambda_{(1)})+\ov A(-\sigma \ov B+\tau\lambda_{(2)})=\tau(\lambda_{(1)}\ov B+\lambda_{(2)}\ov A).
 \end{equation}

Similarly we derive the following equality by eliminating $k$ from $(\ref{3})$ and  $(\ref{4})$:
  \begin{equation}\label{6}
    B\Gamma_{(4)}+A\Gamma_{(5)}= B(kA+\tau\lambda_{(4)})+ A(-kB+\tau\lambda_{(5)})=\tau(\lambda_{(4)} B+\lambda_{(5)} A).
  \end{equation}

Combining  $(\ref{5*})$ and  $(\ref{6})$, we obtain an identity in
an open dense subset whose closure contains $0\in {\mathbb C}^2$,
which will be fundamentally used for the proof of Theorem \ref{003}:
(Hence, in the real analytic case, the identity holds in a
neighborhood of $0$. And in the smooth category, it holds as a germ
at $0$)
\begin{equation}\label{core}
(\ov B\Gamma_{(1)}+\ov A\Gamma_{(2)})(\lambda_{(4)} B+\lambda_{(5)}
A)= (B\Gamma_{(4)}+A\Gamma_{(5)})(\lambda_{(1)}\ov
B+\lambda_{(2)}\ov A).
\end{equation}
For  convenience,  we introduce the following notation:
\begin{equation}\label{7}
X_1:=\ov B\Gamma_{(1)}+\ov A\Gamma_{(2)},\,\,X_2:=\lambda_{(4)} B+\lambda_{(5)} A,\,\,Y_1:=B\Gamma_{(4)}+A\Gamma_{(5)},\,\,Y_2:=\lambda_{(1)}\ov B+\lambda_{(2)}\ov A.
\end{equation}
Then (\ref{core}) can be   written as:
\begin{equation}\label{8}X_1X_2=Y_1Y_2.\end{equation}

\section{Quadratic flattening  in the 3-dimensional case}
Let $M$ be  defined by an equation of the form in (\ref{429eq1}).
Recall that $q(z,\ov{z})=2\Re({{z}\cdot {\mathcal A}\cdot
{z}^t})+{z}\cdot {\mathcal B}\cdot \ov{{z}}^t$ is the quadratic term
in the defining function $(\ref{429eq1})$. When ${\mathcal B}=0$,
$M$ near $0$ is already quadratically flattened. Hence, we assume
that ${\mathcal B}\not = 0$. Then the CR points must have CR
dimension one.  By a result of Coffman in [Co], after a
 holomorphic change of coordinates, the pair $\{{\mathcal A}, {\mathcal
B}\}$ can be transformed  into one of the  forms $\{{\mathcal A}',
{\mathcal B}'\}$  listed below. Namely, there is a biholomorphic
change of coordinates such that in the new coordinates $(z',w')$,
$M$ is defined by an equation of the form: $w'=2\Re({{z'}\cdot
{\mathcal A}'\cdot {z'}^t})+{z'}\cdot {\mathcal B}'\cdot
\ov{{z}'}^t+o(|z'|^2)$ where  the pair $\{{\mathcal A}', {\mathcal
B}'\}$  takes one of the following  forms:

% to be $\mathcal A^\prime,\mathcal B^\prime$ as
%follows, respectively:
\begin{equation}\label{1a}
(1a). \,\,\,\,\, \mathcal B^\prime=\left(\begin{matrix}1&0\\
0&e^{i\theta}
\end{matrix}\right),\,\, 0<\theta<\pi;\,\,\,  \mathcal A^\prime=\left(\begin{matrix}a&b\\b&d
\end{matrix}\right),\,\,a>0,d>0.\quad\quad\quad\quad\quad\quad\quad\quad\quad\quad
\end{equation}
\begin{equation}\label{1b}
 (1b).\,\, \,\,\, \mathcal B^\prime=\left(\begin{matrix}1&0\\
 0&e^{i\theta}
 \end{matrix}\right),\,\, 0<\theta<\pi;\,\,\,  \mathcal A^\prime=\left(\begin{matrix}0&b\\b&d
 \end{matrix}\right),\,\,b\geq 0,d\geq 0.\quad\quad\quad\quad\quad\quad\quad\quad\quad\quad
\end{equation}
\begin{equation}\label{1c}
(1c). \,\,\, \,\,\mathcal B^\prime=\left(\begin{matrix}1&0\\
0&e^{i\theta}
\end{matrix}\right),\,\, 0<\theta<\pi;\,\,\,  \mathcal A^\prime=\left(\begin{matrix}a&b\\b&0
\end{matrix}\right),\,\,a>0,b\geq 0.\quad\quad\quad\quad\quad\quad\quad\quad\quad\quad
\end{equation}
\begin{equation}\label{2a}
(2a). \,\,\,\,\, \mathcal B^\prime=\left(\begin{matrix}0&1\\
\tau&0
\end{matrix}\right),\,\,\,\,\,\, 0<\tau<1;\,\,\,  \mathcal A^\prime=\left(\begin{matrix}a&b\\b&d
\end{matrix}\right),\,\,b> 0,|a|=\frac{1}{2}.\quad\quad\quad
\quad\quad\quad\quad\quad\quad
\end{equation}
\begin{equation}\label{2b}
(2b). \,\,\,\,\, \mathcal B^\prime=\left(\begin{matrix}0&1\\
\tau&0
\end{matrix}\right),\,\,\,\,\,\, 0<\tau<1;\,\,\,  \mathcal A^\prime=\left(\begin{matrix}0&b\\b&d
\end{matrix}\right),\,\,b> 0,|d|=\frac{1}{2}.\quad\quad\quad
\quad\quad\quad\quad\quad\quad
\end{equation}
\begin{equation}\label{2c}
(2c). \,\,\,\,\, \mathcal B^\prime=\left(\begin{matrix}0&1\\
 \tau&0
 \end{matrix}\right),\,\,\,\,\,\, 0<\tau<1;\,\,\, \mathcal A^\prime=\left(\begin{matrix}0&b\\b&0
 \end{matrix}\right),\,\,b> 0.\quad\quad\quad
 \quad\quad\quad\quad\quad\quad\quad\quad\quad\,\,\,
\end{equation}
\begin{equation}\label{2d}
(2d).\,\, \,\,\, \mathcal B^\prime=\left(\begin{matrix}0&1\\
\tau&0
\end{matrix}\right),\,\,\,\,\,\, 0<\tau<1;\,\,\,  \mathcal A^\prime=\left(\begin{matrix}\frac{1}{2}&0\\0&d
\end{matrix}\right).\quad\quad\quad
\quad\quad\quad\quad\quad\quad\quad\quad\quad\quad\quad\quad\,\,\,
\end{equation}
\begin{equation}\label{2e}
 (2e). \,\,\,\,\, \mathcal B^\prime=\left(\begin{matrix}0&1\\
 \tau&0
 \end{matrix}\right),\,\,\,\,\,\, 0<\tau<1;\,\,\,  \mathcal A^\prime=\left(\begin{matrix}0&0\\0&\frac{1}{2}
 \end{matrix}\right).
\quad\quad\quad
\quad\quad\quad\quad\quad\quad\quad\quad\quad\quad\quad\quad\,\,\,
\end{equation}
\begin{equation}\label{2f}
 (2f). \,\,\,\,\, \mathcal B^\prime=\left(\begin{matrix}0&1\\
 \tau&0
 \end{matrix}\right),\,\,\,\,\,\, 0<\tau<1;\,\,\,  \mathcal A^\prime=\left(\begin{matrix}0&0\\0&0
 \end{matrix}\right).
\quad\quad\quad
\quad\quad\quad\quad\quad\quad\quad\quad\quad\quad\quad\quad\,\,\,
\end{equation}
\begin{equation}\label{3a}
 (3a). \,\,\,\,\, \mathcal B^\prime=\left(\begin{matrix}0&1\\
 1&i
 \end{matrix}\right);\,\,\,\,\,\,   \mathcal A^\prime=\left(\begin{matrix}a&b\\b&d
 \end{matrix}\right), a>0,\ b\in {\mathbb R}.
\quad\quad\quad
\quad\quad\quad\quad\quad\quad\quad\quad\quad\quad\quad\quad
\end{equation}
\begin{equation}\label{3b}
(3b). \,\,\,\,\, \mathcal B^\prime=\left(\begin{matrix}0&1\\
1&i
\end{matrix}\right);\,\,\,\,\,\,   \mathcal
A^\prime=\left(\begin{matrix}0&b\\b&d
\end{matrix}\right), b>0,\ d\in{\mathbb R}.
\quad\quad\quad
\quad\quad\quad\quad\quad\quad\quad\quad\quad\quad\quad\quad
\end{equation}
\begin{equation}\label{3c}
(3c).\,\, \,\,\, \mathcal B^\prime=\left(\begin{matrix}0&1\\
1&i
\end{matrix}\right);\,\,\,\,\,\,   \mathcal A^\prime=\left(\begin{matrix}0&0\\0&d
\end{matrix}\right), d\geq 0.
\quad\quad\quad
\quad\quad\quad\quad\quad\quad\quad\quad\quad\quad\quad\quad\quad\quad\quad
\end{equation}
\begin{equation}\label{4a}
 (4a). \,\,\,\,\, \mathcal B^\prime=\left(\begin{matrix}0&1\\
 0&0
 \end{matrix}\right);\,\,\,\,\,\,   \mathcal A^\prime=\left(\begin{matrix}a&b\\b&\frac{1}{2}
 \end{matrix}\right), b>0.
\quad\quad\quad
\quad\quad\quad\quad\quad\quad\quad\quad\quad\quad\quad\quad\quad\quad\quad
\end{equation}
\begin{equation}\label{4b}
 (4b). \,\,\,\,\, \mathcal B^\prime=\left(\begin{matrix}0&1\\
 0&0
 \end{matrix}\right);\,\,\,\,\,\,   \mathcal A^\prime=\left(\begin{matrix}\frac{1}{2}&b\\b&0
 \end{matrix}\right), b>0.
\quad\quad\quad
\quad\quad\quad\quad\quad\quad\quad\quad\quad\quad\quad\quad\quad\quad\quad\,
\end{equation}
\begin{equation}\label{4c}
 (4c). \,\,\,\,\, \mathcal B^\prime=\left(\begin{matrix}0&1\\
 0&0
 \end{matrix}\right);\,\,\,\,\,\,   \mathcal A^\prime=\left(\begin{matrix}0&b\\b&0
 \end{matrix}\right), b>0.
\quad\quad\quad
\quad\quad\quad\quad\quad\quad\quad\quad\quad\quad\quad\quad\quad\quad\quad
\end{equation}
\begin{equation}\label{4d}
(4d). \,\,\,\,\, \mathcal B^\prime=\left(\begin{matrix}0&1\\
0&0
\end{matrix}\right);\,\,\,\,\,\,   \mathcal A^\prime=\left(\begin{matrix}a&0\\0&\frac{1}{2}
\end{matrix}\right), a\geq 0.
\quad\quad\quad
\quad\quad\quad\quad\quad\quad\quad\quad\quad\quad\quad\quad\quad\quad\,\,\,
\end{equation}
\begin{equation}\label{4e}
 (4e). \,\,\,\,\, \mathcal B^\prime=\left(\begin{matrix}0&1\\
 0&0
 \end{matrix}\right);\,\,\,\,\,\,   \mathcal A^\prime=\left(\begin{matrix}\frac{1}{2}&0\\0&0
 \end{matrix}\right).
\quad\quad\quad\quad\quad\quad
\quad\quad\quad\quad\quad\quad\quad\quad\quad\quad\quad\quad\quad\quad\,\,\,
\end{equation}
\begin{equation}\label{4f}
(4f). \,\,\,\,\, \mathcal B^\prime=\left(\begin{matrix}0&1\\
0&0
\end{matrix}\right);\,\,\,\,\,\,
\mathcal A^\prime=\left(\begin{matrix}0&0\\0&0
\end{matrix}\right).
\quad\quad\quad\quad\quad\quad
\quad\quad\quad\quad\quad\quad\quad\quad\quad\quad\quad\quad\quad\quad\,
\end{equation}
\begin{equation}\label{5}
(5). \,\,\,\,\,\, \mathcal B^\prime=\left(\begin{matrix}1&0\\
0&1
\end{matrix}\right);\,\,\,\,\,\,
\mathcal A^\prime=\left(\begin{matrix}\lambda_1&0\\0&\lambda_2
\end{matrix}\right), \,\,0\leq\lambda_1\leq\lambda_2.
\quad\quad\quad\quad\quad\quad
\quad\quad\quad\quad\quad\quad\quad\,\,
\end{equation}
\begin{equation}\label{6a}
(6a). \,\,\,\,\, \mathcal B^\prime=\left(\begin{matrix}1&0\\
0&-1
\end{matrix}\right);\,\,\,\,\,\,
\mathcal A^\prime=\left(\begin{matrix}\lambda_1&0\\0&\lambda_2
\end{matrix}\right), \,\,0\leq\lambda_1\leq\lambda_2.
\quad\quad\quad\quad\quad\quad
\quad\quad\quad\quad\quad\quad\,\,
\end{equation}
\begin{equation}\label{6b}
(6b). \,\,\,\,\, \mathcal B^\prime=\left(\begin{matrix}1&0\\
0&-1
\end{matrix}\right);\,\,\,\,\,\,
\mathcal A^\prime=\left(\begin{matrix}0&\lambda\\\lambda&0
\end{matrix}\right), \,\,\lambda>0.
\quad\quad\quad\quad\quad\quad
\quad\quad\quad\quad\quad\quad\quad\quad\quad\quad
\end{equation}
\begin{equation}\label{6c}
(6c). \,\,\,\,\, \mathcal B^\prime=\left(\begin{matrix}1&0\\
0&-1
\end{matrix}\right);\,\,\,\,\,\,
\mathcal A^\prime=\left(\begin{matrix}\frac{1}{2}&\frac{1}{2}\\\frac{1}{2}&\frac{1}{2}
\end{matrix}\right).
\quad\quad\quad\quad\quad\quad
\quad\quad\quad\quad\quad\quad\quad\quad\quad\quad
\quad\quad\quad\,
\end{equation}
\begin{equation}\label{7a}
(7a). \,\,\,\,\, \mathcal B^\prime=\left(\begin{matrix}0&1\\
1&0
\end{matrix}\right);\,\,\,\,\,\,
\mathcal A^\prime=\left(\begin{matrix}0&b\\b&\frac{1}{2}
\end{matrix}\right),\,\, b>0.
\quad\quad\quad\quad\quad\quad
\quad\quad\quad\quad\quad\quad\quad\quad\quad\quad\quad
\end{equation}
\begin{equation}\label{7b}
(7b). \,\,\,\,\, \mathcal B^\prime=\left(\begin{matrix}0&1\\
1&0
\end{matrix}\right);\,\,\,\,\,\,
\mathcal A^\prime=\left(\begin{matrix}\frac{1}{2}&0\\0&d
\end{matrix}\right),\,\, {\rm Im} d>0.
\quad\quad\quad\quad\quad\quad
\quad\quad\quad\quad\quad\quad\quad\quad\quad\quad
\end{equation}
\begin{equation}\label{8-9}
(8). \,\,\,\,\, \mathcal B^\prime=\left(\begin{matrix}1&0\\
0&0
\end{matrix}\right).\,\,\,\,\,\,  (9). \,\,\,\,\, \mathcal B^\prime=\left(\begin{matrix}0&0\\
0&0
\end{matrix}\right).\,\,\,\,\,\,\quad\quad\quad\quad\quad\quad\quad\,\,\,\quad\quad\quad\quad\quad\quad\quad\quad\quad\quad\,
\end{equation}

We remark that we made the change on the  forms in (3a) and (3b) for
the ${\mathcal B}'$ part (compared with  the original list in [Co]).
(See [p950, Co]). In the cases of $(5)$-$(9)$, the manifold is
already quadratically flattened. Therefore  we restrict our
discussion in this and the next sections to the cases in
$(1)$-$(4)$. We will prove the following theorem in $\S 4$.

\begin{theorem}\label{04} Let $M$ be a real codimension two submanifold in
    ${\mathbb C}^{3}$ with $0\in M$ a  CR singular point. Suppose that
    $M$ near $p=0$ is  defined by an equation as in (\ref{429eq1}). Assume that $M$ is CR non-minimal at its CR
    points. Then $M$ can be quadratically flattened, possibly except in the cases when the normal form for the pair of the  matrices $\{\mathcal A,\mathcal B\}$
    %as in defining function (\ref{001}) and
    %(\ref{0011})
    takes   the  following form:
    $$\ (i).\  {\mathcal A^\prime}=\left(\begin{matrix}0&0\\0&0
    \end{matrix}\right),{\mathcal B^\prime}=\left(\begin{matrix}0&1\\
        0&0
    \end{matrix}\right);\rm or\ (ii).\ {\mathcal A^\prime}=\left(\begin{matrix}0&\frac{1}{2}\\ \frac{1}{2}&0
    \end{matrix}\right),{\mathcal B^\prime}=\left(\begin{matrix}0&1\\
    0&0
    \end{matrix}\right);$$ $$\rm or\  (iii).\ {\mathcal A^\prime}=\left(\begin{matrix}0&0\\0&\frac{1}{2}
    \end{matrix}\right),{\mathcal B^\prime}=\left(\begin{matrix}0&1\\
    0&0
    \end{matrix}\right).$$ In particular, if $\mathcal B$ is
    non-degenerate, namely $\hbox{det}(\mathcal B)\not =0$, then $M$ is always quadratically flattenable.
\end{theorem}
    The following examples show that Theorem \ref{04} is optimal.

\noindent {\it Example 3.1:} ([Hu2]):  Let $M\subset\mathbb C^3$
with coordinates $(z_1,z_2,w)$ be defined by $w=z_1\ov
z_2+z_1z_2+\ov z_1\ov z_2.$ It is easy to see that $0\in M$ is a  CR
singular point and that the quadratic term in the defining function
takes the following normal form:
 $$\mathcal A^\prime=\left(\begin{matrix}0&\frac{1}{2}\\ \frac{1}{2}&0\end{matrix}\right),\quad\mathcal B^\prime=\left(\begin{matrix}0&1\\
0&0
\end{matrix}\right).$$
Now, following [Hu2], we verify that $M$ is CR non-minimal at its CR
points.  For a CR point $p_0\in M,$ if $p_0\in\{z_2=w=0\},M$ must be
non-minimal at $p_0$.  Otherwise $ p_0\notin\{z_2=w=0\}$ and we
define $L=(z_1+\ov z_1)\frac{\partial}{\partial
z_1}-z_2\frac{\partial}{\partial z_2}+(\ov z_1 z_2+\ov z_2 z_1+\ov
z_2\ov z_1)\frac{\partial}{\partial w}$.  Write $h=-w+z_1\ov
z_2+z_1z_2+\ov z_1\ov z_2.$  Then we have $L(h)=L(\ov h)=0$ along
$M.$ Hence $L$ is a holomorphic tangent vector field along $M$,
which is non-vanishing at any CR point of $M.$  Define
$\chi=(z_1+\ov z_1)|z_2|^2,$ which is real valued.  Then one
computes that $L(\chi)=0.$  As shown in [Hu2], $\{h=0,\ov
h=0,\chi=0\}$ defines a submanifold $X_{p_0}$ in $M$ of real
dimension 3 near $p_0.$  Since $L$ is tangent to $X_{p_0},$
$X_{p_0}$ has to be a CR submanifold with CR dimension one near
$p_0.$ By definition, $M$ is non-minimal at $p_0.$ However $(M,0)$
can not be quadratically flattened.
\medskip

{\it Example 3. 2:}  Let $M\subset\mathbb C^3$ with coordinates
$(z_1,z_2,w)$ be defined by $w=z_1\ov
z_2+\frac{1}{2}z_2^2+\frac{1}{2}\ov z_2^2.$  It is easy to see that
$0\in M$ is a  CR singular point and that the quadratic term in the
defining function takes the following normal form:
$$\mathcal A^\prime=\left(\begin{matrix}0&0\\ 0&\frac{1}{2}\end{matrix}\right),\quad\mathcal B^\prime=\left(\begin{matrix}0&1\\
0&0
\end{matrix}\right).$$
Define $L=(z_2+\ov z_1)\frac{\partial}{\partial z_1}+(z_2\ov z_2+\ov
z_1\ov z_2)\frac{\partial}{\partial w}$ and $\chi=|z_2|^2$.  Write
$h=-w+z_1\ov z_2+\frac{1}{2}z_2^2+\frac{1}{2}\ov z_2^2.$ Then one
computes that $L(h)=L(\ov h)=L(\chi)=0.$ Through the similar
argument as above, we verify that $M$ is CR non-minimal at its CR
points. However $(M,0)$ can not be quadratically flattened.
%As an immediate corollary, we have the following:
%\begin{corollary}\label{05} Let $M$ be a real codimension two submanifold in
%    ${\mathbb C}^{3}$ with $0\in M$ a  CR singular point. Suppose that
%    $M$ near $p=0$ is  defined by an equation as in (\ref{001}) and
%    assume that the CR points in $M$ are non-minimal.
%    If  $\det (B)\not =0$, then $(M,0)$ can be quadratically
%    flattened.
%\end{corollary}
\medskip

{\it Example 3.3:}  Let $M\subset\mathbb C^3$ with coordinates
$(z_1,z_2,w)$ be defined by $w=z_1\ov z_2.$  It is easy to see that
$0\in M$ is a  CR singular point and that the quadratic term in the
defining function takes the following normal form:
$$\mathcal A^\prime=\left(\begin{matrix}0&0\\ 0&0\end{matrix}\right),\quad\mathcal B^\prime=\left(\begin{matrix}0&1\\
0&0
\end{matrix}\right).$$
Define $L=\ov z_1\frac{\partial}{\partial z_1}+\ov z_1\ov
z_2\frac{\partial}{\partial w}$ and $\chi=|z_2|^2$.  Write
$h=-w+z_1\ov z_2.$ Then one computes that $L(h)=L(\ov h)=L(\chi)=0.$
Through a similar argument as above, one verifies that $M$ is CR
non-minimal at its CR points. However $(M,0)$ can not be
quadratically flattened.

\section{Proof of Theorem \ref{04}}
In this section, we give a  proof of Theorem $\ref{04}$. Assume that
$M$ is defined by (\ref{429eq1}) with $p=0\in M$ a CR singular
point. We adapt the notations we have set up so far. We assume that
$\{{\mathcal A}, {\mathcal B}\}$ already takes the normal form
listed in the above section. Our argument is through a computation
based on the fundamental identity (\ref{core})  obtained in $\S 2$.
Since ${\mathcal B}\not =0$, it is apparent that any CR point in $M$
has CR dimension one. Hence, the fundamental identity (\ref{core})
or (\ref{8}) can be applied.
%based on a case by case argument in terms of the Coiffman normal
%form for the pair $\{{\mathcal A},{\mathcal B}\}$. We assume that
%$\{{\mathcal A},{\mathcal B}\}$ is already in the normal form stated
%in $\S 3$. We assume also that a generic point in $M$ is a CR
%non-minimal point.
%We will show that most of the
%cases in the classification  of Coffman do not appear, as the normal
%form for the quadratic terms in the defining function in the setting
%of our theorem, except Case $(4f)$ and some special cases of Case
%$(4c)$ and Case $(4d)$.

\subsection{Case $(1a)$} We start by assuming that the pair $\{\mathcal A,\mathcal
B\}$ takes the normal form in $(1a)$. A  direct computation  gives the
following:
\begin{equation*}\begin{array}{lll}
&w=(az_1^2+ a\- {z_1}^2+2bz_1z_2+2\-b\-{z_1}\-{z_2}+dz_2^2+d\-{z_2}^2)+z_1\-{z_1}+{\rm{cos}}\theta|z_2|^2+i({\rm{sin}}\theta|z_2|^2)+o(|z|^2),\\
&G=(az_1^2+ a\- {z_1}^2+2bz_1z_2+2\-b\-{z_1}\-{z_2}+dz_2^2+d\-{z_2}^2)+z_1\-{z_1}+{\rm{cos}}\theta|z_2|^2+o(|z|^2),\\
&E={\rm{sin}}\theta|z_2|^2+o(|z|^2),\,\,\,G_1=2az_1+2bz_2+\ov z_1+o(|z|),\,\,\,G_2={\rm {cos}}\theta \ov z_2+2bz_1+2dz_2+o(|z|),\\
&E_1=o(|z|)\,,\,\,\,E_2={\rm {sin}}\theta \ov z_2+o(|z|).
\end{array}\end{equation*}
We further compute (\ref{325eq2}), (\ref{325eq3}) and (\ref{325eq4})
to get the following:
\begin{equation*}\begin{array}{lll}
&A=e^{-i\theta}\ov z_2+2bz_1+2dz_2+o(|z|),\,\,\,\,B=2az_1+2bz_2+\ov z_1+o(|z|),\\
%&C=-2i(2az_1+2bz_2+\ov z_1){\rm {sin}}\theta\ov z_2 .\\
%&C=-2i(2az_1+2bz_2+\ov z_1){\rm {sin}}\theta\ov z_2 .\\
&\lambda_{(1)}=-(2az_1+2bz_2+\ov z_1)e^{i\theta}+o(|z|),\,\,\,\,\lambda_{(2)}=-(e^{-i\theta} \ov z_2+2bz_1+2dz_2)+o(|z|),\\
&\lambda_{(4)}=(2 a\ov z_1+2\- b\ov z_2+z_1)e^{-i\theta}+o(|z|),\,\,\,\,\lambda_{(5)}=(e^{i\theta}z_2+2\ov b\ov z_1+2d\ov z_2)+o(|z|),\\
&\Gamma_{(1)}=-2a\ov z_2+(4b^2-4ad)e^{i\theta}z_2+2b\ov z_1e^{i\theta}+o(|z|),\\
&\Gamma_{(2)}=-2be^{-i\theta}\ov z_2+(4ad-4b^2)z_1+2d\ov z_1+o(|z|),\\
&\Gamma_{(4)}=2be^{-i\theta}z_1-(4abe^{-i\theta}+4\-bd)\-z_1+(4de^{-i\theta}-2de^{i\theta})z_2+(2e^{-2i\theta}-4|b|^2e^{-i\theta}-4d^2)\-z_2+o(|z|),\\
&\Gamma_{(5)}=(2ae^{-i\theta}-4ae^{i\theta})z_1+(4a^2e^{-i\theta}+4|b|^2-2e^{i\theta})\ov z_1-2be^{i\theta}z_2+(4a\ov be^{-i\theta}+4bd)\ov z_2+o(|z|).
\end{array}\end{equation*}
Substituting the above quantities into (\ref{7}), we get
$X_1,X_2,Y_1$ and $Y_2$:
\begin{equation*}\begin{array}{lll}
X_1=&(2be^{i\theta}+2\-b(4ad-4b^2))z_1\-z_1+(-2a+2d(4ad-4b^2))z_1\-z_2
+(4abe^{i\theta}+4\-bd)\-z_1\-z_1\\
&+(2a(4b^2-4ad)e^{i\theta}+2de^{i\theta})\-z_1z_2+(-4a^2+4d^2+4|b|^2e^{i\theta}-4|b|^2e^{-i\theta})\-z_1\-z_2\\
&+(2\-b(4b^2-4ad)e^{i\theta}-2b)z_2\-z_2+(-4a\-b-4bde^{-i\theta})\-z_2\-z_2+o(|z|^2),\\
X_2=&(2ae^{-i\theta})z_1z_1+(4|b|^2+4a^2e^{-i\theta}+e^{-i\theta})z_1\-z_1+(2be^{i\theta}+2b^{-i\theta})z_1z_2
+(2ae^{-i\theta})\-z_1\-z_1\\
&+(4bd+4a\-be^{-i\theta})z_1\-z_2
+(4\-bd+4abe^{-i\theta})\-z_1z_2+(4\-be^{-i\theta})\-z_1\-z_2
+(2de^{i\theta})z_2z_2\\
&+(4d^2+1+4|b|^2e^{-i\theta})z_2\-z_2+(2de^{-i\theta})\-z_2\-z_2+o(|z|^2),\\
Y_1=&(8abe^{-i\theta}-8abe^{i\theta})z_1z_1+(8b|b|^2-8a\-bd+2be^{-i\theta}-4be^{i\theta})z_1\-z_1+(-4abe^{-i\theta}-4\-bd)\-z_1\-z_1\\
&+(4b^2e^{-i\theta}-4b^2e^{i\theta}+12ade^{-i\theta}-12ade^{i\theta})z_1z_2
+(6ae^{-2i\theta}-4a+8b^2d-8ad^2)z_1\-z_2
\\
&+(8a^2de^{-i\theta}+4de^{-i\theta}-6de^{i\theta}-8ab^2e^{-i\theta})\-z_1z_2
+(4a^2e^{-2i\theta}-4d^2+2e^{-2i\theta}-2)\-z_1\-z_2\\
&+(8bde^{-i\theta}-8bde^{i\theta})z_2z_2
+(8a\-bde^{-i\theta}+4be^{-2i\theta}-8b|b|^2e^{-i\theta}-2b)z_2\-z_2\\
&+(4a\-be^{-2i\theta}+4bde^{-i\theta})\-z_2\-z_2+o(|z|^2),  \\
Y_2=&(-2ae^{i\theta})z_1z_1+(-4a^2e^{i\theta}-4|b|^2-e^{i\theta})z_1\-z_1+(-4be^{i\theta})z_1z_2+(-4a\-be^{i\theta}-4bd)z_1\-z_2\\
&+(-2ae^{i\theta})\-z_1\-z_1+(-4abe^{i\theta}-4\-bd)\-z_1z_2+(-2\-be^{-i\theta}-2\-be^{i\theta})\-z_1\-z_2+(-2de^{i\theta})z_2z_2\\
&+(-4d^2-1-4|b|^2e^{i\theta})z_2\-z_2+(-2de^{-i\theta})\-z_2\-z_2+o(|z|^2).
\end{array}\end{equation*}
Now we substitute the above formulas into (\ref{8}):
$X_1X_2=Y_1Y_2\,.$ Comparing  the coefficients for the  $\-z_1^4$
terms, we have:
$$(4abe^{i\theta}+4\-bd)(2ae^{-i\theta})=(-4abe^{-i\theta}-4\-bd)(-2ae^{i\theta}).$$
Hence, we see that $8a\-bd(e^{i\theta}-e^{-i\theta})=0.$

Because $a>0,d>0,0<\theta<\pi$, we must have $b=0.$  Substituting $b=0$
back into $X_1,X_2,Y_1$ and $Y_2,$ we obtain:
\begin{equation}\begin{array}{lll}\label{11}
X_1=&(4d^2-4a^2)\-z_1\-z_2+(2de^{i\theta}-8a^2de^{i\theta})\-z_1z_2+(-2a+8ad^2)z_1\-z_2+o(|z|^2),   \\
X_2=&(2ae^{-i\theta})z_1z_1+(4a^2e^{-i\theta}+e^{-i\theta})z_1\-z_1+(2ae^{-i\theta})\-z_1\-z_1+(2de^{i\theta})z_2z_2+(4d^2+1)z_2\-z_2\\
&+(2de^{-i\theta})\-z_2\-z_2+o(|z|^2),\\
Y_1=&(12ade^{-i\theta}-12ade^{i\theta})z_1z_2+(6ae^{-2i\theta}-4a-8ad^2)z_1\-z_2+(8a^2de^{-i\theta}+4de^{-i\theta}-6de^{i\theta})\-z_1z_2\\
&+(4a^2e^{-2i\theta}-2+2e^{-2i\theta}-4d^2)\-z_1\-z_2+o(|z|^2), \\
Y_2=&(-2ae^{i\theta})z_1z_1+(-4a^2e^{i\theta}-e^{i\theta})z_1\-z_1+(-2ae^{i\theta})\-z_1\-z_1+(-2de^{i\theta})z_2z_2+(-1-4d^2)z_2\-z_2\\
&+(-2de^{-i\theta})\-z_2\-z_2+o(|z|^2).
\end{array}\end{equation}

Now compare the coefficients for the $z_2^3\-z_1$ terms in
($\ref{8}$) after being substituted  by  ($\ref{11}$). We then have:
$$(2de^{i\theta}-8a^2de^{i\theta})(2de^{i\theta})=(8a^2de^{-i\theta}+4de^{-i\theta}-6de^{i\theta})(-2de^{i\theta}).$$
Thus, we see that $(16a^2d^2+8d^2)(1-e^{2i\theta})=0.$  However this can not
hold since $a>0,d>0$ and $0<\theta<\pi$.  Thus, we proved that,
under the assumptions in Theorem \ref{04}, the pair $\{\mathcal
A,\mathcal B\}$ can not take the normal form in (1a).

\subsection{Case $(1b)$}
We now assume that the pair $\{\mathcal A,\mathcal B\}$ takes the
normal form in $(1b)$.  We modify the computation we just did for
Case $(1a)$ by substituting $a=0$ into $X_1,X_2,Y_1$ and $Y_2$ and
by requiring that $b\geq 0,d\geq 0$.  We then  have:
\begin{equation*}\begin{array}{lll}
X_1=&(2be^{i\theta}-8b^3)z_1\-z_1+(-8b^2d)z_1\-z_2+(4bd)\-z_1\-z_1+(2de^{i\theta})\-z_1z_2+(8b^3e^{i\theta}-2b)z_2\-z_2\\
&+(4b^2e^{i\theta}-4b^2e^{-i\theta}+4d^2)\-z_1\-z_2+(-4bde^{-i\theta})\-z_2\-z_2+o(|z|^2),\\
X_2=&(4b^2+e^{-i\theta})z_1\-z_1+(2be^{i\theta}+2be^{-i\theta})z_1z_2+(4bd)z_1\-z_2+(4bd)\-z_1z_2+(4be^{-i\theta})\-z_1\-z_2\\
&+(2de^{i\theta})z_2z_2+(4d^2+1+4b^2e^{-i\theta})z_2\-z_2+(2de^{-i\theta})\-z_2\-z_2+o(|z|^2),  \\
Y_1=&(8b^3+2be^{-i\theta}-4be^{i\theta})z_1\-z_1+(4b^2e^{-i\theta}-4b^2e^{i\theta})z_1z_2+(8b^2d)z_1\-z_2+(-4bd)\-z_1\-z_1\\
&+(4de^{-i\theta}-6de^{i\theta})\-z_1z_2
+(2e^{-2i\theta}-4d^2-2)\-z_1\-z_2+(8bde^{-i\theta}-8bde^{i\theta})z_2z_2\\
&+(4be^{-2i\theta}-2b-8b^3e^{-i\theta})z_2\-z_2+(4bde^{-i\theta})\-z_2\-z_2+o(|z|^2),   \\
Y_2=&(-e^{-i\theta}-4b^2)z_1\-z_1+(-4be^{i\theta})z_1z_2+(-4bd)\-z_1z_2+(-2be^{i\theta}-2be^{-i\theta})\-z_1\-z_2\\
&+(-4bd)z_1\-z_2+(-2de^{i\theta})z_2z_2+(-4b^2e^{i\theta}-1-4d^2)z_2\-z_2+(-2de^{-i\theta})\-z_2\-z_2+o(|z|^2).
\end{array}\end{equation*}

Substituting $X_1,X_2,Y_1$ and $Y_2$ into  (\ref{8}) and comparing
the coefficients for the  $z_2^4$ terms, we have:
$$(8bde^{-i\theta}-8bde^{i\theta})(-2de^{i\theta})=0.$$
Since $0<\theta<\pi,$ we either have $b=0$ or $d=0$.  If $b=0$, we
have:
\begin{equation*}\begin{array}{lll}
X_1&=(2de^{i\theta})\-z_1z_2+4d^2\-z_1\-z_2+o(|z|^2),\\
X_2&=e^{-i\theta}z_1\-z_1+(2de^{i\theta})z_2z_2+(4d^2+1)z_2\-z_2+(2de^{-i\theta})\-z_2\-z_2+o(|z|^2),   \\
Y_1&=(4de^{-i\theta}-6de^{i\theta})\-z_1z_2+(2e^{-2i\theta}-4d^2-2)\-z_1\-z_2+o(|z|^2),\\
Y_2&=-e^{-i\theta}z_1\-z_1+(-2de^{i\theta})z_2z_2+(-1-4d^2)z_2\-z_2+(-2de^{-i\theta})\-z_2\-z_2+o(|z|^2).
\end{array}\end{equation*}

Comparing the coefficients for the $z_1\-z_1^2\-z_2$ terms in (\ref{8}), we obtain:
$$e^{-i\theta}4d^2=-e^{-i\theta}(2e^{-2i\theta}-4d^2-2).$$
This contradicts with the fact $0<\theta<\pi$.  If $d=0$,
\begin{equation*}\begin{array}{lll}
X_1=&(2be^{i\theta}-8b^3)z_1\-z_1+(4b^2(e^{i\theta}-e^{-i\theta}))\-z_1\-z_2+(8b^3e^{i\theta}-2b)z_2\-z_2+o(|z|^2),\\
X_2=&(2be^{i\theta}+2be^{-i\theta})z_1z_2+(4b^2+e^{-i\theta})z_1\-z_1+(4be^{-i\theta})\-z_1\-z_2+(1+4b^2e^{-i\theta})z_2\-z_2+o(|z|^2),\\
Y_1=&(4b^2e^{-i\theta}-4b^2e^{i\theta})z_1z_2+(8b^3-4be^{i\theta}+2be^{-i\theta})z_1\-z_1+(2e^{-2i\theta}-2)\-z_1\-z_2\\
&+(-2b+4be^{-2i\theta}-8b^3e^{-i\theta})z_2\-z_2+o(|z|^2),\\
Y_2=&(-e^{i\theta}-4b^2)z_1\-z_1+(-4be^{i\theta})z_1z_2+(-4b^2e^{i\theta}-1)z_2\-z_2+(-2be^{i\theta}-2be^{-i\theta})\-z_1\-z_2+o(|z|^2).
\end{array}\end{equation*}
Comparing the coefficients for the $z_1^2\-z_1^2$ terms in (\ref{8}), we have:
$$(2be^{i\theta}-8b^3)(4b^2+e^{-i\theta})=(8b^3-4be^{i\theta}+2be^{-i\theta})(-e^{i\theta}-4b^2).$$
Therefore $4b(1-e^{2i\theta})=0.$ Since $0<\theta<\pi$, we conclude
$b=0$ and thus we are in the setting of   the previous case, which
is impossible as just shown. Thus, we proved that, under the
assumptions in Theorem \ref{04}, the pair $\{\mathcal A,\mathcal
B\}$ can not take the normal form in $(1b)$.

\subsection{Case $(1c)$}
Assume that the pair $\{\mathcal A,\mathcal B\}$ takes the normal
form in $(1c)$.  We modify our argument  in $(1a)$ by  substituting
$d=0$ into $X_1,X_2,Y_1, Y_2$ and by requiring that $b\geq 0,a> 0$.
We have:
\begin{equation*}\begin{array}{lll}
X_1=&(2be^{i\theta}-8b^3)z_1\-z_1+(-2a)z_1\-z_2+(4abe^{i\theta})\-z_1\-z_1+(-4a^2+4b^2e^{i\theta}-4b^2e^{-i\theta})\-z_1\-z_2\\
&+(8ab^2e^{i\theta})\-z_1z_2+(8b^3e^{i\theta}-2b)z_2\-z_2+(-4ab)\-z_2\-z_2+o(|z|^2),\\
X_2=&(2ae^{-i\theta})z_1z_1+(2be^{i\theta}+2be^{-i\theta})z_1z_2+(4b^2+4a^2e^{-i\theta}+e^{-i\theta})z_1\-z_1+(4abe^{-i\theta})z_1\-z_2\\
&+(4abe^{-i\theta})z_2\-z_1+(2ae^{-i\theta})\-z_1\-z_1+(4be^{-i\theta})\-z_1\-z_2+(1+4b^2e^{-i\theta})z_2\-z_2+o(|z|^2),\\
Y_1=&(8abe^{-i\theta}-8abe^{i\theta})z_1z_1+(8b^3-4be^{i\theta}+2be^{-i\theta})z_1\-z_1+(4b^2e^{-i\theta}-4b^2e^{i\theta})z_1z_2\\
&+(6ae^{-2i\theta}-4a)z_1\-z_2+(-4abe^{-i\theta})\-z_1\-z_1+(-8ab^2e^{-i\theta})\-z_1z_2+(4a^2e^{-2i\theta}-2+2e^{-2i\theta})\-z_1\-z_2\\
&+(-2b+4be^{-2i\theta}-8b^3e^{-i\theta})z_2\-z_2+(4abe^{-2i\theta})\-z_2\-z_2+o(|z|^2),\\
Y_2=&(-2ae^{i\theta})z_1z_1+(-4abe^{i\theta})z_1\-z_2+(-4be^{i\theta})z_1z_2+(-4a^2e^{i\theta}-e^{i\theta}-4b^2)z_1\-z_1\\
&+(-2ae^{i\theta})\-z_1\-z_1+(-4abe^{i\theta})\-z_1z_2+(-2be^{i\theta}-2be^{-i\theta})\-z_1\-z_2+(-4b^2e^{i\theta}-1)z_2\-z_2+o(|z|^2).
\end{array}\end{equation*}

Compare the coefficients for the $\-z_2^3z_2$ terms in (\ref{8}). We
get
$$(-4ab)(1+4b^2e^{-i\theta})=(4abe^{-2i\theta})(-4b^2e^{i\theta}-1).$$
Therefore, $4ab(e^{2i\theta}-1)=0.$ Since $0<\theta<\pi$, we
conclude that either $a=0$ or $b=0.$  The case for  $a=0$ can be
included in Case $(1b)$, which has been shown to be  impossible.
When $b=0$, we have:
\begin{equation*}\begin{array}{lll}
X_1=&(-2a)z_1\-z_2+(-4a^2)\-z_1\-z_2+o(|z|^2),\\
X_2=&(2ae^{-i\theta})z_1z_1+(4a^2e^{-i\theta}+e^{-i\theta})z_1\-z_1+(2ae^{-i\theta})\-z_1\-z_1+z_2\-z_2+o(|z|^2),   \\
Y_1=&(6ae^{-2i\theta}-4a)z_1\-z_2+(4a^2e^{-2i\theta}-2+2e^{-2i\theta})\-z_1\-z_2+o(|z|^2),\\
Y_2=&(-2ae^{i\theta})z_1z_1+(-4a^2e^{i\theta}-e^{i\theta})z_1\-z_1+(-2ae^{i\theta})\-z_1\-z_1-z_2\-z_2+o(|z|^2).
\end{array}\end{equation*}
Comparing the coefficients for the $z_1^3\-z_2$ in (\ref{8}), we get
$$(-2a)(2ae^{-i\theta})=(6ae^{-2i\theta}-4a)(-2ae^{i\theta}).$$
Hence, $8a^2(1-e^{2i\theta})=0.$ Since $0<\theta<\pi$, $a=0$, which
is reduced to the case discussed above.  Thus, we proved that, under
the assumptions in Theorem \ref{04}, the pair $\{\mathcal A,\mathcal
B\}$ can not take the normal form in $(1c)$.

\subsection{Case $(2a)$}
 Assume that the pair $\{\mathcal A,\mathcal
 B\}$ takes the normal form in $(2a)$.  Notice that in Case $(2a)-(2f)$, $b$ is required to be a real number. A direct computation  gives the  following:
\begin{equation*}\begin{array}{lll}
&w=(az_1^2+ \-a\- {z_1}^2+2bz_1z_2+2b\-{z_1}\-{z_2}+dz_2^2+\-d\-{z_2}^2)+z_1\-{z_2}+\tau\-z_1z_2+o(|z|^2),\\
&G=(az_1^2+ \-a\- {z_1}^2+2bz_1z_2+2b\-{z_1}\-{z_2}+dz_2^2+\-d\-{z_2}^2)+\frac{1+\tau}{2}(z_1\-{z_2}+\-z_1z_2)+o(|z|^2),\\
&E=\frac{1-\tau}{2i}(z_1\-{z_2}-\-z_1z_2)+o(|z|^2),\\
&G_1=2az_1+2bz_2+\frac{1+\tau}{2}\-z_2+o(|z|),\,\,G_2=2bz_1+2dz_2+\frac{1+\tau}{2}\-z_1+o(|z|),\\
&E_1=\frac{1-\tau}{2i}\ov z_2+o(|z|)\,,\,\,E_2=\frac{1-\tau}{2i}(-\ov z_1)+o(|z|).
\end{array}\end{equation*}
We further compute  (\ref{325eq2}), (\ref{325eq3}) and
(\ref{325eq4}) to derive the following:
\begin{equation*}\begin{array}{lll}
&A=2bz_1+2dz_2+\-z_1+o(|z|),\,\,\,B=2az_1+2bz_2+\tau\ov z_2+o(|z|),\\
%&C=-2i(2az_1+2bz_2+\ov z_1){\rm {sin}}\theta\ov z_2 .\\
%&C=-2i(2az_1+2bz_2+\ov z_1){\rm {sin}}\theta\ov z_2 .\\
&\lambda_{(1)}=2bz_1+2dz_2+\-z_1+o(|z|),\,\,\,\lambda_{(2)}=(2az_1+2bz_2+\tau\ov z_2)\tau+o(|z|),\\
&\lambda_{(4)}=-(2b\-z_1+2\-d\-z_2+z_1)+o(|z|),\,\,\,\lambda_{(5)}=-(2\-a\-z_1+2b\-z_2+\tau z_2)\tau+o(|z|),\\
&\Gamma_{(1)}=(4b^2-4ad)z_1+2b\-z_1-2\tau d\-z_2+o(|z|),\\
&\Gamma_{(2)}=(4ad\tau-4b^2\tau) z_2+(2a\tau)\-z_1-2b\tau^2\-z_2+o(|z|),\\
&\Gamma_{(4)}=-2bz_1+(4b^2+4\-ad\tau-2)\-z_1+(2d\tau^2-4d)z_2+(4b\-d+4bd\tau)\-z_2+o(|z|),\\
&\Gamma_{(5)}=(4a\tau^2-2a)z_1+(-4ab-4\-ab\tau)\ov z_1+(2b\tau^2)z_2+(2\tau^3-4b^2\tau-4a\-d)\ov z_2+o(|z|).\\
\end{array}\end{equation*}
Substituting the above quantities into (\ref{7}), we have
\begin{equation*}\begin{array}{lll}
X_1=&(2\-a(4b^2-4ad)+2a\tau)z_1\-z_1+(4\-ab+4ab\tau)\-z_1\-z_1+(2b(4b^2-4ad)-2b\tau^2)z_1\-z_2\\
&+(2b\tau+2b(4ad\tau-4b^2\tau))\-z_1z_2
+(4\tau a\-d-4\tau\-ad+4b^2-4b^2\tau^2)\-z_1\-z_2\\
&+(2\-d(4ad\tau-4b^2\tau)-2\tau^2d)z_2\-z_2+(-4\tau bd-4b\-d\tau^2)\-z_2\-z_2+o(|z|^2),    \\
X_2=&(-2a)z_1z_1+(-4\-ab\tau-4ab)z_1\-z_1+(-2b\tau^2-2b)z_1z_2+(-4b^2\tau-4a\-d-\tau)z_1\-z_2\\
&+(-2\-a\tau)\-z_1\-z_1+(-4\-ad\tau-\tau^2-4b^2)\-z_1z_2+(-4b\tau)\-z_1\-z_2+(-2d\tau^2)z_2z_2\\
&+(-4bd\tau-4b\-d)z_2\-z_2+(-2\-d\tau)\-z_2\-z_2+o(|z|^2),    \\
Y_1=&(8ab\tau^2-8ab)z_1z_1+(4a\tau^2-6a+8|a|^2d\tau-8\-ab^2\tau)z_1\-z_1+(4b\-d\tau+4bd\tau^2)\-z_2\-z_2  \\
&+(4b\tau^3-8b^3\tau+8abd\tau-2b\tau)z_1\-z_2+(-4ab-4\-ab\tau)\-z_1\-z_1+(8bd\tau^2-8bd)z_2z_2\\
&+(2b(4b^2+4\-ad\tau-2)+2b\tau^2-2d(4ab+4\-ab\tau))\-z_1z_2+(2\tau^3-2\tau+4\-ad\tau^2-4a\-d)\-z_1\-z_2\\
&+(6d\tau^3-8a|d|^2+8b^2\ov d-4d\tau)z_2\-z_2+(12ad\tau^2-12ad+4b^2\tau^2-4b^2)z_1z_2+o(|z|^2),\\
Y_2=&(2a\tau)z_1z_1+(4\-ab+4ab\tau)z_1\-z_1+(4b\tau)z_1z_2+(4b^2+4a\-d\tau+\tau^2)z_1\-z_2+(2\-a)\-z_1\-z_1\\
&+(4\-ad+\tau+ 4b^2\tau)\-z_1z_2+(2b+2b\tau^2)\-z_1\-z_2+(2d\tau)z_2z_2+(4bd+4b\-d\tau)z_2\-z_2\\
&+(2\-d\tau^2)\-z_2\-z_2+o(|z|^2).
\end{array}\end{equation*}
Now comparing the coefficients for the $z_1^4$ terms in (\ref{8}), we have:
$$(8ab\tau^2-8ab)(2a\tau)=0.$$
Since $a,b\neq 0$ and $0<\tau<1,$ the above equation has no
solution. Thus, we proved that, under the assumptions in Theorem
\ref{04}, the pair $\{\mathcal A,\mathcal B\}$ can not take the
normal form in $(2a)$.

\subsection{Case $(2b)$}
 We now assume that the pair $\{\mathcal A,\mathcal
 B\}$ takes the normal form in $(2b)$.  Modifying the formulas  in  Case $(2a)$ by substituting $a=0$ into $X_1,X_2,Y_1,Y_2$ and by
 requiring that $b,d\neq 0$ and $0<\tau<1$, we have the
following:
\begin{equation*}\begin{array}{lll}
X_1=&(8b^3-2b\tau^2)z_1\-z_2 +(2b\tau-8b^3\tau)\-z_1z_2
+(4b^2-4b^2\tau^2)\-z_1\-z_2+(-8\-db^2\tau-2\tau^2d)z_2\-z_2\\
&+(-4\tau bd-4b\-d\tau^2)\-z_2\-z_2+o(|z|^2),\\
X_2=&(-2b\tau^2-2b)z_1z_2+(-4b^2\tau-\tau)z_1\-z_2
+(-\tau^2-4b^2)\-z_1z_2+(-4b\tau)\-z_1\-z_2\\
&+(-2d\tau^2)z_2z_2+(-4bd\tau-4b\-d)z_2\-z_2+(-2\-d\tau)\-z_2\-z_2+o(|z|^2),  \\
Y_1=&(4b^2\tau^2-4b^2)z_1z_2 +(4b\tau^3-8b^3\tau-2b\tau)z_1\-z_2
+(2\tau^3-2\tau)\-z_1\-z_2
+(8b^3-4b+2b\tau^2)\-z_1z_2\\
&+(8bd\tau^2-8bd)z_2z_2
+(6d\tau^3+8b^2\-d-4d\tau)z_2\-z_2+(4b\-d\tau+4bd\tau^2)\-z_2\-z_2+o(|z|^2),    \\
Y_2=&(4b\tau)z_1z_2+(4b^2+\tau^2)z_1\-z_2
+(\tau+ 4b^2\tau)\-z_1z_2+(2b+2b\tau^2)\-z_1\-z_2+(2d\tau)z_2z_2+(4bd\\
&+4b\-d\tau)z_2\-z_2+(2\-d\tau^2)\-z_2\-z_2+o(|z|^2).
\end{array}\end{equation*}
Now substitute the above formulas for $X_1,X_2,Y_1$ and $Y_2$ into (\ref{8}). By comparing the coefficients for the $z_2^4$ terms, we have:
$$(8bd\tau^2-8bd)2d\tau=0,$$
which is impossible because $b,d\neq 0$ and $0<\tau<1.$  Thus, we proved that,
under the assumptions in Theorem \ref{04}, the pair $\{\mathcal
A,\mathcal B\}$ can not take the normal form in $(2b)$.

\subsection{Case $(2c)$}
 Assume that the pair $\{\mathcal A,\mathcal
 B\}$ takes the normal form in $(2c)$.  Modifying the formulas  in the Case $(2a)$ and $(2b)$ by substituting $a=0,d=0$ into $X_1,X_2,Y_1,Y_2$, we have the
 following:
\begin{equation*}\begin{array}{lll}
X_1=&(8b^3-2b\tau^2)z_1\-z_2 +(2b\tau-8b^3\tau)\-z_1z_2
+(4b^2-4b^2\tau^2)\-z_1\-z_2+o(|z|^2),    \\
X_2=&(-2b\tau^2-2b)z_1z_2+(-4b^2\tau-\tau)z_1\-z_2
+(-\tau^2-4b^2)\-z_1z_2+(-4b\tau)\-z_1\-z_2+o(|z|^2),\\
Y_1=&(4b^2\tau^2-4b^2)z_1z_2 +(4b\tau^3-8b^3\tau-2b\tau)z_1\-z_2
+(2\tau^3-2\tau)\-z_1\-z_2\\
&+(8b^3-4b+2b\tau^2)\-z_1z_2+o(|z|^2),\\
Y_2=&(4b\tau)z_1z_2+(4b^2+\tau^2)z_1\-z_2 +(\tau+
4b^2\tau)\-z_1z_2+(2b+2b\tau^2)\-z_1\-z_2+o(|z|^2).
\end{array}\end{equation*}
Comparing the coefficients for the $z_1^2z_2^2$ terms in the identity
$X_1X_2=Y_1Y_2$, we have:
$$(4b^2\tau^2-4b^2)4b\tau=0.$$
This is not possible because $b\neq 0$ and $0<\tau<1$. Thus, under
the assumptions in Theorem \ref{04}, the pair $\{\mathcal A,\mathcal
B\}$ can not take the normal form in $(2c)$.
\subsection{Case $(2d),(2e),(2f)$}
Use the computation we derived in the Case $(2a)$ and substitute
$b=0$ into $X_1,X_2,Y_1,Y_2$. We have the
following:
\begin{equation}\begin{array}{lll}\label{ww}
X_1=&(2a\tau-8|a|^2d)z_1\-z_1
+(4\tau a\-d-4\tau\-ad)\-z_1\-z_2+(8a|d|^2\tau-2\tau^2d)z_2\-z_2+o(|z|^2),\\
X_2=&(-2a)z_1z_1+(-4a\-d-\tau)z_1\-z_2+(-2\-a\tau)\-z_1\-z_1
+(-4\-ad\tau-\tau^2)\-z_1z_2+(-2d\tau^2)z_2z_2\\
&+(-2\-d\tau)\-z_2\-z_2+o(|z|^2),   \\
Y_1=&(4a\tau^2-6a+8|a|^2d\tau)z_1\-z_1+(12ad\tau^2-12ad)z_1z_2
+(2\tau^3-2\tau+4\-ad\tau^2-4a\-d)\-z_1\-z_2\\
&+(6d\tau^3-8a|d|^2-4d\tau)z_2\-z_2+o(|z|^2),\\
Y_2=&(2a\tau)z_1z_1+(4a\-d\tau+\tau^2)z_1\-z_2+(2\-a)\-z_1\-z_1
+(4\-ad+\tau
)\-z_1z_2+(2d\tau)z_2z_2+(2\-d\tau^2)\-z_2\-z_2+o(|z|^2).

\end{array}\end{equation}
Comparing the coefficients for the $\-z_1^3z_1$ terms in the
identity (\ref{8}), we have
$$(2a\tau-8|a|^2d)(-2\-a\tau)=(4a\tau^2-6a+8|a|^2d\tau)(2\-a).$$
Hence $12|a|^2(\tau^2-1)=0,$ from which it follows that $a=0$ by
taking into consideration  the fact $0<\tau<1$.

Substituting $a=0$ into the above expressions (\ref{ww}), we have:
\begin{equation*}\begin{array}{lll}\label{2d}
X_1=&-2\tau^2dz_2\-z_2+o(|z|^2),\\
X_2=&-\tau z_1\-z_2-\tau^2\-z_1z_2+(-2d\tau^2)z_2z_2+(-2\-d\tau)\-z_2\-z_2+o(|z|^2),\\
Y_1=&(2\tau^3-2\tau)\-z_1\-z_2+(6d\tau^3-4d\tau)z_2\-z_2+o(|z|^2),\\
Y_2=&\tau^2z_1\-z_2+\tau \-z_1z_2+(2d\tau)z_2z_2+(2\-d\tau^2)\-z_2\-z_2+o(|z|^2).
\end{array}\end{equation*}
Compare the coefficients for the $z_1\-z_1\-z_2^2$ terms. Then
$(2\tau^3-2\tau)\tau^2=0$, which is impossible.  Thus, we proved
that, under the assumptions in Theorem \ref{04}, the pair
$\{\mathcal A,\mathcal B\}$ can not take the normal form in
$(2d),(2e),$ and $(2f)$.

\subsection{Case (3)}
We study the  cases in  $(3a)$, $(3b)$ and $(3c)$  in this
subsection. Notice that in (3a),  (3b) and (3c), $b$ is required to
be a real number. By a  computation,  we derive that
\begin{equation*}\begin{array}{lll}
&w=(az_1^2+ \-a\- {z_1}^2+2bz_1z_2+2\-b\-{z_1}\-{z_2}+dz_2^2+\-d\-{z_2}^2)+z_1\-{z_2}+\-z_1z_2+iz_2\-z_2+o(|z|^2),\\
&G=(az_1^2+ \-a\- {z_1}^2+2bz_1z_2+2\-b\-{z_1}\-{z_2}+dz_2^2+\-d\-{z_2}^2)+z_1\-{z_2}+\-z_1z_2+o(|z|^2),\\
&E=z_2\-z_2+o(|z|^2),\,\,\,G_1=2az_1+2bz_2+\-z_2+o(|z|),\,\,\,G_2=2bz_1+2dz_2+\-z_1+o(|z|),\\
&E_1=o(|z|),\,\,\,E_2=\-z_2+o(|z|).
\end{array}\end{equation*}
We substitute the above quantities into (\ref{325eq2}),\
(\ref{325eq3}) and (\ref{325eq4}), to derive the following:
\begin{equation*}\begin{array}{lll}
&A=2bz_1+2dz_2+\-z_1-i\-z_2+o(|z|),\,\,\,B=2az_1+2bz_2+\-z_2+o(|z|),\\
%&C=-2i(2az_1+2bz_2+\ov z_1){\rm {sin}}\theta\ov z_2 .\\
%&C=-2i(2az_1+2bz_2+\ov z_1){\rm {sin}}\theta\ov z_2 .\\
&\lambda_{(1)}=(2b-2ai)z_1+(2d-2bi)z_2+\-z_1-2i\-z_2+o(|z|)\\
&\lambda_{(2)}=2az_1+2bz_2+\ov z_2+o(|z|),\\
&\lambda_{(4)}=-(2b+2\-ai)\-z_1-(2\-d+2bi)\-z_2-z_1-2iz_2+o(|z|)\\
&\lambda_{(5)}=-(2\-a\-z_1+2b\-z_2+z_2)+o(|z|),\\
&\Gamma_{(1)}=(4b^2-4ad)z_1+(2b-2ai)\-z_1+i(4b^2-4ad)z_2+(-2a-2d)\-z_2+o(|z|),\\
&\Gamma_{(2)}=(4ad-4b^2) z_2+2a\-z_1-(2ai+2b)\-z_2+o(|z|),\\
&\Gamma_{(4)}=(8ai-2b)z_1+(4b^2+4\-abi+4\-ad-2)\-z_1+(12bi-2d)z_2+(4b\-d+4bd+4b^2i+6i)\-z_2+o(|z|),\\
&\Gamma_{(5)}=2az_1+(-4ab-4\-ab-4|a|^2i)\ov z_1+(2b-4ai)z_2+(2-4b^2-4a\-d-4abi)\ov z_2+o(|z|).
\end{array}\end{equation*}
Further we substitute the above quantities into (\ref{7}) to obtain:
\begin{equation*}\begin{array}{lll}
X_1=&(8\-ab^2-8|a|^2d+2a)z_1\-z_1+(8b^3-8abd-2ai-2b)z_1\-z_2+(4\-ab+4ab-4|a|^2i)\-z_1\-z_1\\
&+(8\-ab^2i-8|a|^2di+2b+8abd-8b^3)\-z_1z_2+(4a\-d-4\-ad-4|a|^2-8bai)\-z_1\-z_2\\
&+(8ib^3-8iabd+8a|d|^2-8\-db^2-2d-2bi)z_2\-z_2+(-4ba-4bd-4\-dai-4\-db)\-z_2\-z_2+o(|z|^2),  \\
X_2=&(-2a)z_1z_1+(-4\-ab-4ab-4i|a|^2)z_1\-z_1+(-4b-4ai)z_1z_2\\
&+(-2\-a)\-z_1\-z_1
+(-4\-ad-1-4b^2-4\-abi)\-z_1z_2+(-4b)\-z_1\-z_2+(-2d-4bi)z_2z_2\\
&+(-4b^2-4a\-d-4iab-1)z_1\-z_2+(-4bd-4b\-d-4b^2i-i)z_2\-z_2+(-2\-d)\-z_2\-z_2+o(|z|^2), \\
Y_1=&(16a^2i)z_1z_1+(-8\-ab^2-2a+8|a|^2d)z_1\-z_1+(4b\-d+4bd+4i+8b^2i+4a\-di-4ab)\-z_2\-z_2\\
&+(2b-8b^3+8abd+18ai)z_1\-z_2
+(4\-ad-4a\-d-4|a|^2+8\-abi)\-z_1\-z_2+(32abi)z_1z_2\\
&+(-4ab-4\-ab-4|a|^2i)\-z_1\-z_1+(8b^3-8abd-2b-8d|a|^2i+8\-ab^2i-4ai)\-z_1z_2\\
&+(24b^2i-8adi)z_2z_2+(2d-8a|d|^2+8b^2\-d-4a+22bi-8abdi+8b^3i)z_2\-z_2+o(|z|^2),   \\
Y_2=&(2a)z_1z_1+(4\-ab+4ab-4|a|^2i)z_1\-z_1+(4b)z_1z_2+(4b^2+4a\-d+1-4abi)z_1\-z_2\\
&+(4\-ad+1+ 4b^2-4\-abi)\-z_1z_2+(4b-4\-ai)\-z_1\-z_2+(2d)z_2z_2
+(2\-a)\-z_1\-z_1\\
&+(4bd+4b\-d-i-4b^2i)z_2\-z_2 +(2\-d-4bi)\-z_2\-z_2+o(|z|^2).
\end{array}\end{equation*}
Substituting the above formulas into (\ref{8}) and comparing the
coefficients for the $z_1^4$ terms, we have: $32a^3i=0.$ Hence $a=0$
and thus, the pair $\{\mathcal A,\mathcal B\}$ can not take the
normal form in  $(3a)$.

Substituting $a=0$ into
$X_1,X_2,Y_1,Y_2$, we have:
\begin{equation}\label{a=0}\begin{array}{lll}
        X_1=&(8b^3-2b)z_1\-z_2+(8ib^3-8\-db^2-2d-2bi)z_2\-z_2+(-4bd-4\-db)\-z_2\-z_2\\
        &+(2b-8b^3)\-z_1z_2+o(|z|^2),  \\
        X_2=&(-4b)z_1z_2+(-1-4b^2)\-z_1z_2+(-4b)\-z_1\-z_2+(-2d-4bi)z_2z_2\\
        &+(-4b^2-1)z_1\-z_2+(-4bd-4b\-d-4b^2i-i)z_2\-z_2+(-2\-d)\-z_2\-z_2+o(|z|^2), \\
        Y_1=&(4b\-d+4bd+4i+8b^2i)\-z_2\-z_2+(2b-8b^3)z_1\-z_2+(8b^3-2b)\-z_1z_2\\
        &+(24b^2i)z_2z_2+(2d+8b^2\-d+22bi+8b^3i)z_2\-z_2+o(|z|^2),\\
        Y_2=&(4b)z_1z_2+(4b^2+1)z_1\-z_2+(1+ 4b^2)\-z_1z_2+(4b)\-z_1\-z_2+(2d)z_2z_2\\
        &+(4bd+4b\-d-i-4b^2i)z_2\-z_2+(2\-d-4bi)\-z_2\-z_2+o(|z|^2).
    \end{array}\end{equation}
Substituting the above formulas into (\ref{8}) and comparing the coefficients for the $z_2^3\-z_2$ terms, we have:
\begin{equation*}
(8ib^3-8\-db^2-2d-2bi)(-2d-4bi)=(24b^2i)(4bd+4b\-d-i-4b^2i)+(2d+8b^2\-d+22bi+8b^3i)2d
\end{equation*}
Noticing that, in $(3b)$ and $(3c)$, $b$ and $d$ are both real numbers, we can simplify the above equality to derive the following:
\begin{equation*}
(64b^4+32b^2)+i(192b^3d+32bd)=0.
\end{equation*}
Hence $b=0.$   Substituting $b=0$ further into $(\ref{a=0})$, we have
\begin{equation*}\begin{array}{lll}
X_1=&-2dz_2\-z_2+o(|z|^2),\\
X_2=&-z_1\-z_2-\-z_1z_2+(-2d)z_2z_2+(-i)z_2\-z_2+(-2\-d)\-z_2\-z_2+o(|z|^2),    \\
Y_1=&(2d)z_2\-z_2+(4i)\-z_2\-z_2+o(|z|^2),  \\
Y_2=&z_1\-z_2+\-z_1z_2+(2d)z_2z_2+(-i)z_2\-z_2
+(2\-d)\-z_2\-z_2+o(|z|^2).
\end{array}\end{equation*}
Substituting the above formulas into (\ref{8}) and comparing the coefficients for the $z_1\-z_2^3$, we get $4i=0,$ which is a contradiction.

%Therefore the pair
%$\{\mathcal A, \mathcal B\}$ in $(3c)$ can not be taken  as  a
%normal form under our assumption in Theorem \ref{04}.

%In $(3b)$, $a=d=0$ and $b>0$. Substituting $a=d=0$ into $X_1,X_2,Y_1,Y_2$ and
%(\ref{8}), we have:
%\begin{equation*}\begin{array}{lll}
%X_1=&(8b^3-2b)z_1\-z_2+(2b-8b^3)\-z_1z_2+(8ib^3-2bi)z_2\-z_2+o(|z|^2),\\
%X_2=&(-4b)z_1z_2+(-4b^2-1)z_1\-z_2+(-1-4b^2)\-z_1z_2+(-4b)\-z_1\-z_2+(-4bi)z_2z_2\\
%&+(-4b^2i-i)z_2\-z_2+o(|z|^2), \\
%Y_1=&(2b-8b^3)z_1\-z_2+(8b^3-2b)\-z_1z_2+(24b^2i)z_2z_2
%+(22bi+8b^3i)z_2\-z_2\\
%&+(4i+8b^2i)\-z_2\-z_2+o(|z|^2),    \\
%Y_2=&(4b)z_1z_2+(4b^2+1)z_1\-z_2+(1+4b^2)\-z_1z_2+(4b)\-z_1\-z_2\\
%&+(-i-4b^2i)z_2\-z_2 +(-4bi)\-z_2\-z_2+o(|z|^2).
%\end{array}\end{equation*}
%Comparing the coefficients for the $z_2\ov z_2^3$ terms in (\ref{8}), we have $(4+8b^2)(1+4b^2)+4b(22b+8b^3)=0$, which is a contradiction; for $b>0$.
Therefore, we proved that, under the
assumptions in Theorem \ref{04}, the pair $\{\mathcal A,\mathcal
B\}$ can not take the normal form in $(3)$.

\subsection{Case $(4)$}
We modify the computation we did for $(2a)$ by substituting $\tau=0$
into $X_1,X_2,Y_1,Y_2$. Further noticing that $b$ and $d$ are real
numbers, we obtain:
\begin{equation*}\begin{array}{lll}
X_1=&(8\-ab^2-8|a|^2d))z_1\-z_1+(4\-ab)\-z_1\-z_1+(8b^3-8abd))z_1\-z_2+(4b^2)\-z_1\-z_2+o(|z|^2),\\
X_2=&(-2a)z_1z_1+(-4ab)z_1\-z_1+(-2b)z_1z_2+(-4ad)z_1\-z_2
+(-4b^2)\-z_1z_2+(-4bd)z_2\-z_2+o(|z|^2),   \\
Y_1=&(-8ab)z_1z_1+(-6a)z_1\-z_1+(-12ad-4b^2)z_1z_2 +(-4ab)\-z_1\-z_1
+(-4ad)\-z_1\-z_2\\
&+(8b^3-4b-8abd))\-z_1z_2+(-8bd)z_2z_2
+(-8ad^2+8b^2d)z_2\-z_2+o(|z|^2),   \\
Y_2=&(4\-ab)z_1\-z_1+(4b^2)z_1\-z_2+(2\-a)\-z_1\-z_1
+(4\-ad)\-z_1z_2+(2b)\-z_1\-z_2+(4bd)z_2\-z_2+o(|z|^2).
\end{array}\end{equation*}
Comparing the coefficients for the $\-z_1z_2^3$ terms in
$X_1X_2=Y_1Y_2$, we have: $(-8bd)(4\-ad)=0.$ If $a=0$, then we
substitute it into $X_1,X_2,Y_1,Y_2$ to derive:
\begin{equation*}\begin{array}{lll}
X_1=&(8b^3)z_1\-z_2+(4b^2)\-z_1\-z_2+o(|z|^2),\\
X_2=&(-2b)z_1z_2
+(-4b^2)\-z_1z_2+(-4bd)z_2\-z_2+o(|z|^2),   \\
Y_1=&(-4b^2)z_1z_2+(8b^3-4b)\-z_1z_2+(-8bd)z_2z_2
+(8b^2d)z_2\-z_2+o(|z|^2),  \\
Y_2=&(4b^2)z_1\-z_2+(2b)\-z_1\-z_2+(4bd)z_2\-z_2+o(|z|^2).
\end{array}\end{equation*}
Comparing the coefficients for the $\-z_1^2z_2\-z_2$ terms in
$X_1X_2=Y_1Y_2$, we derive that
$$(-4b^2)(4b^2)=(8b^3-4b)(2b).$$  Hence $b=0$ or $b=\frac{1}{2}.$ If
$a=0,b=0$, then $d=\frac{1}{2}$ or $d=0$ by making use of the normal form.  If $a=0,b=\frac{1}{2},$ by
comparing the coefficients for the $z_2^2\ov z_2^2$ terms, we have $(8b^2d)(4bd)=0.$ Thus, $d=0.$ Therefore we have
the following possibilities:
\begin{enumerate}
    \item $a=0,b=0,d=0\,\,;$
    \item $a=0,b=0,d=\frac{1}{2}\,\,;$
    \item $a=0,b=\frac{1}{2},d=0.$
\end{enumerate}
If $a\neq 0,b=0$, then we have:
\begin{equation*}\begin{array}{lll}
X_1=&(-8|a|^2d)z_1\-z_1+o(|z|^2),\\
X_2=&(-2a)z_1z_1+(-4ad)z_1\-z_2+o(|z|^2),\\
Y_1=&(-6a)z_1\-z_1+(-12ad)z_1z_2 +(-4ad)\-z_1\-z_2
+(-8ad^2)z_2\-z_2+o(|z|^2),\\
Y_2=&(2\-a)\-z_1\-z_1 +(4\-ad)\-z_1z_2+o(|z|^2).
\end{array}\end{equation*}
Comparing the coefficients for the $\-z_1^3z_1$ terms in
$X_1X_2=Y_1Y_2$, we have  $(-6a)(2\-a)=0,$ which is impossible.

 If $a\neq 0,b\neq 0, d=0$, then by the normal form in $\S 3$, we have
$a=\frac{1}{2},b>0.$ Therefore,
\begin{equation*}\begin{array}{lll}
X_1=&(4b^2)z_1\-z_1+(2b)\-z_1\-z_1+(8b^3)z_1\-z_2+(4b^2)\-z_1\-z_2+o(|z|^2),\\
X_2=&-z_1z_1+(-2b)z_1\-z_1+(-2b)z_1z_2+(-4b^2)\-z_1z_2+o(|z|^2),\\
Y_1=&(-4b)z_1z_1+(-3)z_1\-z_1+(-4b^2)z_1z_2 +(-2b)\-z_1\-z_1
+(8b^3-4b)\-z_1z_2+o(|z|^2),    \\
Y_2=&(2b)z_1\-z_1+(4b^2)z_1\-z_2+\-z_1\-z_1
+(2b)\-z_1\-z_2+o(|z|^2).
\end{array}\end{equation*}
Comparing the coefficients for the $\-z_1z_1^3$ terms in
$X_1X_2=Y_1Y_2$, we have  $(4b^2)(-1)=(-4b)(2b)$, which contradicts
 the fact $b>0$.

\medskip
Combing all the above, we finally achieved a proof for  Theorem
$\ref{04}$. \ $\endpf$

\section{Proof of Theorem \ref{003}}

In this section, we give a proof for Theorem \ref{003}\,.  Recall
that $M$ is assumed to be a real codimension two submanifold in
${\mathbb C}^{n+1}$ with $0\in M$ a CR singular point.
 Near $p=0$, $M$ is  defined by an equation as in (\ref{001}) and
(\ref{0011}). Let $M$ be CR non-minimal at its CR points.
Furthermore, we assume that the matrix $\mathcal B$ is
non-degenerate.

 It is clear that, under the above setting, $M$ is quadratically flattenable if and only if there exists an $n$ by $n$ invertible matrix $P$ and a $\mu\neq 0$ such that
\begin{equation}\label{matrix}
\wt{\mathcal B}=\frac{1}{\mu} P\cdot \mathcal B\cdot \ov{P}^t\
\end{equation} is Hermitian. (See (\ref{00111}) in $\S 2$).  We will prove Theorem \ref{003} by slicing $M$ with some special
three-dimensional complex submanifolds and thus reducing the proof
to that of Theorem \ref{04}.

{\it  Proof of Theorem \ref{003}:}  By Schur's Lemma, there is  a
unitary matrix $U$ such that $U\cdot\mathcal B\cdot \ov U^t$ is an
upper triangular matrix:
\begin{equation}
    \mathcal B^{\prime}=U\cdot\mathcal B\cdot \ov U^t=\left(\begin{matrix}
    \mu_1&\cdots&\cdots&\cdots&\cdots\\
    &\mu_2&\cdots&\cdots&\cdots\\
    &&\ddots&b_{ij}&\cdots\\
    &&&\ddots&\cdots\\
    \text{0}&&&&\mu_n\\
    \end{matrix}\right).
\end{equation}
Since we assumed that $\hbox{det}(\mathcal B)\not =0$,
$\mu_1,\mu_2,\cdots,\mu_n\neq 0.$  Choose $\mu=\mu_1, P=U$ to be as
 in (\ref{matrix}).  Without lost of generality,  we can just
assume $\mathcal B$ is  $\mathcal B^{\prime}$ with the above form
and $\mu_1=1,\mu_2,\cdots,\mu_n\neq 0$.

Now if $0\in M$ is not a quadratically flattenable CR singular
point, then either we have some $\mu_j\notin\mathbb R$ or
$\mu_2,\cdots,\mu_n\in \mathbb R$ but some element $b_{ij}(j>i)\neq
0$.

In the first situation, we let $E=\{z_k=0,k\neq1,j\}$ and
$M^{\prime}=M\cap E.$   Notice that $M^{\prime}$ is a real
codimension two submanifold in $\mathbb C^3$ with $0$ as a CR
singular point and satisfies all the assumptions made in Theorem
\ref{04}. Moreover,   the defining function for $M^{\prime}$ in
$\mathbb C^3$ (with coordinates $(z_1,z_j,w)$) takes the form that
is  similar to (\ref{429eq1}):
\begin{equation*}
w=q(z,\ov{z})+p(z,\ov{z})+iE(z,\ov{z}),
\end{equation*}
\begin{equation*}
{\text {where}\,\, }q(z,\ov{z})=2\Re(
h^{(2)}(z))+(z_1,z_j)\cdot\mathcal B_1\cdot (\ov z_1,\ov
z_j)^t,\,\,\text{and}\,\, \mathcal B_1=\left(\begin{matrix}
1&b_{1j}\\
0&\mu_2\\
\end{matrix}\right).
\end{equation*}
On the one hand, by applying Theorem \ref{04}, we know $M^{\prime}$
is quadratically flattenable at the point $0$; for  $
\hbox{det}(\mathcal B_1)\not =0$.  However, we can not find
$\mu\neq0$ and  an invertible 2 by 2 matrix $P$ such that
\begin{equation*} \wt{\mathcal B}_1=\frac{1}{\mu} P\cdot \mathcal
B_1\cdot \ov{P}^t,
\end{equation*} is Hermitian due to the fact $\mu_2\notin\mathbb R.$
This contradicts Theorem  \ref{04}.

We turn  to the second situation where $\mu_2,\cdots,\mu_n\in
\mathbb R$ but some element $b_{ij}(j>i)\neq 0$.  In this case we
let $E=\{z_k=0,k\neq i,j\}$ and $M^{\prime}=M\cap E$.  Here
$M^{\prime}$ satisfies all the assumptions made in Theorem \ref{04}
as well.  Similar to (\ref{429eq1}) and the discussion above, we
derive the defining function for $M^{\prime}$ in $\mathbb C^3$ (with
coordinates $(z_i,z_j,w)$) as
\begin{equation*}
w=q+O(3),
\end{equation*}
\begin{equation*}
{\text {where}\,\, }q=2\Re( h^{(2)}(z_i,z_j))+(z_i,z_j)\cdot\mathcal
B_2\cdot (\ov z_i,\ov z_j)^t,\,\,\text{and}\,\, \mathcal
B_2=\left(\begin{matrix}
\mu_i&b_{ij}\\
0&\mu_j\\
\end{matrix}\right).
\end{equation*}
By the same token,  we conclude that $M^{\prime}$ is quadratically
flattenable at the point $0$; for $ \hbox{det}(\mathcal B_2)\not
=0$. However for this  matrix $\mathcal B_2,$ it is clear that there
do not exist $\mu\neq0$ and an invertible 2 by 2 matrix $P$ such
that
\begin{equation*} \wt{\mathcal B}_2=\frac{1}{\mu} P\cdot \mathcal
B_2\cdot \ov{P}^t, \end{equation*} is Hermitian, due to the fact
$\mu_i,\mu_j,b_{ij}\neq 0.$  This is again a  contradiction. As a
conclusion, we complete the proof  Theorem \ref{003}.

\begin{remark}:  It is clear from our proof,  Theorem \ref{003} and
Theorem \ref{04}  hold when the submanifold $M$ is just assumed to
be ${\mathbb C}^3$-smooth.
\end{remark}
\section{Holomorphic flattening: a geometric approach}
\subsection{A general approach}  In this section, we apply Theorem \ref{003} and the work
in [HY3] to give a proof of Theorem \ref{005} when the set of CR
singular points has real dimension less than $(2n-2)$. Our approach
here is more along the lines of the geometric argument used in
[HY3]. The great benefit of this argument, compared with the formal
argument in [HY4],  is that we do not need to know the precise
structure of the quadratic normal form which is almost impossible to
obtain when $n+1>3$. The reason we want to have the set of CR
singular points has real dimension less than $(2n-2)$ is because we
need to find a good elliptic complex tangency for the sliced
manifolds.
%The observation of elliptic directions in the case of many
%non-degenerate CR singular points will also be used here.

We  recall the definition of elliptic directions first introduced in
the paper of Huang-Yin [HY3] (see [Theorem 2.3, HY3] and also the
papers [Bur1-2] and  [LNR]):
\begin{definition}\label{0.0}
Let $M$ be a smooth  submanifold with $0\in M$ being its CR singular
point. Suppose that  there is a holomorphic change of coordinates
preserving the origin
% of the special form $(z',w'):=(z,
%w+O(|zw|+|w|^2+|z|^3))$
such that  $M$ in the new coordinates (which, for simplicity, we
still write as $(z,w)\in {\mathbb C}^n\times {\mathbb C}$) is
defined by an equation of the form:
\begin{equation} \label{00-00}
w=G(z,\-{z})+iE(z,\-{z})=O(|z|^2),\
(G+iE)(z_1,0,\-{z_1},0)=|z_1|^2+\ld_1(z_1^2+\-{z}_1^2)+o(|z_1|^2).
\end{equation}
Here the constant $\ld_1$ is such that $0\le \ld_1 <\frac{1}{2}$. We
say that $M$ has an elliptic direction along the  $z_1$-direction in
the new coordinates.
\end{definition}
 Suppose $M$ in
Definition \ref{0.0} has an elliptic direction along the
$z_1$-direction. By the stability of elliptic complex tangency, then
$M_t:=M\cap \{(z_1, z''=t,w)\}$ is also an elliptic Bishop surface
in the affine plane $\{(z_1, z''=t,w)\}$ for each fixed
$t=(t_2,\cdots,t_n)\approx 0$ and has a unique elliptic CR singular
point, when regarded  as a real surface in the affine complex plane
defined above. Denote by $p(t)$ the elliptic CR singular point in
$M_t$, which is regarded as a real surface in a complex affine
plane. We need the following result whose proof is contained in
[HY3, Theorems 2.1, 2.2, 2.3] for our purpose here:
\begin{theorem} (Huang-Yin [HY3])\label{1.1} Let $M$ be a real analytic submanifold of
codimension two with $0\in M$ a CR singular point. Suppose that $M$
has an elliptic direction along the $z_1$-direction and let $p(t)$
be the elliptic CR singular point of $M_t$ defined above. Suppose
for some sufficiently small $|t|$, there is a real-analytic
Levi-flat hypersurface $H_{p(t)}$ containing $p(t)$ such that all
small holomorphic disks attached to $M$ near $p(t)$ stay inside
$H_{p(t)}$. Then $(M,0)$ can be holomorphically flattened. Namely,
there is a holomorphic change of coordinates preserving the origin
such that in the new coordinates, $M$ can be defined by a function
near the origin of the form: $w=\rho(z,\ov{z})$ with
$\Im{\rho}\equiv 0$.
\end{theorem}
Applying Theorem \ref{1.1}  and a holomorphic change of coordinates,
we have  the following corollary:
% that to get the holomorphic
%flattening of $(M,0)$
\begin{corollary} \label{11.11}
%Applying Theorem \ref{1.1}  and a holomorphic change of coordinates,
Let $(M,0)$ be a real analytic real codimension two submanifold with
$0\in M$ a CR singular point. Then $(M,0)$ can be   holomorphically
flattened if the following two conditions hold: (A) $M$ has an
elliptic direction $\vec{c}=(c_1,\cdots,c_n)\not =0$ at $0$, namely,
$M\cap \{z_j=c_j\xi,\ \xi\in {\mathbb C}\}$, when regarded as a real
surface in the $(\xi,w)$-plane,  has an elliptic complex tangency at
$0$; and (B) there is a real analytic Levi-flat hypersurface
$H_{p_{c,a}(\t)}$ passing through $p_{c,a}(\t)$  for some
sufficiently small $|\t|$, which contains all small holomorphic
disks attached to $M$ near $p_{c,a}(\t)$.  Here $p_{c,a}(\t)$ is the
elliptic complex tangent point of $M_{a,c,\t}:=M\cap
\{z_j=\vec{c}\xi+\tau_1\vec{a}_1+\cdots+ \tau_{n-1}\vec{a}_{n-1}\}$
when regarded as a real surface in the complex plane with
coordinates $(\xi,w)$, where $\{\vec{c},
\vec{a}_1,\cdots,\vec{a}_{n-1}\}$ forms a linear independent system
of vectors in ${\mathbb C}^n$ and
%=(a_1,\cdots,a_n)$ is a certain vector not proportional to the
%vector $c$ and
$\t=(\tau_1,\cdots,\tau_{n-1})\in {\mathbb C}^{n-1}$ is the
parameter.
\end{corollary}

Now we assume that $(M,p)$ is as in Theorem \ref{005}. After a
holomorphic  change of coordinates, we  assume that $p=0$ and $M$ is
defined by  a real analytic function whose quadratic term is as in
(\ref{001}) and (\ref{0011}). The rest of this section is to verify
that the hypotheses in Corollary  \ref{11.11} hold when $n+1\ge 4$
or when $(M,0)$ is not equivalent to a surface $(M',0)$ whose
quadratic term is either as in (\ref{P1}) or as in (\ref{M1}). We
first modify an argument in [HY3] to give a proof of the existence
of a Levi-flat piece as in Corollary \ref{11.11} assuming the
existence of  a generic {\it good elliptic complex point}. The
existence of these {\it good points} will be verified later this
section:

\begin{proposition} \label{3.3} Let $M\subset {\mathbb C}^{n+1}$ be  a real analytic real  codimension
two submanifold with $0\in M$ a CR singular point  and with all its
CR points non-minimal. Assume that $M$ is elliptic along the
$\vec{c}=(c_1,\cdots,c_n)$-direction. Let  $\{\vec{a}_j\in {\mathbb
C}^n\}_{j=1}^{n-1}$  be such that $\{\vec{c},
\vec{a}_1,\cdots,\vec{a}_{n-1}\}$ are linearly independent. Write,
for a parameter $\t=(\t_1,\cdots,\t_{n-1})$ with $|\t|<<1$,
 $p_{c,a}(\t)$ for   the elliptic complex tangent   point  of
 $M\cap \{z=\vec{c}\xi+\vec{a}_1\t_1+\cdots+\vec{a}_{n-1}\tau_{n-1}\}$,
  when regarded as a real surface in the $(\xi,w)$-complex plane.
 Assume that for some  $\t \approx 0$, $p_{c,a}(\t)$ is a CR point of  $M$, and $M$ does not contain
 any
 complex analytic variety  of dimension $(n-1)$ passing through
 $p_{c,a}(\t)$. Then $(M,0)$ can be holomorphically flattened.
 \end{proposition}

 \begin{remark} \label{22.221} Before we proceed to the proof of Proposition \ref{3.3}, we mention that  it is
crucial to have some $p_{c,a}(\t)$ to be a CR point of $M$. For
instance, the 4-manifold $M$ defined by $w=z_1\ov{z_2}$ in Example
3.3 has an elliptic direction $\vec{c}=(1,1)$. However, the CR
singular point of $M$ is of real dimension two, defined by $z_1=0$.
Hence, no matter how $\vec{a}$ is chosen, $p_{c,a}(\tau)$ will be a
CR singular point of $M$. Hence, the assumption for the existence of
$H_{p_{c,a}(\t)}$ can not hold. Indeed, as we know in Example 3.3,
$M$ is not even quadratically flattenable at $0$, though all CR
points of $M$ are non-minimal.
\end{remark}

{\it Proof of Proposition \ref{3.3}}: Assume, without loss of
generality, that $c=(0,\cdots,0,1)$ is the $z_n$-direction. By
Theorem \ref{1.1}, we need only to show that there is a Levi-flat
hypersurface $H_{p_{c,a}(\tau)}$ passing through $p_{c,a}(\tau)$ as
in the proposition  such that any small holomorphic disks attached
to $M$ near $p_{c,a}(\tau)$ are contained inside $H_{p_{c,a}}(\t)$.
First, since $p_{c,a}(\t)$ is a CR point of $M$, by the assumption,
 $M$ is non-minimal for any point near  $p_{c,a}(\tau)$ and for
each point $q(\in M)\approx p_{c,a}$ there is a CR submanifold $X_q$
of CR dimension $(n-1)$ through $q$. Since we assumed that there is
no complex analytic variety of dimension $(n-1)$ passing through
$p_{c,a}(\t)$, $X_p$ is of hypersurface type and is of  finite type
in the sense of Bloom-Graham [BG]. Hence, we  conclude  that
$X_q\subset M$ is a CR submanifold of hypersurface type of CR
dimension $(n-1)$, that is also of Bloom-Graham finite type. Notice
that for each $q\approx p_{c,a}(t)$ sufficiently close to $0$, the
tangent vector fields of type $(1,0)$  of $X_q$ near $q$ are spanned
by $\{L_1,\cdots, L_{n-1}\}$ with
\begin{equation}\begin{split}\label{325eq2}
L_j=&(G_n-iE_n)\frac{\p}{\p z_j}-(G_j-iE_j)\frac{\p}{\p
z_n}+2i(G_nE_j-G_jE_n)\frac{\p}{\p w}
\end{split}\end{equation}
where $G, E$ are similarly defined as in  (\ref{429eq1}) and (
\ref{82eq1}). Also, a certain fixed iterated Lie bracket $T$ from
elements in $\{L_1,\cdots, L_{n-1}, \ov{L_1},\cdots,\ov{L_{n-1}}\}$ is a
non-zero tangent vector field of $M$ near $p_{c,a}(\t)$ such that
the span of the vector fields   $\{\Re{L_j},\Im{L_j}, \Re{T}\}$
defines a real analytic distribution of real codimension one near
$p_{c,a}(\t)$.
 After
a linear (holomorphic) change of coordinates, we assume that
$p_{c,a}(\t)=0$ and $T_0M=\{ y_1=v=0\}$ and
$T_0^{(1,0)}M=\{z_1=w=0\}$, here $z_j=x_j+\sqrt{-1}y_j, w=u+\sqrt{-1}v$. Performing a linear transformation, we
also assume that $\frac{\p}{\p x_1}|_0$ is tangent to the CR
foliation at $0$ and $\frac{\p}{\p u}|_0$ is transversal to the CR
foliation of $M$ near $q=0$. Notice here that $\frac{\p}{\p x_1}|_0$ and $\frac{\p}{\p u}|_0$ are tangential to $M$ at 0. Now,   from the foliation theory, we
can find a real valued real analytic function $t(Z)$ defined over a
certain neighborhood $U_0$ of $p_{c,a}(\t)=0$ in $M$ such that for
each $t_0\in I_{\d_0}=(-\d_0, \d_0)$ with a certain small
$0<\d_0<<1$, $M_{t_0}=\{Z\in U_0, t(Z)=t_0\}$ is a connected real
analytic CR submanifold of hypersurface  type of CR dimension
$(n-1)$. Moreover $dt|_{U_0}\not = 0$. We assume that $0\in
M_{t_0=0}$. Define $\Psi: U_0\rightarrow {\mathbb C}^n\times
{\mathbb R}$ by sending $Z=(z_1,z_2,...z_n,w)\in U_0$ to
$\Psi(Z)=(z_1,z_2,...,z_n,t(Z))$.  After shrinking $U_0$ and $\d_0$ if
needed, we can assume that $\Psi$ is a real analytic embedding.
Write ${M^*}_t$ for $\Psi(M_t)$ for
 each $t\in I_{\d_0}$.  Since each component of $\Psi$ is the restriction of a holomorphic function over $M_t$,
  $\Psi$ is a CR diffeomorphism from $M_t$ to
 $M^*_t$. Write $M^*=\Psi(M)$. Then $M^*$ is a real analytic
 hypersurface in ${\mathbb C}^n\times {\mathbb R}$ and $\Psi$ is a
 real analytic differomorphism from a neighborhood  of $p_{c,a}(\t)$ in $M$
   to a certain neighborhood of  $p^*:=\Psi(p_{c,a}(\t))=({p^*}',0)$ in  $M^*$. Also
   $\Phi:=\Psi^{-1}$ is  a CR diffeomorphism when restricted to each
   $M ^*_t$ near $p^*=\Psi(p_{c,a}(\t))$. Now,  the real analytic CR function $\Phi$ extends
   holomorphically to a certain fixed neighborhood of ${p^*}'$ in
   ${\mathbb C}^n$ for each fixed $t$. Also the holomorphic
   extension depends real analytically on its parameter $t$ from
   the way the holomorphic extension is constructed as demonstrated  by the following lemma:
\begin{lemma}\label{22.22} Let $M_t\subset {\mathbb C}^{n}$ be a
real analytic family of real analytic hypersurfaces with a real
parameter $t$ defined by $\Im{w}=\rho(z,\ov{z},\Re{w};t)$. Here
$\rho(z,\ov{z},\Re{w};t) $ is a real analytic function in
$(z,\ov{z},\Re{w},t)$ in  a neighborhood of the origin with
$\rho(0,0,0;0)=0,\ d\rho|_0=0$. Suppose that $f$ is a real analytic
function in $(z,w,t)$ near the origin and is a CR function when
restricted to each $M_t$ near $0$. Then $f$ extends to a real
analytic function in a neighborhood of $(0',0)\in {\mathbb
C}^n\times {\mathbb R}$ which is holomorphic for each fixed small
$t$.
\end{lemma}
{\it Proof of Lemma \ref{22.22}}: We first notice that after
applying a holomorphic transformation of the form:
$(z,w)=H(z',w';t)$, where $H(z',w';t)$ is real analytic in
$(z',w',t)$  in a small neighborhood of the origin and is
biholomorphic for each fixed $t$, and by applying an implicit
function theorem (see [\S3, BR], for instance), we can assume that
$M_t$ is defined near the origin by an equation of the form:
$\Im{w}=\rho(z,\ov{z},\Re{w};t)$ or $\ov{w}=\wt{\rho}(z,\ov{z},w;t)$.
Here, $\rho,\ \wt{\rho}$ can be expanded to a Taylor series at the
origin in $(z,\ov{z},\Re{w},t)$ (or $(z,\ov{z},w,t)$, respectively)
and $\rho|_0=\wt{\rho}|_0=0, d\rho_0=d\wt{\rho}|_0=0,\
\rho(z,0,\Re{w};t)= \rho(0,\ov{z},\Re{w};t)=0,\
\wt{\rho}(z,0,w;t)=\wt{\rho}(0,\ov{z},w;t)=0$. Now, let
$f(z,\ov{z},w,\ov{w};t)$ be a real analytic function near $0$ and is
CR when restricted to each $M_t$. Notice that the CR vector field of
$M_t$ is given by $\ov{L_j}=\frac{\partial}{\partial
\ov{z}_j}+\wt{\rho}_{\ov{z_j}}\frac{\partial}{\partial \ov{w}}$,\
$j=1,\cdots, n-1$. Define
$\wt{f}(z,\ov{z},w,t)=f(z,\ov{z},w,\wt{\rho}(z,\ov{z},w;t);t)$. Since
$f$ is CR along each $M_t$, we obtain
 $$\frac{\partial \wt{f}}{\partial
\ov{z}_j}=\frac{\partial f}{\partial \ov{z}_j}+\frac{\partial
f}{\partial \ov{w}}\wt{\rho}_{\ov{z_j}}=0,\ j=1,\cdots, (n-1)\ \
\hbox{along each } \ M_t.$$ Notice that $\wt{\rho}_{\ov{z_j}}$ is
holomorphic in $(z,\ov{z},w,t)$. $\frac{\partial \wt{f}}{\partial
\ov{z}_j}= F_j(z,\ov z,
w,\ov{w},t)|_{\ov{w}=\wt{\rho}(z,\ov{z},w;t)}$ for a certain $F_j$
holomorphic in its variables. Hence we see, in the same way, that
$\frac{\partial^2 \wt{f}}{\partial \ov{z}_j{\partial\ov{z_k}}}=0$
along $M_t$. Inductively,  we see that
  for each multi-indices
$\a=(\a_1,\cdots,\a_{n-1})$ and $\b=(\b_1,\cdots,\b_{n-1})$ with
$|\b|\ge 1$, we have ${L^{\a}}\frac{\p^{|\b|} \wt f}{\p \-
z^{\b}}=0$ along each $M_t$. Evaluating this at $(z,w)=(0,\Re{w})$
and making use of the normalization condition of $\wt{\rho}$, we
obtain that
$$\frac{\partial^{|\a|+|\b|} \wt f}{\partial
    {z}^\a \partial
\ov{z}^\b}(0,0,\Re{w};t)=0,\ \ (0,\Re{w};t)\in M_t,\  |\b|\ge 1.$$
Hence, from the facts that $\frac{\partial^{|\a|+|\b|} \wt f}{\partial
    {z}^\a \partial
    \ov{z}^\b}(0,0,\Re{w};t)=0$ and it is holomorphic in $w$, it follows that $\frac{\partial^{|\a|+|\b|} \wt f}{\partial
    {z}^\a \partial
    \ov{z}^\b}(0,0,w;t)\equiv 0$
for $|\b|\ge 1$. From this, one can easily conclude that in the
Taylor expansion of $\wt{f}(z,\ov{z},w;t)$ in $(z,\ov{z},w,t)$ at
$0$, there are no $\ov{z}$-terms. Hence, $\wt{f}$ gives the desired
holomorphic extension of $f$ to a certain  fixed neighborhood of the
origin and is real analytic  as a function  in the joint variables
$(z,w,t)$. \ \ $\endpf$

\medskip
   We now continue the proof of the proposition. Complexifying
   $t$, we then get a holomorphic extension of $\Psi$ to
 a neighborhood of $p^*$ in    ${\mathbb C}^n\times {\mathbb C}$, that
 is biholomorphic near $p^*$. Then define $X_{p}$ to be the
 biholomorphic image of ${\mathbb C}^n\times {\mathbb R}$ near
 $p^*$. Then $X_{p}$ is a Levi flat hypersurface containing a small
 piece of $M$ near
 $p_{c,a}(\t)$, which certainly  contains all small holomorphic disks attached
 $M$ near $p_{c,a}(\t)$. By Theorem \ref{1.1}, we thus complete the proof
 of Proposition \ref{3.3}. $\endpf$

\subsection{The case of $n+1\ge 4$}

For the rest of this section, we will verify the hypotheses stated
in Proposition  \ref{3.3} and thus complete the proof of Theorem
\ref{005} except for one case which  will be handled in the next
section.

We now assume that $p\in M$ is a non-degenerate CR singular point.
By Theorem \ref{003}, there is a holomorphic change of variables
such that in the new coordinates, $p=0$ and $M$  is defined by an
equation as in (\ref{001}). Moreover ${\mathcal B}$ is a
non-degenerate Hermitian matrix. Hence, after a linear change of
variables, we can further assume that ${\mathcal
B}=\hbox{diag}(1,\cdots, 1, -1,\cdots, -1)$. Therefore, $M$ is
defined by a real analytic equation of the form:
\begin{equation}\label {123}
w=F(z,\-{z})=\sum_{j=1}^{\ell}|z_j|^2-\sum_{j=\ell+1}^{n}|z_j|^2+2\Re\left(\sum_{j,k=1}^n
a_{jk}z_jz_k\right)+O(|z|^3). \end{equation}
 Here we can always
assume, without loss of generality, that $\ell\ge n/2$. Write
${\mathcal S}$ for the CR singular points of $M$. Then over
${\mathcal S}$, we have $\frac{\partial F}{\partial \ov{z}}=0$. From
the implicit function theorem, we then can solve that
$z=\phi(\ov{z})$. Hence, ${\mathcal S}$ is a totally real analytic
variety of real dimension at most $n$. We next give the following
lemma:

\begin{lemma}\label{5.5}Suppose that $M$ is a real codimension two smooth submanifold in $\mathbb C^{n+1}$
 with coordinates $(w,z_1,\cdots,z_n)$, that is  defined by
 $w=\sum_{j=1}^{l}|z_j|^2-\sum_{j=l+1}^n|z_j|^2+2\Re{({z}\cdot {\mathcal A}\cdot {z}^t)}+E$ with $E=O(|z|^3)$.
 Then there is  a small neighborhood $U\subset M$ of the origin
such that there is no $r$ dimensional complex analytic subvariety
contained in $U$ for any $r>[\frac{n}{2}]$.
\end{lemma}
{\it Proof of Lemma \ref{5.5}}:   If for any small neighborhood of
$0\in M$,  there is a $r-$dimensional complex analytic variety
contained in M, then we can choose a sequence of points
$\{p_i\}_{i=1}^{\infty}\subset M$  converging to $0$ such that at
each $p_i$, there is a $r-$dimensional complex submanifold $\mathcal
V_{p_i}\subset M$ through $p_i.$ We denote the tangent space of
$\mathcal V_{p_i}$ at $p_i$ by $H_{p_i}$.  Write $M^*$ for the real
hypersurface defined by  $\rho:= -\frac{1}{2}w-\frac{1}{2}\ov
w+\sum_{j=1}^{l}|z_j|^2-\sum_{j=l+1}^n|z_j|^2+2\Re{({z}\cdot
{\mathcal A}\cdot {z}^t)}+\Re{E}$. Then $M\subset M^*$ and the Levi
form of $M^*$ is given by $L(Y,Y)=\partial\ov\partial\rho(Y,\ov{Y})$
when restricted to the holomorphic tangent space of $M^*$. Notice
that the Levi form is zero when restricted to $H_{p_i}$.  We then
extract a subsequence $\{p^{\prime}_{j}\}_{j=1}^{\infty}$ of
$\{p_i\}_{i=1}^{\infty}$ such that $p'_j\rightarrow 0$ and
$H_{p'_j}\rightarrow H_0\subset T^{(1,0)}_0M^*.$ Here we view
$H_{p_j}$ and $H_0$ as $r-$dimensional complex vector spaces in
$\mathbb C^{n}$.
% through
%the natural embedding $M\hookrightarrow\mathbb C^{n+1}.$
%It is easy
%to see $H_0\subset T^{(1,0)}_0M.$
By continuity,  the Levi form $L$ is zero when restricted to $H_0.$
On the other hand, the naturally associated matrix for the  Levi
form of $M^*$ at the origin is given by
\begin{equation*}
\mathcal L=\left(\begin{matrix}
I_l&\\
&&-I_{n-l}\\
\end{matrix}\right).
\end{equation*}
Write a tangent vector as a column vector and represent  a basis of
$H_0$ by  an $n\times r$ matrix $H$. Then
$$\ov H^t\mathcal L H=0.$$ There is a certain invertible $r\times r$
matrix $P$ such that $H^{\prime}=HP$ has the following form:
\begin{equation*}
H^{\prime}=\left[
\begin{array}{c;{2pt/2pt}c}
A & B \\ \hdashline[2pt/2pt] C & 0
\end{array}
\right],
\end{equation*}
where $C$ is an $(n-l)\times k$ matrix with rank $k$ (hence $k\leq
n-l$).  Then since $\ov H^t\mathcal L H=0,$ we have
%\begin{equation*}
%\ov B^t \left(\begin{matrix}
%0&\\
%&I_r\\
%\end{matrix}\right)B=0.
%\end{equation*}
$\ov{B}^t\cdot B=0$. Hence $B=0$. This shows that $r= k$ and thus
$r=k\le n-l\le \frac{n}{2}.$ $\endpf$
% Then either $p=q$ or $B=[1,0,\cdots,0]^t.$
%However $H_0\subset TM_0$ and $\frac{\partial}{\partial w}\notin
%TM_0.$ We thus have $p=q\leq n-r.$ $\endpf$
\begin{corollary} \label{6.6} When $n+1\geq 4$, there is no $(n-1)-$dimensional complex analytic subvariety that is contained in a  small neighborhood of $0\in M.$
\end{corollary}

We now can complete the proof of Theorem \ref{005} for $n+1\ge 4$.
In this case, a simple  algebra shows that $M$ always has an
elliptic direction. (See [Lemma 3.1, LNR], for instance). Indeed,
setting $z_j=0$ for $j>\ell$ and diagonalizing the harmonic
quadratic part in (\ref{123}),   we need only to show that an
$M\subset {\mathbb C}^3$ defined by an equation of the form
$$w=|z_1|^2+|z_2|^2+\ld_1(z_1^2+\ov{z}_1^2)+\ld_2(z_2^2+\ov{z}_2^2)+O(|z|^3),\
\ld_1,\ld_2\ge 0$$ has an elliptic direction passing through the
origin. Indeed, if $\ld_1=\ld_2$, $M$ is elliptic along the
direction $\vec{c}=(1, \sqrt{-1})$. Otherwise, assume that
$\ld_1<\ld_2$. Then the elliptic direction can be simply taken as
$\vec{c}=(\ld_2, \sqrt{-1}\ld_1).$ Now when we move along the
transversal directions to the elliptic direction, we get a real
$2(n-1)$-families of elliptic Bishop surfaces. And the elliptic
complex tangent points thus obtained form a smooth manifold of real
dimension $2(n-1)>\max\{n, 2(n-2)\}$. (See the last equation on page
394 of [HY3]).  Hence, a generic elliptic singular point stays at a
CR point of $M$ through which there is no complex analytic
subvariety of complex dimension $(n-1)$ contained in $M$. By
Proposition \ref{3.3}, we see the proof of Theorem \ref{005} when
$n+1\ge 4$. $\endpf$

\subsection{The case of $n+1$=3}
 We next study the existence of good elliptic points  when $n+1=3$. We
start with the following:
\begin{lemma}\label{132}
    Let  $M$ be a $4-$dimensional real analytic submanifold of $\mathbb C^3$  defined by
    $$w=2\Re{({z}\cdot {\mathcal A}\cdot {z}^t)}+{z}\cdot
    {\mathcal B}\cdot \ov{{z}}^t+O(|z|^3).$$  Assume  that $\{{\mathcal A},{\mathcal B}\}$ take one of the following forms from part of the list given   in $\S 2$:

\noindent (A).
$ \mathcal B=\left(\begin{matrix}1&0\\
0&1
\end{matrix}\right),\,\,\,\,\,\,
\mathcal A=\left(\begin{matrix}\lambda_1&0\\0&\lambda_2
\end{matrix}\right), \,\,\lambda_1,\lambda_2\geq 0$ and $\lambda_1\neq\frac{1}{2}$ or
$\lambda_2\neq\frac{1}{2}\,\,\,\,(\ref{5})$;

\noindent (B).
$ \mathcal B=\left(\begin{matrix}1&0\\
0&-1
\end{matrix}\right),\,\,\,\,
\mathcal A=\left(\begin{matrix}\lambda_1&0\\0&\lambda_2
\end{matrix}\right), \,\,$
$0\leq\lambda_1\leq\lambda_2$ and $\lambda_1<\frac{1}{2}$ or
$\frac{1}{2}\leq\lambda_1<\lambda_2\,\, (\ref{6a})$;

\noindent (C).
$ \mathcal B=\left(\begin{matrix}1&0\\
0&-1
\end{matrix}\right),\,\,\,\,\,\,
\mathcal A=\left(\begin{matrix}0&\lambda\\\lambda&0
\end{matrix}\right), \,\,\lambda>0,\quad\quad\quad\quad\quad\quad\quad\quad\quad\quad\quad\quad\,\, $ $(\ref{6b})$;

\noindent (D).
$ \mathcal B=\left(\begin{matrix}1&0\\
0&-1
\end{matrix}\right),\,\,\,\,\,\,
\mathcal
A=\left(\begin{matrix}\frac{1}{2}&\frac{1}{2}\\\frac{1}{2}&\frac{1}{2}
\end{matrix}\right),\quad\quad\quad\quad\quad\quad\quad\quad\quad\quad\quad\quad\quad\quad\quad\,\,$  $(\ref{6c})$;

\noindent (E).
$\mathcal B=\left(\begin{matrix}0&1\\
1&0
\end{matrix}\right),\,\,\,\,\,\,
\mathcal A=\left(\begin{matrix}0&b\\b&\frac{1}{2}
\end{matrix}\right),\,\, b>0,\quad\quad\quad\quad\quad\quad\quad\quad\quad\quad\quad\quad\quad\,\,$  $(\ref{7a})$;

\noindent (F).
$\mathcal B=\left(\begin{matrix}0&1\\
1&0
\end{matrix}\right),\,\,\,\,\,\,
\mathcal A=\left(\begin{matrix}\frac{1}{2}&0\\0&d
\end{matrix}\right),\,\, {\rm Im} d>0,\quad\quad\quad\quad\quad\quad\quad\quad\quad\quad\quad\quad$  $(\ref{7b})$;

\noindent Then $(M,0)$ has an elliptic direction $\vec{c}$ at $0$.
Also there is a direction $\vec{a}$ transversal to $\vec{c}$ such
that for the  family of complex affine planes $V(\tau)\subset
\mathbb C^3:=\{(z,w): z=\vec{c}\xi +\vec{a}\t,\ \xi, \t \in {\mathbb
C}\}$, being parametrized by $\tau (\approx 0)\in\mathbb C$ with
$0\in V(0)$, the elliptic complex tangent  point $P(\t)$ of
$V(\tau)\cap M$, for a generic $\t\approx 0$, is a CR point in $M$.
Moreover, there is no holomorphic curve contained in $M$ passing
through $P(\t)$ for a generic $\t$.
\end{lemma}
Assuming this lemma,  when the quadratic normal form of $M$ is not
of the type in (\ref{M1}),  then the proof of Theorem \ref{005} for
$n+1=3$ follows immediately from Proposition \ref{3.3}, the
classification of the quadratic terms as listed in $\S 3$ and Lemma
\ref{132}. We now proceed to the proof of Lemma \ref{132}.

 {\it Proof of Lemma \ref{132}:} We will give the proof of
the lemma based on a case by case argument  in terms of  the normal
forms  listed in the lemma.

 Case (A) : $w=|z_1|^2+|z_2|^2+\lambda_1(z_1^2+\ov
z_1^2)+\lambda_2(z_2^2+\ov z_2^2)+E(z,\ov z)$  with $\lambda_1,\lambda_2\geq 0$ and
$\lambda_1\neq\frac{1}{2}$ or $\lambda_2\neq\frac{1}{2}$.
 It is contained in a strongly pseudoconvex hypersurface and thus there is no non-trivial holomorphic curve contained in
 $M$. The set ${\mathcal S}$  of the complex tangent points in this case are  of at most real
 dimension one. Indeed, the set of CR singular points in $M$ is
 defined by the equation:
\begin{equation*}
0=\frac{\partial w}{\partial \ov z_1}=z_1+\lambda_1\ov z_1+o(|z|),
\quad 0=\frac{\partial w}{\partial \ov z_2}=z_2+\lambda_2\ov
z_2+o(|z|)
\end{equation*}
Suppose that $\ld_1\not = 1/2$. Then one can solve $x_1, y_1$ from
the first equation in terms of $x_2,y_2$. Here
$x_j+\sqrt{-1}y_j=z_j$. Substituting into the second equation, we
can always  solve  $x_2$ in terms of $y_2$ and thus ${\mathcal S}$
has real dimension at most one. In fact, if $\ld_2$ is not $1/2$
neither, then ${\mathcal S}$ has an isolated point at $0$.

 As mentioned in the above subsection,  there is an elliptic direction
 through the origin in this case. Now, we choose  a holomorphic family of  affine  complex
 planes transversal to the elliptic direction. Then the resulting  elliptic complex tangent points form a real dimensional two subset. Hence, we can easily find the {\it good points}
 as required in Proposition  \ref{3.3}.

Case (B):  Now, the manifold is defined by a real analytic equation
of the form: $$w=|z_1|^2-|z_2|^2+\lambda_1(z_1^2+\ov
z_1^2)+\lambda_2(z_2^2+\ov z_2^2)+E(z,\ov z),\ \ \hbox{with}\
E=O(|z|^3).$$ We first assume that $0\leq\lambda_1<\frac{1}{2}$.
Then $\vec{c}=(1,0)$ is an elliptic direction and we  construct the
family $V(\tau)$ by intersecting $M$ with the affine complex plane:
$z_1=\xi,z_2=\tau$, where $\tau$ is a complex parameter. (Namely,
$\vec{c}=(1,0)$ and $\vec{a}=(0,1)$). Then ${\mathcal S}$ is defined
by the following system of equations:
\begin{equation}\label{B1}\left\{
\begin{matrix}
0=\frac{\partial w}{\partial \ov z_1}=&z_1+2\lambda_1\ov z_1+o(|z|)\\
0=\frac{\partial w}{\partial \ov z_2}=&-z_2+2\lambda_2\ov z_2+o(|z|)\\
0=\frac{\partial\ov  w}{\partial  z_1}=&\ov z_1+2\lambda_1 z_1+o(|z|)\\
0=\frac{\partial \ov w}{\partial  z_2}=&-\ov z_2+2\lambda_2 z_2+o(|z|)\\
\end{matrix}\right.\quad\quad.
\end{equation}
Since $\lambda_1\neq\frac{1}{2}$,  as  in case (A), we can solve for
$x_1,y_1$ in terms of $x_2,y_2$ through the first and the third
equation. Substituting back, we can get at least $x_2$ in terms of
$y_2$.
%the rank of the above matrix is at least three.
Therefore, $\mathcal S$ has real dimension at most one.
% Thus the
%Jacobian matrix with respect to $(z_1,\ov{z_1}, z_2,\ov{z_2})$ at
%the origin is given by
%\begin{equation}\label{jac}\left(
%\begin{matrix}
%1&2\lambda_1&0&0\\
%2\lambda_1&1&0&0\\
%0&0&-1&2\lambda_2\\
%0&0&2\lambda_2&-1\\
%\end{matrix}\right).
%\end{equation}
%Since $\lambda_1\neq\frac{1}{2}$,  the rank of the above matrix is
%at least three.  Therefore by implicit function theorem, $\mathcal
%S$ has real dimension  at most one.

At points where there is a smooth holomorphic curve  passing
through, it holds that $[L,\ov L]\subset{\rm Span}\{L,\ov L\}.$
Hence,  there are   constants $\alpha$ and $\beta$ such that $[L,\ov
L]=\alpha L+\beta \ov L.$ By $(\ref{325eq2})$ and $(\ref{325eq3})$,
we have
\begin{equation*}\begin{array}{lll}
A=G_2=-\ov z_2+2\lambda_2z_2+o(|z|),\quad B=G_1=\ov z_1+2\lambda_1 z_1+o(|z|),  \quad C=o(|z|), \\
\lambda_{(1)}=G_1,\quad\lambda_{(2)}=-G_2,\quad\lambda_{(4)}=-G_{\ov 1},\quad\lambda_{(5)}=G_{\ov 2},\quad\lambda_{(3)}=o(|z|),\quad\lambda_{(6)}=o(|z|).\\
\end{array}\end{equation*}
Substituting the above into  $[L,\ov L]=\alpha L+\beta \ov L,$  we
have
$$G_1\frac{\partial}{\partial\ov z_1}-G_2\frac{\partial}{\partial\ov z_2}-G_{\ov 1}\frac{\partial}{\partial z_1}+G_{\ov 2}\frac{\partial}{\partial z_2}
=\alpha(G_{ 2}\frac{\partial}{\partial
z_1}-G_{1}\frac{\partial}{\partial z_2})+\beta(G_{\ov
2}\frac{\partial}{\partial\ov z_1}-G_{\ov
1}\frac{\partial}{\partial\ov z_2}).$$ In particular, $G_1G_{\ov
1}=G_2G_{\ov 2}$. On the other hand, the set of the  elliptic CR
singular points of the Bishop surfaces  $V(\tau)$  is given by the
equation $\frac{\partial w}{\partial\ov \xi}=0.$  We substitute
$z_1=\xi,z_2=\tau$ in $\frac{\partial w}{\partial\ov \xi}=0$ and
$G_1G_{\ov 1}=G_2G_{\ov 2}$ to get
\begin{equation*}\left\{
\begin{matrix}
&(1+4\lambda_1^2)|\xi|^2+2\lambda_1(\xi^2+\ov \xi^2)=(1+4\lambda_2^2)|\tau|^2-2\lambda_2(\tau^2+\ov \tau^2)+o({|\xi^2|+|\tau|^2})\\
&\ov\xi+2\lambda_1\xi+o(\sqrt{|\xi^2|+|\tau|^2})=0\\
&\xi+2\lambda_1\ov\xi+o(\sqrt{|\xi^2|+|\tau|^2})=0\\
\end{matrix}\right.
\end{equation*}
Since $\ld_1\neq 1/2$, from the last two equations, we can solve
$\xi=\phi(\tau,\ov \tau)$ with $\phi(\tau,\ov \tau)=o(\tau).$  Then
the first equation can not hold for a generic $\tau$ by comparing
the second order terms in its  Taylor expansion at the origin.

We now assume that $\frac{1}{2}\leq\lambda_1<\lambda_2$.  Then
$\vec{c}=(1,i\sqrt{\frac{\lambda_1}{\lambda_2}})$ is an elliptic
direction and we construct the family $V(\tau)$ by setting
$z_1=\xi,z_2=i\sqrt{\frac{\lambda_1}{\lambda_2}}\xi+\tau$. (Namely
$\vec{a}=(0,1)$).  Then $\mathcal S$ is still defined by the  system
of equations in (\ref{B1}).
% whose Jacobian at the origin is given by (\ref{jac}).
Since $\lambda_2\neq\frac{1}{2}$,
%the rank of the Jacobian matrix
%is at least three.  Therefore
 by implicit function theorem,
$\mathcal S$ has an isolated point at $0$.

At points where there is a smooth holomorphic curve  passing
through, following the argument given above, it holds that
$G_1G_{\ov 1}=G_2G_{\ov2}$ as well. Also, the set of the elliptic CR
singular points of the Bishop surfaces $V(\tau)$ is given by the
equation $\frac{\partial w}{\partial\ov \xi}=0.$ We substitute
$z_1=\xi,z_2=i\sqrt{\frac{\lambda_1}{\lambda_2}}\xi+\tau$ into
$\frac{\partial w}{\partial\ov \xi}=0$ and $G_1G_{\ov 1}=G_2G_{\ov
2}$ to get
\begin{equation*}\left\{
\begin{matrix}
&(1+4\lambda_2^2)\{|\tau|^2+\frac{\lambda_1}{\lambda_2}|\xi^2|+i\sqrt{\frac{\lambda_1}{\lambda_2}}(-\tau\ov\xi+\ov\tau\xi)\}-2\lambda_2\{\tau^2+\ov \tau^2+2i\sqrt{\frac{\lambda_1}{\lambda_2}}(\tau\xi-\ov\tau\ov\xi)\}\\
&-(1+4\lambda_1^2)|\xi|^2+o({|\xi^2|+|\tau|^2})=0\\
&(1-\sqrt{\frac{\lambda_1}{\lambda_2}})\xi+i\sqrt{\frac{\lambda_1}{\lambda_2}}\tau-2i\sqrt{\lambda_1\lambda_2}\ov\tau+o(\sqrt{|\xi^2|+|\tau|^2})=0\\
&(1-\sqrt{\frac{\lambda_1}{\lambda_2}})\ov\xi-i\sqrt{\frac{\lambda_1}{\lambda_2}}\ov\tau+2i\sqrt{\lambda_1\lambda_2}\tau+o(\sqrt{|\xi^2|+|\tau|^2})=0\\
\end{matrix}\right.
\end{equation*}
Since $\ld_1<\ld_2$, from the last two equations, we can solve
$\xi=\phi(\tau,\ov \tau)$ with $\phi(\tau,\ov
\tau)=\frac{-i\sqrt{\lambda_1}}{\sqrt{\lambda_2}-\sqrt{\lambda_1}}\tau+\frac{2i\sqrt{\lambda_1}\lambda_2}{\sqrt{\lambda_2}-\sqrt{\lambda_1}}\ov\tau+o(\tau).$
Substituting $\xi=\phi$  into the first equation and comparing  the
coefficient for the $\tau^2$ term, we get
$$\frac{1}{\sqrt{\lambda_2}-\sqrt{\lambda_1}}\{4\lambda_1^2\lambda_2(\sqrt{\lambda_2}+\sqrt{\lambda_1})+(\lambda_1+2\lambda_2)(\sqrt{\lambda_2}-\sqrt{\lambda_1})+2\lambda_1
\sqrt{\lambda_2}(1+4\lambda_2^2)\}$$ which is not zero due to the
fact that  the second factor is strictly positive. Hence, we see
that the first equation does not hold for a generic $\tau.$ Now, as
before, the set ${\mathcal P}$ of the elliptic CR singular points
associated with the family form a real analytic subset of dimension
two, while the set of  CR singular points of $M$ is isolated at $0$
and for a generic point in ${\mathcal P}$ there is no holomorphic
curve in $M$ passing through a generic point in ${\mathcal P}$.
Hence a generic point in ${\mathcal P}$ satisfies the property
stated in the lemma.

\medskip
Case (C): In this case, $M$ is defined by
$w=|z_1|^2-|z_2|^2+\lambda(z_1z_2+\ov z_1\ov z_2)+E(z,\ov z).$ Then
$\vec{c}=(1,0)$ is an elliptic direction and we construct the family
$V(\tau)$ by setting $z_1=\xi,z_2=\tau$. Then $\mathcal S$ is
defined by the following system of equations:
\begin{equation*}\left\{
\begin{matrix}
0=\frac{\partial w}{\partial \ov z_1}=&z_1+\lambda\ov z_2+o(|z|)\\
0=\frac{\partial w}{\partial \ov z_2}=&-z_2+\lambda\ov z_1+o(|z|)\\
0=\frac{\partial\ov  w}{\partial  z_1}=&\ov z_1+\lambda z_2+o(|z|)\\
0=\frac{\partial \ov w}{\partial  z_2}=&-\ov z_2+\lambda z_1+o(|z|)\\
\end{matrix}\right.\quad\quad.
\end{equation*}
Thus the Jacobian matrix at the origin is given by
\begin{equation*}\left(
\begin{matrix}
1&0&0&\lambda\\
0&\lambda&-1&0\\
0&1&\lambda&0\\
\lambda&0&0&-1\\
\end{matrix}\right).
\end{equation*}
It is clear that the  matrix is invertible.  Therefore by the
implicit function theorem, $\mathcal S$ consists of a single point
near 0.

At points where there is a smooth holomorphic curve  passing
through, following the  argument given for Case (B), we have
\begin{equation*}\begin{array}{lll}
A=G_2=-\ov z_2+\lambda z_1+o(|z|),\quad B=G_1=\ov z_1+\lambda z_2+o(|z|),  \quad C=o(|z|), \\
\lambda_{(1)}=G_1,\quad\lambda_{(2)}=-G_2,\quad\lambda_{(4)}=-G_{\ov 1},\quad\lambda_{(5)}=G_{\ov 2},\quad\lambda_{(3)}=o(|z|),\quad\lambda_{(6)}=o(|z|),\\
\end{array}\end{equation*}
and
%$$G_1\frac{\partial}{\partial\ov z_1}-G_2\frac{\partial}{\partial\ov z_2}-G_{\ov 1}\frac{\partial}{\partial z_1}+G_{\ov 2}\frac{\partial}{\partial z_2}=\alpha(G_{ 2}\frac{\partial}{\partial z_1}-G_{1}\frac{\partial}{\partial z_2})+\beta(G_{\ov 2}\frac{\partial}{\partial\ov z_1}-G_{\ov 1}\frac{\partial}{\partial\ov z_2}).$$
%In particular,
$G_1G_{\ov 1}=G_2G_{\ov 2}$.

Similarly, the intersection  of the set of the elliptic CR singular
points of the Bishop surfaces  $V(\tau)$ with smooth points of
holomorphic curves contained in $M$  is given by the following
system:
%equation $\frac{\partial w}{\partial\ov \xi}=0.$  We
%substitute $z_1=\xi,z_2=\tau$ in $\frac{\partial w}{\partial\ov
%\xi}=0$ and $G_1G_{\ov 1}=G_2G_{\ov 2}$ to get
\begin{equation*}\left\{
\begin{matrix}
&(1-\lambda^2)|\tau|^2=(1-\lambda^2)|\xi|^2+2\lambda(\tau\xi+\ov \tau\ov\xi)+o({|\xi^2|+|\tau|^2})\\
&\xi+\lambda\ov\tau+o(\sqrt{|\xi^2|+|\tau|^2})=0\\
&\ov\xi+\lambda\tau+o(\sqrt{|\xi^2|+|\tau|^2})=0\\
\end{matrix}\right.
\end{equation*}
From the last two equations, we can solve  $\xi=\phi(\tau,\ov \tau)$
with $\phi(\tau,\ov \tau)=-\lambda\ov\tau+o(\tau).$  Substituting
this into the first equation, we have
$(1+\lambda^2)^2|\tau|^2+o(|\tau|^2)=0.$  Thus the system can not
hold for  $0\not =|\tau|\approx 0.$ It follows that a generic
elliptic CR singular point of $V(\tau)$ is a CR point of $M$ which
no holomorphic curves contained in $M$ can pass through.

\medskip
Case (D): Now, $M$ is defined by
$w=|z_1|^2-|z_2|^2+\frac{1}{2}(z_1^2+\ov z_1^2+z_2^2+\ov
z_2^2)+(z_1z_2+\ov z_1\ov z_2)+E(z,\ov z)$ with $E=O(|z|^3)$.  Then
$\vec{c}=(1,-1+\epsilon)$ with $0<\epsilon<<1$ is an elliptic
direction and we construct the family $V(\tau)$  by letting
$z_1=\xi,z_2=(-1+\epsilon)\xi+\tau$.  Then $\mathcal S$ is defined
by  the following system of equations:
\begin{equation*}\left\{
\begin{matrix}
0=\frac{\partial w}{\partial \ov z_1}=&z_1+\ov z_1+\ov z_2+o(|z|)\\
0=\frac{\partial w}{\partial \ov z_2}=&-z_2+\ov z_2+\ov z_1+o(|z|)\\
0=\frac{\partial\ov  w}{\partial  z_1}=&z_1+\ov z_1+ z_2+o(|z|)\\
0=\frac{\partial \ov w}{\partial  z_2}=&-\ov z_2+ z_2+ z_1+o(|z|)\\
\end{matrix}\right.\quad\quad.
\end{equation*}
Its Jacobian matrix at the origin is given by
\begin{equation*}\left(
\begin{matrix}
1&1&0&1\\
0&1&-1&1\\
1&1&1&0\\
1&0&1&-1\\
\end{matrix}\right),
\end{equation*}
which is invertible.  Therefore by the implicit function theorem,
$\mathcal S$ consists of a single point near 0.

At points where there is a smooth holomorphic curve  passing
through, it holds that
%$[L,\ov L]\subset{\rm Span}\{L,\ov L\}.$
%Following the argument given  above, we have
\begin{equation*}\begin{array}{lll}
A=G_2=-\ov z_2+z_2+z_1+o(|z|),\quad B=G_1=\ov z_1+z_1+ z_2+o(|z|),  \quad C=o(|z|), \\
\lambda_{(1)}=G_1,\quad\lambda_{(2)}=-G_2,\quad\lambda_{(4)}=-G_{\ov 1},\quad\lambda_{(5)}=G_{\ov 2},\quad\lambda_{(3)}=o(|z|),\quad\lambda_{(6)}=o(|z|),\\
\end{array}\end{equation*}
and
%$$G_1\frac{\partial}{\partial\ov z_1}-G_2\frac{\partial}{\partial\ov z_2}-G_{\ov 1}\frac{\partial}{\partial z_1}+G_{\ov 2}\frac{\partial}{\partial z_2}=\alpha(G_{ 2}\frac{\partial}{\partial z_1}-G_{1}\frac{\partial}{\partial z_2})+\beta(G_{\ov 2}\frac{\partial}{\partial\ov z_1}-G_{\ov 1}\frac{\partial}{\partial\ov z_2}).$$
%In particular,
$G_1G_{\ov 1}=G_2G_{\ov 2}$. Also, the intersection of the set of
the elliptic CR singular points of the Bishop surfaces $V(\tau)$
with smooth points of holomorphic curves contained in $M$ is given
by the following system:
% On the other hand, the elliptic
%points of the Bishop surfaces $V(\tau)\cap M$ is given by the
%equation $\frac{\partial w}{\partial\ov \xi}=0.$  We substitute
%$z_1=\xi,z_2=-\frac{1}{2}\xi+\tau$ into $\frac{\partial
%w}{\partial\ov \xi}=0$ and $G_1G_{\ov 1}=G_2G_{\ov 2}$ to get
\begin{equation*}\left\{
\begin{matrix}
&(\frac{1}{2}\xi+\tau)^2+(\frac{1}{2}\ov\xi+\ov\tau)^2+|\xi|^2=|-\frac{1}{2}\xi+\tau|^2+o({|\xi^2|+|\tau|^2})\\
&3\xi+\ov\xi+2\tau+2\ov\tau+o(\sqrt{|\xi^2|+|\tau|^2})=0\\
&3\ov\xi+\xi+2\ov\tau+2\tau++o(\sqrt{|\xi^2|+|\tau|^2})=0\\
\end{matrix}\right.
\end{equation*}
From the last two equations, we can solve  $\xi=\phi(\tau,\ov \tau)$
with $\phi(\tau,\ov
\tau)=-\frac{1}{2}\ov\tau-\frac{1}{2}\tau+o(\tau).$  Substituting
$\xi=\phi$ into the first equation, we have
$\frac{1}{16}(9\ov\tau^2+9\tau^2-30|\tau|^2)+o(|\tau|^2)=0.$  Thus
the system can not hold for a generic $\tau$.  It follows that a
generic elliptic CR singular point of $V(\tau)$ stays at a CR point
of $M$ where no holomorphic curves contained in $M$ can pass.

\medskip
Case (E) : In this case, $M$ is defined by $w=z_1\ov z_2+z_2\ov
z_1+2b(z_1z_2+\ov z_1\ov z_2)+\frac{1}{2}(z_2^2+\ov
z_2^2)+O(|z|^3).$ Then $\vec{c}=(1,-4b)$ is an elliptic direction
and we construct the family $V(\tau)$  by setting
$z_1=\xi,z_2=-4b\xi+\tau$. Then $\mathcal S$ is defined by the
following system of equations:
\begin{equation*}\left\{
\begin{matrix}
0=\frac{\partial w}{\partial \ov z_1}=&z_2+2b\ov z_2+o(|z|)\\
0=\frac{\partial w}{\partial \ov z_2}=&z_1+2b\ov z_1+\ov z_2+o(|z|)\\
0=\frac{\partial\ov  w}{\partial  z_1}=&\ov z_2+2bz_2+o(|z|)\\
0=\frac{\partial \ov w}{\partial  z_2}=&\ov z_1+ 2bz_1+ z_2+o(|z|)\\
\end{matrix}\right.\quad\quad.
\end{equation*}
Thus the Jacobian matrix at the origin is given by
\begin{equation*}\left(
\begin{matrix}
0&0&1&2b\\
1&2b&0&1\\
0&0&2b&1\\
2b&1&1&0\\
\end{matrix}\right),
\end{equation*}
of which the rank is at least three.  Therefore, by the implicit
function theorem, $\mathcal S$ has real dimension at most one.

Similarly, at the points where there is a holomorphic curve passing,
%it holds that $[L,\ov L]\subset{\rm Span}\{L,\ov L\}.$
%Therefore,
we have
\begin{equation*}\begin{array}{lll}
G_2=\ov z_1+z_2+2bz_1+o(|z|),\quad G_1=\ov z_2+2bz_2+o(|z|),
G_1G_{\ov 2}=-G_2G_{\ov 1}.
%G_2,\quad\lambda_{(2)}=G_1,\quad\lambda_{(4)}=-G_{\ov 2},\quad\lambda_{(5)}=-G_{\ov 1},\quad\lambda_{(3)}=o(|z|),\quad\lambda_{(6)}=o(|z|),\\
\end{array}\end{equation*}
%and
%$$G_2\frac{\partial}{\partial\ov z_1}+G_1\frac{\partial}{\partial\ov z_2}-G_{\ov 2}\frac{\partial}{\partial z_1}-G_{\ov 1}\frac{\partial}{\partial z_2}=\alpha(G_{ 2}\frac{\partial}{\partial z_1}-G_{1}\frac{\partial}{\partial z_2})+\beta(G_{\ov 2}\frac{\partial}{\partial\ov z_1}-G_{\ov 1}\frac{\partial}{\partial\ov z_2}).$$
%In particular,
%$G_1G_{\ov 2}=-G_2G_{\ov 1}$.

At the intersection of  the elliptic points of the Bishop surfaces
with smooth points of a holomorphic curve, we have
%$V(\tau)\cap M$ is given by the equation $\frac{\partial
%w}{\partial\ov \xi}=0.$  We substitute $z_1=\xi,z_2=-4b\xi+\tau$ in
%$\frac{\partial w}{\partial\ov \xi}=0$ and $G_1G_{\ov 2}=-G_2G_{\ov
%1}$ to get
\begin{equation*}\left\{
\begin{matrix}
&\Re((-4b\ov\xi+\tau-8b^2\xi+2b\tau)(\xi-2b\ov\xi+\ov\tau))=o({|\xi^2|+|\tau|^2})\\
&-8b\xi+\tau-2b\ov\tau+o(\sqrt{|\xi^2|+|\tau|^2})=0\\
&-8b\ov\xi+\ov\tau-2b\tau+o(\sqrt{|\xi^2|+|\tau|^2})=0\\
\end{matrix}\right.\quad
\end{equation*}
From the last two equations, we can solve $\xi=\phi(\tau,\ov \tau)$
with $\phi(\tau,\ov
\tau)=\frac{1}{8b}\tau+\frac{-1}{4}\ov\tau+o(\tau).$  Substituting
$\xi=\phi$ into the first equation, we have
$\frac{1}{4}\{(4b^3+4b^2+b+2)(\tau^2+\ov\tau^2)+(4b^2+5b+\frac{1}{4b})|\tau|^2\}+o(|\tau|^2)=0.$
Thus the system can not hold for a generic $\tau.$  It follows that
a generic elliptic CR singular point of $V(\tau)$ stays at a CR
point of $M$ where no holomorphic curves contained in $M$ can pass.

\medskip
Case (F): We now have $w=z_1\ov z_2+z_2\ov z_1+\frac{1}{2}(z_1^2+\ov
z_1^2)+(dz_2^2+\ov d \ov z_2^2)+E(z,\ov z).$ Then $\vec{c}=(1,C)$ is
an elliptic direction, where $C$ is a complex number such that
$\frac{1}{2}+dC^2=0$.  We construct the family $V(\tau)$  by letting
$z_1=\xi,z_2=C\xi+\tau$.  Then $\mathcal S$ is defined by the
following system of equations:
\begin{equation*}\left\{
\begin{matrix}
0=\frac{\partial w}{\partial \ov z_1}=&z_2+\ov z_1+o(|z|)\\
0=\frac{\partial w}{\partial \ov z_2}=&z_1+2\ov d\ov z_2+o(|z|)\\
0=\frac{\partial\ov  w}{\partial  z_1}=&z_1+\ov z_2+o(|z|)\\
0=\frac{\partial \ov w}{\partial  z_2}=&\ov z_1+ 2dz_2+o(|z|)\\
\end{matrix}\right.\quad\quad.
\end{equation*}
Thus the Jacobian matrix at the origin is given by
\begin{equation*}\left(
\begin{matrix}
0&1&1&0\\
1&0&0&2\ov d\\
1&0&0&1\\
0&1&2d&0\\
\end{matrix}\right),
\end{equation*}
which is invertible due to the fact that ${\rm Im}d>0$.  Therefore
by the implicit function theorem, $\mathcal S$ consists of a single
point near 0.

Similar to the previous cases, we have
\begin{equation*}\begin{array}{lll}
A=G_2=\ov z_1+2dz_2+o(|z|),\quad B=G_1=z_1+\ov z_2+o(|z|),  \quad C=o(|z|), \\
\lambda_{(1)}=G_2,\quad\lambda_{(2)}=G_1,\quad\lambda_{(4)}=-G_{\ov 2},\quad\lambda_{(5)}=-G_{\ov 1},\quad\lambda_{(3)}=o(|z|),\quad\lambda_{(6)}=o(|z|),\\
\end{array}\end{equation*}
and
$$G_2\frac{\partial}{\partial\ov z_1}+G_1\frac{\partial}{\partial\ov z_2}-G_{\ov 2}\frac{\partial}{\partial z_1}-G_{\ov 1}\frac{\partial}{\partial z_2}=\alpha(G_{ 2}\frac{\partial}{\partial z_1}-G_{1}\frac{\partial}{\partial z_2})+\beta(G_{\ov 2}\frac{\partial}{\partial\ov z_1}-G_{\ov 1}\frac{\partial}{\partial\ov z_2}).$$
In particular, $G_1G_{\ov 2}=-G_2G_{\ov 1}$.

On the other hand, the elliptic points of the Bishop surfaces
$V(\tau)$ is given by the equation $\frac{\partial w}{\partial\ov
\xi}=0.$  We substitute $z_1=\xi,z_2=C\xi+\tau$ into the equations:
$\frac{\partial w}{\partial\ov \xi}=0$ and $G_1G_{\ov 2}=-G_2G_{\ov
1}$ to get
\begin{equation*}\left\{
\begin{matrix}
&\Re(\xi^2+\xi(\ov C\ov\xi+\ov\tau)+2\ov d\xi(\ov C\ov\xi+\ov\tau)+2\ov d(\ov C\ov\xi+\ov\tau)^2)=o({|\xi^2|+|\tau|^2})\\
&(C+\ov C)\ov\xi+\ov\tau+2 dC\tau+o(\sqrt{|\xi^2|+|\tau|^2})=0\\
&(C+\ov C)\xi+\tau+2\ov d\ov C\ov\tau+o(\sqrt{|\xi^2|+|\tau|^2})=0\\
\end{matrix}\right.\quad
\end{equation*}
From the last two equations, we can solve  $\xi=\phi(\tau,\ov \tau)$
with $\phi(\tau,\ov \tau)=\frac{-1}{C+\ov C}\tau+\frac{-2\ov d\ov
C}{C+\ov C}\ov\tau+o(\tau).$  Substituting $\xi=\phi$ in the first
equation, we derive the coefficient of  $\tau^2$ to be
$\frac{1-2d}{(C+\ov C)^2}(1+2d|C|^2),$ which shows that the above
system can not hold for a generic $\tau$.  Hence a generic elliptic
CR singular point of $V(\tau)\cap M$ stays at a CR point of $M$
where no holomorphic curves contained in $M$ can pass.

Combining what we did above, we complete the proof of Lemma
\ref{132}. Thus we completed the proof of Theorem \ref{005} except
when $n+1=3$ and the quadratic normal form of $M$ is as in
(\ref{M1}).\ $\endpf$

\section{ Flattening parabolic CR singular points: a formal argument approach}
\subsection {A reduction to a formal argument}
We now proceed to the proof of the Theorem \ref{005} when
$(M,p)\subset \mathbb C^3$ has both  Bishop invariants parabolic at
$p.$ Namely, we assume that, after a holomorphic change of
variables, $p$ is mapped to $0$ and $M$ is defined by an equation of
the form:
\begin{equation}\label{def}
w=F(z,\ov z)=|z_1|^2+|z_2|^2+\frac{1}{2}(z_1^2+z_2^2+\ov z_1^2+\ov
z_2^2)+O(|z|^3).
\end{equation}
In this case, $M$ still has an elliptic direction $\vec{c}=(1,i)$.
And the same proof as in the previous section shows that $(M,0)$ can
be holomorphically flattened if the set of CR singular points of $M$
has real dimension at most $1$.

 However,  the set  $S$ of the CR singular points of $M$ might have real
dimension two. For instance, if $M$ is defined by
$w=q^{(2)}(z,\ov{z})(1+O(1))$ with
$q^{(2)}=|z_1|^2+|z_2|^2+\frac{1}{2}(z_1^2+z_2^2+\ov z_1^2+\ov
z_2^2)$, then the set of the CR singular points near the origin is
defined by $\Re{z_1}=\Re{z_2}=0$, which is of real dimension two. If
this happens, when we consider
$M_{a,c,\tau}=M\cap\{(z_1,z_2)=(1,i)\xi+\vec{a} \tau \}$, regarded
as a surface in the $(\xi,w)$ plane, for any vector
$\vec{a}=(a_1,a_2)$ not proportional to $\vec{c}$ and for any
$\tau\approx 0,$ the elliptic CR singular points $p_{c,a}(\tau)$ of
$M_{a,c,\tau}$ will also be the CR singular points of $M$ as a
4-manifold  in $\mathbb C^3.$ Thus the construction of the Levi-flat
hypersurface $H_{p_{c,a}(\tau)}$ in the last section will not apply
here. Indeed, one would suspect that $(M,0)$ might be not
flattenable as suggested in Remark \ref{22.221}.
%Indeed, we need to employ a
%different method to prove the flatness of $M$ at $p=0.$
Fortunately, the quadratic term in the defining equation of $M$ in
(\ref{def})  now is in the simplest symmetric form, which made the
formal normal form theory developed in [HY4] disposable here. And we
are still  able to produce a positive flattening result. (Other
related work related to formal arguments in CR analysis and geometry
can be found, for instance, in [CM] [LM1-2] [KZ].)

We will follow the strategy employed in [HY4]. However, since the
formal argument in [HY4] has to exclude the case when both
generalized Bishop invariants are parabolic, we need  some new ideas
to handle the current situation. In what follows, we will be very
brief for those arguments already contained in [HY4].

%. we now, describe the proof in this setting in detail.
We first choose $\vec{c}=(1,i)$ and $\vec{a}=(0,1)$. As in [HY4],
let $\wh M_{p_{a,c}(\tau)}$ be the Levi-flat submanifold bounded by
$M_{a,c,\tau}$  as constructed in Huang-Krantz [HK].   Let $\wh
M_{a,c}=\cup_{||\tau||<<1}\wh M_{p_{c,a}(\tau)}.$ Then $\wh M_{a,c}$
is a real analytic hypersurface in $\mathbb C^3$ with $M$ near $0$
as part of its real analytic boundary [HY3].  Now, suppose we can
flatten $(M,0) $ to order $N$.  Then as remarked in $(6.1)$ in
[HY4], the Levi-form of $\wh M_{a,c}$ vanishes to order at least
$\frac{N}{2}-3$ at $0.$  Suppose for $N^{\prime}>N,$ with a
holomorphic change of coordinates of the form $$\Phi:z^{\prime}=z,
w^{\prime}=w+O(|(z,w)|)^2,$$ $M$ can be further flattened to order
$N^{\prime}$, namely, in the new coordinates $M$ near $0$ is now
defined by an equation of the form as in (\ref{def}) with
$\Im{F}=O(N^{\prime}+1).$ Notice that with such a transform,
$M_{p_{a,c}(\tau)}$ is still mapped to
$M^{\prime}_{p_{a,c}(\tau)}=\Phi(M)\cap\{z^{\prime}=c\xi+a\tau\}$
and $p^{\prime}_{a,c}(\tau)=\Phi(p_{a,c}(\tau)).$ $\Phi(\wh
M_{p_{a,c}(\tau)})$ is now a Levi flat submanifold foliated by
attached holomorphic disks shrinking down to
$p^{\prime}_{a,c}(\tau).$ By the unique result of Kenig-Webster
[KW], $\Phi((\wh M_{a,c},0))=(\wh {M^{\prime}_{a,c}},0)$. Thus the
Levi form of $\wh {M^{\prime}_{a,c}}$ has vanishing order at least
$\frac{N^{\prime}}{2}-3$.  Since the vanishing order of the Levi
form is a holomorphic invariant, by the analyticity of $\wh
{M^{\prime}}_{p_{a,c}(\tau)}$, we see that $\wh M_{a,c}$ is in fact
Levi flat itself.  Thus we will complete the proof of  the
flattening of $M$ near $0$ if we can verify  the following:
\begin{theorem}\label{formal}Let M be a smooth 4-manifold in $\mathbb C^3$ defined by (\ref{def}).  Assume all CR points of $M$ are non-minimal.  Then for any $N$, there is a holomorphic change of coordinates of the form:
    $$z^{\prime}=z, w^{\prime}=w+B(z,w)=w+O(|(z,w)|)^2,$$
which flattens $(M,0)$ to order $N$
\end{theorem}
For the proof of Theorem \ref{formal}, we basically follow the
approach  developed in Huang-Yin [HY4].  However,  new ideas are
needed as the argument in [HY4] had to exclude the case when all
generalized Bishop invariants are  parabolic.
\subsection{Proof of Theorem \ref{formal} }
Before proceeding to the proof of Theorem \ref{formal}, we recall
results from [HY4] which will be needed here. We let $H$ be a
real-valued  homogeneous polynomial in $(z,\ov z)$ of order $m.$
Define
\begin{equation}\label{defi1}
\begin{split}
&\Phi=(z_2+\ov z_2)H_{\ov{1}}-(z_1+\ov z_1)H_{\ov{2}} {\,\,\,\rm {and}}\\
&\Psi=(z_2+\ov z_2)^2 \Phi_1-(z_2+\ov z_2)(z_1+\ov
z_1)\Phi_2+(z_1+\ov z_1)\cdot \Phi,\
\\
\end{split}
\end{equation}
where for a function $\chi(z,\ov z),$ as before, we write
$\chi_{\alpha}=\frac{\partial\chi}{\partial z_{\alpha}},\chi_{\ov
\alpha}=\frac{\partial\chi}{\partial\ov z_{\alpha}}$. Consider the
linear partial differential equation in $H$: (See [Appendix 7, HY4])
\begin{equation}\label{fundamental}
(z_2+\ov z_2)\Psi_1-(z_1+\ov z_1)\Psi_2=0.\\
\end{equation}
We follow [HY4] to introduce the following notation: For any
homogeneous polynomial $\chi$ of degree $N$ we write
$$\chi=\sum_{t+s+r+h=N}\chi_{[tsrh]}
z_1^{s}z_2^t\ov z_1^{h} \ov z_2^{r},$$ where by convention,
$\chi_{[tsrh]}=0$ if one of the indices is a negative integer.
Notice that in (\ref{defi1}) $H,\Phi$ and $\Psi$ are homogeneous
polynomials of degree $m,m$ and $m+1$ respectively.

 From (\ref{defi1}), we get:
 \begin{equation}\label{A1}
 \begin{split}
 \Phi_{[tsrh]}=&(h+1)H_{[ts(r-1)(h+1)]}+(h+1)H_{[(t-1)sr(h+1)]}\\
 &-(r+1)H_{[t(s-1)(r+1)h]}-(r+1)H_{[ts(r+1)(h-1)]},
 \end{split}
 \end{equation}
 and
 \begin{equation}\label{A2}
 \begin{split}
 \Psi_{[tsrh]}=&(s+1)\big\{\Phi_{[t(s+1)(r-2)h]}+2\Phi_{[(t-1)(s+1)(r-1)h]}
 +\Phi_{[(t-2)(s+1)rh]}\big\}\\
 &-(t+1)\Phi_{[(t+1)s(r-1)(h-1)]}
 -t\Phi_{[tsr(h-1)]}-(t+1)\Phi_{[(t+1)(s-1)(r-1)h]}\\
 & - t\Phi_{[t(s-1)rh]}+\Phi_{[tsr(h-1)]}+\Phi_{[t(s-1)rh]}.
 \end{split}
 \end{equation}

 Collecting the coefficients of
 $z_2^tz_1^{s}\ov{z_2}^{r}\ov{z_1}^h$ for $t\geq 0$, $s\geq 0$,
 $r\geq 0$ and $h=m+1-t-s-r\geq 0$ in (\ref{fundamental}), we get
 \begin{equation}\label{iden}
 (s+1)\Psi_{[(t-1)(s+1)rh]}+(s+1)\Psi_{[t(s+1)(r-1)h]}-(t+1)\Psi_{[(t+1)(s-1)rh]}-(t+1)\Psi_{[(t+1)sr(h-1)]}=0,
 \end{equation} or
 \begin{equation*}
 \begin{split}
 (s+1)\big\{& \Psi_{[(t-1)(s+1)rh]}+\Psi_{[t(s+1)(r-1)h]}\big\}=
 \mathcal{F}\{(\Psi_{[t's'r'h']})_{s'+h'\leq s+h-1,s'\leq s,h'\leq
    h}\}.
 \end{split}
 \end{equation*}
 Here, for a set of complex numbers (or polynomials) $\{a_j,
 b\}_{j=1}^{k}$, we say $b\in {\cal F}\{a_1,\cdots, a_k\}$ if
 $b=\sum_{j=1}^{k}(c_ja_j+d_j\ov{a_j})$ with $c_j,d_j\in
 {\mathbb{C}}$.
Now, we start in (\ref{iden}) with $r=0, t=m+1-(s+h)$ and notice
that $\Psi_{[(m-s-h)(s+1)(-1)h]}=0$. Letting $r=1,2,\cdots $, we
inductively get
 \begin{equation}\label{A3}
 \Psi_{[tsrh]}=\mathcal F\{(\Psi_{[t's'r'h']})_{s'+h'\leq s+h-2,s'\leq s,h'\leq
    h}\}, {\rm\,\, for\,\, }  s\geq1, t,r,h\geq0  {\,\,\rm and\,\, } t+s+r+h=m+1.
 \end{equation}

 Making use of (\ref{A2}) and (\ref{A3}), we get, for $s\geq 1,r,t,h\geq 0,r+s+t+h=m+1$, the following:
 $$(s+1)\big\{\Phi_{[t(s+1)(r-2)h]}+2\Phi_{[(t-1)(s+1)(r-1)h]}
 +\Phi_{[(t-2)(s+1)rh]}\big\}=\mathcal F\{(\Phi_{[t's'r'h']})_{s'+h'\leq s+h-1,s'\leq s+1,h'\leq
    h}\},$$
 for $s\geq1,$ and $ t,r,h\geq0.$   As above, we then get inductively:
 \begin{equation}\label{A4}
 \Phi_{[tsrh]}=\mathcal F\{(\Phi_{[t's'r'h']})_{s'+h'\leq s+h-2,s'\leq s,h'\leq
    h}\}, {\rm\,\, for\,\, } s\geq2, t,r,h\geq0  {\,\,\rm and\,\, } t+s+r+h=m.
 \end{equation}
 Further substituting (\ref{A1}) into (\ref{A4}), we get, for $s\geq2, h\geq0,$ the following:
 $$(h+1)H_{[ts(r-1)(h+1)]}+(h+1)H_{[(t-1)sr(h+1)]}=\mathcal F\{(H_{[t's'r'h']})_{s'+h'\leq s+h-1,s'\leq s,h'\leq
    h+1}\}.$$
 Hence, by an induction argument, we have for $s\geq2,h\geq1:$
 $$H_{[ts(m-t-s-h)h]}=\mathcal F\{(H_{[t's'(m-t'-s'-h')h']})_{s'+h'\leq s+h-2,s'\leq s,h'\leq
    h}\}.$$
 Applying this property inductively  and noticing that $H_{[tsrh]}=\ov{H_{[rhts]}}$, we get
 \begin{equation}\label{A5}
 H_{[ts(m-t-s-h)h]}=\mathcal F\{(H_{[t'1(m-2-t')1]}),(H_{[t'0(m-t'-i)i]})_{i\leq {\rm max\,}(s,h)}\}.
 \end{equation}

A crucial step for the proof of Theorem \ref{formal} is  to prove
the following uniqueness result:
\begin{proposition}\label{unique}
    Assume $\{H_{[tsrh]}\}$ satisfies the normalization condition in Definition \ref{norm}  and $H$ satisfies the equations (\ref{defi1}) and (\ref{fundamental}).  Suppose $H_{[t1(m-t-2)1]}=H_{[t0(m-t)0]}\equiv0.$  Then $H\equiv0.$
\end{proposition}

Before giving a proof of Proposition \ref{unique}, we recall a
special case of an initial normalization condition introduced in
[HY4]. (See [Theorem 4.2, HY4])

\begin{definition}\label{norm}
For a homogeneous real-valued polynomial  $E(z,\-{z})$ of degree
$m_0\ge 3$ $(z=(z_1,z_2))$, we say that it satisfies the initial
normalization condition if we have
\begin{equation}\label{ng}
  E_{(s_1e_1+s_2e_2,0)}=E_{(te_2+t_1e_1+t_2e_2,se_2)}=0\ \text{for}\ t\geq s,\ s_1+s_2=m_0,\ t_1+t_2>0,
  t+t_1+t_2+s=m_0,
\end{equation}
and we have the following condition  in terms of the remainder of
$m_0$ when divided by $6$:
%\begin{enumerate}
%%%%%%%%%%%%%%%%%%%%%%%%%%%%%% -3 %%%%%%%%%%%
% \item[($II_{-3}$)]

\noindent (II$_{-3}$)
 If $m_0=6\hat{m}-3$, then
\begin{equation}\begin{split}\label{n-3}
&E_{(te_2,se_2)}=0\ \text{for}\  4\hat{m}-1\leq t\leq m_0-1,\\
&E_{((2t+1)e_2+e_1,(m_0-2t-3)e_2+e_1)}=0\ \text{for}\ 2\hat{m}-2\leq
t\leq 3\hat{m}-3.
\end{split}\end{equation}

%%%%%%%%%%%%%%%%%%%%%%%%%%%%%% -2 %%%%%%%%%%%
%\item[($II_{-2}$)]
(II$_{-2}$) If $m_0=6\hat{m}-2$, then
\begin{equation}\begin{split}\label{n-2}
&E_{(te_2,se_2)}=0\ \text{for}\  4\hat{m}-1\leq t\leq m_0-1,\\
&E_{((2t+1)e_2+e_1,(m_0-2t-3)e_2+e_1)}=0\ \text{for}\ 2\hat{m}-1\leq
t\leq 3\hat{m}-3,\\
&{{\Re}}E_{((4\hat{m}-3)e_2+e_1,(2\hat{m}-1)e_2+e_1)}=0.
\end{split}\end{equation}

%%%%%%%%%%%%%%%%%%%%%%%%%%%%%% -1 %%%%%%%%%%%
% \item[($II_{-1}$)]
(II$_{-1}$)
 If $m_0=6\hat{m}-1$, then
\begin{equation}\begin{split}\label{n-1}
&E_{(te_2,se_2)}=0\ \text{for}\  4\hat{m}\leq t\leq m_0-1,\\
&E_{((2t+1)e_2+e_1,(m_0-2t-3)e_2+e_1)}=0\ \text{for}\ 2\hat{m}-1\leq
t\leq 3\hat{m}-2.
\end{split}\end{equation}

%%%%%%%%%%%%%%%%%%%%%%%%%%%%%% 0 %%%%%%%%%%%
%\item[($II_{0}$)]
(II$_{0}$) If $m_0=6\hat{m}$, then
\begin{equation}\begin{split}\label{n0}
&E_{(te_2,se_2)}=0\ \text{for}\  4\hat{m}+1\leq t\leq m_0-1,\\
&E_{((2t+1)e_2+e_1,(m_0-2t-3)e_2+e_1)}=0\ \text{for}\ 2\hat{m}-1\leq
t\leq 3\hat{m}-2,\\
&{{\Re}}E_{(4\hat{m}e_2,2\hat{m}e_2)}=0.
\end{split}\end{equation}

%%%%%%%%%%%%%%%%%%%%%%%%%%%%%% 1 %%%%%%%%%%%
% \item[($II_{1}$)]
(II$_{1}$)
 If $m_0=6\hat{m}+1$, then
\begin{equation}\begin{split}\label{n1}
&E_{(te_2,se_2)}=0\ \text{for}\  4\hat{m}+1\leq t\leq m_0-1,\\
&E_{((2t+1)e_2+e_1,(m_0-2t-3)e_2+e_1)}=0\ \text{for}\  2\hat{m}\leq
t\leq 3\hat{m}-1.
\end{split}\end{equation}

%%%%%%%%%%%%%%%%%%%%%%%%%%%%%% 2 %%%%%%%%%%%
%\item[($II_{2}$)]
(II$_{2}$) If $m_0=6\hat{m}+2$, then
\begin{equation}\begin{split}\label{n2}
&E_{(te_2,se_2)}=0\ \text{for}\  4\hat{m}+2\leq t\leq m_0-1,\\
&E_{((2t+1)e_2+e_1,(m_0-2t-3)e_2+e_1)}=0\ \text{for}\  2\hat{m}\leq
t\leq 3\hat{m}-1,\\
&{{\Re}}E_{((4\hat{m}+1)e_2,(2\hat{m}+1)e_2)}=0.
\end{split}\end{equation}

%\end{enumerate}

\end{definition}
In the above, we write
$E_{(s_1e_1+s_2e_2,t_1e_1+t_2e_2)}=E_{[s_2s_1t_2t_1]}.$ By an
induction argument,  to conclude Proposition \ref{unique}, it
suffices to prove the following:
\begin{lemma}\label{induction} Assume that $\{H_{[tsrh]}\}$ satisfies the normalization condition as in Definition \ref{norm} and satisfies the condition that
 $H_{[t1(m-t-2)1]}=0$ and  $H_{[t0(m-t)0]}=0.$  Assume that there exists an $h_0\geq
-1$  such that
\begin{equation}\label{11a}
\Psi_{[tsrh]}=\Phi_{[tsrh]}=0\ \text{for}\ h\leq h_0\ \hbox{ and any
} s,\hbox{ and } H_{[tsrh]}=0\ \text{for}\ \max(s,h)\leq h_0+1.
\end{equation}
Then we have
\begin{equation}\label{11c}
\Psi_{[tsrh]}=\Phi_{[tsrh]}=0\ \text{for}\ h\leq h_0+1 \hbox{ and
any } s,\ \hbox{ and } H_{[tsrh]}=0\ \text{for }\ \max(s,h)\leq
h_0+2.
\end{equation}
\end{lemma}

Before proceeding to prove Lemma \ref{induction}, we first introduce
some notation:
\begin{equation}\label{psiphik}
\begin{split}
&\Psi^{(k)}_{[sh]}={\sum}_{t=k}^{m+1-s-h}(-1)^{m+1-t-s-h}(_k^t)\Psi_{[ts(m+1-t-s-h)h]},\\
&\Phi^{(k)}_{[sh]}={\sum}_{t=k}^{m-s-h}(-1)^{m-t-s-h}(_k^t)\Phi_{[ts(m-t-s-h)h]},\\
&H^{(k)}_{[sh]}={\sum}_{t=k}^{m-s-h}(-1)^{m-t-s-h}(_k^t)H_{[ts(m-t-s-h)h]}.
\end{split}
\end{equation}
In what follows, we make the convention that for two integers $k$
and $ t$, we set $(_k^t)=0$ if $t<k$ or if $k\cdot t<0$ and we set
$(_0^t)=1$ for $t\ge 0$. We also recall the convention that
$H_{[tsrh]}=\Phi_{[tsrh]}=\Psi_{[tsrh]}=0$ when one of the indices
is negative. Then the following can also be written as:
\begin{equation}\label{psiphik1}
\begin{split}
&\Psi^{(k)}_{[sh]}={\sum}_{t=-\infty}^{\infty}(-1)^{m+1-t-s-h}(_k^t)\Psi_{[ts(m+1-t-s-h)h]},\\
&\Phi^{(k)}_{[sh]}={\sum}_{t=-\infty}^{\infty}(-1)^{m-t-s-h}(_k^t)\Phi_{[ts(m-t-s-h)h]},\\
&H^{(k)}_{[sh]}={\sum}_{t=-\infty}^{\infty}(-1)^{m-t-s-h}(_k^t)H_{[ts(m-t-s-h)h]}.
\end{split}
\end{equation}

 As in [HY4], we would like
first to transfer the relations among $\Psi$, $\Phi$ and $H$ into
the relations among $\Psi^{(k)}_{[s(h_0+1)]}$,
$\Phi^{(k)}_{[s(h_0+1)]}$ and $H^{(k)}_{[s(h_0+2)]}$.
\begin{lemma}\label{formular}
    Assume that there exists an $h_0\geq
    -1$ such that
    \begin{equation}\label{11aa}
    \Psi_{[tsrh]}=\Phi_{[tsrh]}=0\ \text{for}\ h\leq h_0,\ H_{[tsrh]}=0\
    \text{for}\ \max(s,h)\leq h_0.
    \end{equation}
    Then from (\ref{A1}) and (\ref{A2}), we have for any $k\ge 0$ the following
    \begin{align}
    \Phi^{(k)}_{[s(h_0+1)]}=&(h_0+2)H^{(k-1)}_{[s(h_0+2)]} +(m-s-h_0-k)H^{(k)}_{[(s-1)(h_0+1)]}
    -(k+1)H^{(k+1)}_{[(s-1)(h_0+1)]}, \label{B1}\\
    \Psi^{(k)}_{[s(h_0+1)]}=&(s+1)\Phi^{(k-2)}_{[(s+1)(h_0+1)]}-(k-1)\Phi^{(k)}_{[(s-1)(h_0+1)]}.  \label{B2}
    \end{align}
    Moreover, by (\ref{fundamental}), $\Psi^{(k)}_{[s(h_0+1)]}$ satisfies the following
    equation:
    \begin{equation}\label{B3}
    (s+1)\Psi^{(k-1)}_{[(s+1)(h_0+1)]}=
    (k+1)\Psi^{(k+1)}_{[(s-1)(h_0+1)]}.
    \end{equation}
\end{lemma}
{\it Proof of the Lemma \ref{formular}:} The proofs for  (\ref{B1})
and (\ref{B2}) are the same as that in [Lemma 5.4, HY4].  For the
last one, we modify the method used in [HY4] as follows: Combining
(\ref{iden}) and the assumption that $\Psi_{[tsrh_0]}=0$, we have
$$(s+1)\Psi_{[(t-1)(s+1)r(h_0+1)]}+(s+1)\Psi_{[t(s+1)(r-1)(h_0+1)]}-(t+1)\Psi_{[(t+1)(s-1)r(h_0+1)]}=0.$$
Then using the convention set up before, we obtain the following:
\begin{equation*}
\begin{split}
0&=\sum_{t=-\infty}^{\infty}(-1)^r(_k^t)\{(s+1)\Psi_{[(t-1)(s+1)r(h_0+1)]}+(s+1)\Psi_{[t(s+1)(r-1)(h_0+1)]}-(t+1)\Psi_{[(t+1)(s-1)r(h_0+1)]}\}\\
&=(s+1)\sum_{t=-\infty}^{\infty}\{(_k^{t-1})+(_{k-1}^{t-1})\}(-1)^r\Psi_{[(t-1)(s+1)r(h_0+1)]}+(s+1)\sum_{t=-\infty}^{\infty}(-1)^r(_k^t)\Psi_{[t(s+1)(r-1)(h_0+1)]}\\
&\,\,\,\,\,\,\,\,\,\,\,\,-(k+1)\sum_{t=-\infty}^{\infty}(-1)^r(_{k+1}^{t+1})\Psi_{[(t+1)(s-1)r(h_0+1)]}\\
&=(s+1)\Psi^{(k)}_{[(s+1)(h_0+1)]}+(s+1)\Psi^{(k-1)}_{[(s+1)(h_0+1)]}-(s+1)\Psi^{(k)}_{[(s+1)(h_0+1)]}-(k+1)\Psi^{(k+1)}_{[(s-1)(h_0+1)]}.\\
&=(s+1)\Psi^{(k-1)}_{[(s+1)(h_0+1)]}-(k+1)\Psi^{(k+1)}_{[(s-1)(h_0+1)]}.
\quad\quad\quad\quad\quad\quad\endpf
\end{split}
\end{equation*}
%Now we proceed to prove Lemma \ref{induction}.
{\it Proof of the Lemma \ref{induction}:}  From (\ref{A3}) we have
\begin{equation*}
\begin{split}
\Psi_{[t(2s+1)r(h_0+1)]}=\mathcal{F}\{(\Psi_{[t'(2s'+1)r'(h_0+1)]})_{s'<s},
(\Psi_{[t''s''r''h'']})_{h''\leq h_0} \}\ \text{for}\ s\geq
0,h_0\geq-1, \text{and } t,r\geq 0.
\end{split}
\end{equation*}
In particular, we obtain
$$
\Psi_{[t1r(h_0+1)]}=\mathcal{F}\{(\Psi_{[t's'r'h']})_{h'\leq h_0}
\}=0.
$$
The last equality follows from the assumptions in (\ref{11a}). By an
induction argument, we obtain
\begin{equation}
\Psi_{[t(2s+1)r(h_0+1)]}=0, \ \text{for } s\geq0, h_0\geq-1.
\end{equation}
Combining this with (\ref{B2}), we have
$$(2s+2)\Phi^{(k-2)}_{[(2s+2)(h_0+1)]}-(k-1)\Phi^{(k)}_{[(2s)(h_0+1)]}=0, \text{ for } s\geq 0 \text{ and } h_0\geq-1.$$
Taking $k=0$ in the above equality, by the fact that
$\Phi^{(-2)}_{[(2s+2)(h_0+1)]}=0$, we have
$\Phi^{(0)}_{[(2s)(h_0+1)]}=0.$ Inductively, by setting
$k=2,4,\cdots$, we finally have $\Phi^{(2k)}_{[0(h_0+1)]}=0$ for all
$k\geq 0$.

Combing this with (\ref{B1}), we have $H_{[0(h_0+2)]}^{(2k-1)}=0$.
Now by the  argument from [(5.75), HY4] to [(5.77), HY4], we know
$H_{[t0r(h_0+2)]}=0$. (This is the place part of the initial
normalization conditions are used.) By the assumption $H_{[t1r1]}=0$
and the relation in (\ref{A5}), we get
$H_{[ts(m-t-s-h_0-2)(h_0+2)]}=0$ for $0\leq s\leq h_0+2.$ Then, by
(\ref{A1}), we have $\Phi_{[tsr(h_0+1)]}=0$ for $s\leq h_0+1.$  In
particular, $\Phi_{[t0r(h_0+1)]}=\Phi_{[t1r(h_0+1)]}=0.$ Thus by
(\ref{A4}), it follows that $\Phi_{[tsr(h_0+1)]}=0.$  Similarly by
(\ref{A2}) and (\ref{A3}), we  conclude
$\Psi_{[tsr(h_0+1)]}\equiv0.$  Thus we completed the proof of Lemma
\ref{induction}. $\endpf$

\medskip
 Now, let $M$ be as defined by (\ref{def}).  We assume $M$ is
flattened to order $m-1$.  We need to find a holomorphic change of
coordinates of the form
\begin{equation}\label{Bchange}
z'=z,w=w+B(z,w)=w+O(|z,w|^2).
\end{equation}to flatten $M$ to order $m$.  Write $H$ for the homogeneous polynomial of degree $m$ in the Taylor
expansion of $\Im{F}$ at $0.$ By [Appendix 7, HY4], $H$ satisfies
the linear equations in (\ref{defi1}) and (\ref{fundamental}).  By
 [Theorem
4.2, HY4], there is a unique transformation of the form
$z'=z,w'=w+B(z,w)$, where $B(z,w)=\sum_{|\a|+2j=m}b_{\a j }z^\a w^j$
with $b_{0\frac{m}{2}}=0$ if $m$ is even, which transforms $M$ to a
new manifold with $H=(\Im{F})^{(m)}$, the homogenous part of degree
$m$ in $\Im{F}$, that satisfies the normalization conditions in
Definition \ref{norm}.

{\bf Case I.  $m$ is odd}:  Step I: We first consider a holomorphic
change of coordinates: $(z_1,z_2)\rightarrow(-z_1,z_2).$  It is
clear that $H(z_1,z_2)$ and $H(-z_1,z_2)$ both satisfy
(\ref{defi1}),
 (\ref{fundamental}) and the initial normalization conditions in Definition \ref{norm}. Therefore  $H(z_1,z_2)-H(-z_1,z_2)$ is also a
solution to the linear system (\ref{defi1}) and (\ref{fundamental})
satisfying the initial normalizations in Definition \ref{norm}.
Write
$$H(z_1,z_2)=\sum_{s+h\text{ even}}H_{[tsrh]}z_1^sz_2^t\ov z_1^h\ov
z_2^r+\sum_{s+h\text{ odd}}H_{[tsrh]}z_1^sz_2^t\ov z_1^h\ov z_2^r.$$
$$\text{ Then    } \,\,\,\,H(-z_1,z_2)=\sum_{s+h\text{ even}}H_{[tsrh]}z_1^sz_2^t\ov z_1^h\ov z_2^r-\sum_{s+h\text{  odd}}H_{[tsrh]}z_1^sz_2^t\ov z_1^h\ov z_2^r$$
 $$H(z_1,z_2)-H(-z_1,z_2)=2\sum_{s+h\text{  odd}}H_{[tsrh]}z_1^sz_2^t\ov z_1^h\ov z_2^r$$

By Lemma \ref{unique},  $H(z_1,z_2)-H(-z_1,z_2)\equiv0$ or
equivalently, \begin{equation}\label{oddvanishing}
H(z_1,z_2)=\sum_{s+h\text{ even}}H_{[tsrh]}z_1^sz_2^t\ov z_1^h\ov
z_2^r.
\end{equation}

Step II: In this step, we assume the associated homogeneous
polynomial $H$ of the 4-manifold $M$ already take the form as in
(\ref{oddvanishing}).  We then holomorphically change the
coordinates by $(\wt z_1,\wt z_2):=(z_2,z_1).$  $\wt M$, the image
of $M$ under this map, is defined now by $$\wt w(\wt z_1,\wt
z_2)=\wt F(\wt z_1,\wt z_2)=F(\wt z_2,\wt z_1)=|\wt z_1|^2+|\wt
z_2|^2+\frac{1}{2}(\wt z_1^2+\wt z_2^2+\ov{ \wt z_1}^2+\ov {\wt
z_2}^2)+O(|\wt z|^3).$$ In particular, $\wt H$, the homogeneous
polynomial of degree $m$ in Tayler expansion of  $\Im{\wt F}$ at
$0,$ takes the following form, by the fact that $H_{[tsrh]}=0$ for
$t+r$ even:
$$\wt H(\wt z_1, \wt z_2)=\sum_{s+h\text{ odd}}\wt H_{[tsrh]}\wt z_1^s\wt z_2^t\ov {\wt z_1}^h\ov {\wt z_2}^r=\sum_{s+h\text{ odd}}H_{[sthr]}\wt z_1^s\wt z_2^t\ov {\wt z_1}^h\ov
{\wt z_2}^r.$$  Notice that $\wt H$ may not satisfy the
normalization conditions in Definition \ref{norm}.  Therefore, we
normalize $\wt M$ using a transformation of the following form:
\begin{equation}\label{renorm}
\begin{split}
&z'=\wt z\\
&w'=\wt w+\wt B(\wt z,\wt w)=\wt
w+\sum_{\alpha_1+\alpha_2+2k=m}a_{a_1a_2k}\wt z_1^{\alpha_1}\wt
z_2^{\alpha_2}\wt w^k,\ a_{0\frac{m}{2}}=0\ \hbox{if } m \
\hbox{even}.
\end{split}
\end{equation}
By Theorem 4.2 in [HY4], we know there is a unique $\wt B$ such that
the new $\wt H'$ satisfies the normalization conditions,
(\ref{defi1}) and (\ref{fundamental}).  Also from [HY4], we know
$\wt H'=\wt H+\text{Im} \wt B(\wt z, q^{(2)}(\wt z, \ov{\wt z}))$,
where $q^{(2)}(\wt z, \ov{\wt z})=|\wt z_1|^2+|\wt
z_2|^2+\frac{1}{2}(\wt z_1^2+\wt z_2^2+\ov{ \wt z_1}^2+\ov {\wt
z_2}^2).$

A crucial observation here is that, when $\alpha_1$, the power of
$\wt z_1$ in (\ref{renorm}), is odd, the expansion of $\wt
z_1^{\alpha_1}\wt z_2^{\alpha_2}\wt w^k$ after substituting $\wt
w=q^{(2)}(\wt z, \ov{\wt z})$ only contains terms $\wt z_1^s\wt
z_2^t\ov {\wt z_1}^h\ov {\wt z_2}^r$  with $s+h$ odd; and when
$\alpha_1$ is even, the expansion of $\wt z_1^{\alpha_1}\wt
z_2^{\alpha_2}\wt w^k$ after substituting $\wt w=q^{(2)}(\wt z,
\ov{\wt z})$ only contains terms $\wt z_1^s\wt z_2^t\ov {\wt
z_1}^h\ov {\wt z_2}^r$  with $s+h$ even. Then we split the
polynomial $B$ into two parts depending on whether $\alpha_1$ is odd
or even:
$$\wt B(\wt z,\wt w)=\wt B_1(\wt z,\wt w)+\wt B_2(\wt z,\wt w)=\sum_{\substack{  \alpha_1+\alpha_2+2k=m,\\ \alpha_1\text{ odd }}}a_{a_1a_2k}\wt z_1^{\alpha_1}\wt z_2^{\alpha_2}\wt w^k+\sum_{\substack{  \alpha_1+\alpha_2+2k=m,\\ \alpha_1\text{ even }}}a_{a_1a_2k}\wt z_1^{\alpha_1}\wt z_2^{\alpha_2}\wt w^k.$$
By the uniqueness of $\wt B$, we have $\wt B_2\equiv0$ and thus $\wt
H'_{[tsrh]}=\wt H_{[tsrh]}=0$ for $s+h$ even.  Now, by Lemma
\ref{norm} and  what we just argued, we have $\wt H'\equiv 0.$  Next
for the original 4-manifold $M$ with $H$ as in (\ref{oddvanishing}),
we define the transformation of the form:
\begin{equation}
\begin{split}
&z'= z\\
&w'= w+ B( z, w)= w+\sum_{\alpha_1+\alpha_2+2k=m}a_{a_1a_2k}
z_2^{\alpha_1}z_1^{\alpha_2} w^k,
\end{split}
\end{equation}
where $\{a_{a_1a_2k}\}$ are the same as in (\ref{renorm}), or
equivalently, $B(z_1,z_2,w)=\wt B(z_2,z_1,w)$.  Then after the
transformation,  $H'=H+\text{Im} B(z,q^{(2)}(z,\ov z))\equiv0.$ Thus
we prove Theorem \ref{formal} in the case when $m$ is odd.

{\bf Case II. $m$ is even: } In this case,  after establishing
[(5.47), HY4], the same argument there gives that
$H_{[t1r1]}=H_{[t0r0]}\equiv0$ for even $m$. Hence, by Lemma
\ref{unique}, we complete the argument for the case of $m$ even. Now
we proceed to establish the following identity:
\begin{equation}
H^{(2k-2)}_{[11]}+(m+1-2k)H^{(2k-1)}_{[00]}-2kH^{(2k)}_{[00]}=0.
\end{equation}
which is $(5.47)$ in [HY4] (with
$\theta=0$ and $\xi=1$).

By substituting $h_0=-1,k=2l$ in (\ref{B3}), we have
\begin{equation}\label{C1}
(s+1)\Psi^{(2l-1)}_{[(s+1)0]}= (2l+1)\Psi^{(2l+1)}_{[(s-1)0]}.
\end{equation}

Substituting $l=0,s=2k+1$ in (\ref{C1}), we have
$\Psi^{(1)}_{[(2k)0]}=(2k+2)\Psi^{(-1)}_{[(2k+2)0]}=0.$ By setting
$l=i, s=2k-2i+1$ for $i=1,\cdots,k$, inductively,  we will have
$\Psi^{(2k+1)}_{[00]}=\cdots=\Psi^{(2l+1)}_{[(2k-2l)0]}=\cdots=\Psi^{(1)}_{[(2k)0]}=0.$
Hence we proved that $\Psi^{(2k+1)}_{[00]}=0$ for all $k\geq0.$
Next, by substituting $s=0,h_0=-1, k=2l+1$ in (\ref{B2}), we have
$\Phi_{[10]}^{(2l-1)}=0$ for $l\geq 1.$ Finally, substituting
$s=1,h_0=-1,k=2l-1$ in (\ref{B1}), we have
\begin{equation}\label{lasteq}
H^{(2l-2)}_{[11]}+(m+1-2l)H^{(2l-1)}_{[00]}-2lH^{(2l)}_{[00]}=0,\text{
for } l\geq 1,
\end{equation}
which is exactly (5.47) in [HY4] by setting the $\theta$ and $\xi$
there to be zero and $1$, respectively. Therefore, we proved Theorem
\ref{formal} for $m$ even case.

We thus finally completed the proof of Theorem \ref{005}. $\endpf$
%\section{Proof of Corollary \ref{44.44} and further remarks}

\bigskip
We now complete the proof of Corollary \ref{44.44}. Let $(M,p)$ be
as in Corollary \ref{44.44} with $p$ now a definite non-degenerate
CR singular point. Let $\Phi$ be a biholomorphic map sending a
neighborhood of $p$ to a neighbor of $0$ with $\Phi(p)=0$ such that
$M'=\Phi(M)$ is defined by an equation of the form
$u=|z|^2+\sum_{j=1}^{n}\ld_j(z^2+\ov{z_j}^2)+O(3),\  v=0$. Now, let
$0<\epsilon_0<<1$. Then $\Phi(B_{\epsilon_0}(p))$ is a strongly
pseudoconvex domain containing the origin. Hence the holomophic hull
of $\Phi(B_{\epsilon_0}(p))\cap M'$ near $0$ coincides with  a
neighborhood of $0$ in the set defined  $u\ge
|z|^2+\sum_{j=1}^{n}\ld_j(z^2+\ov{z_j}^2)+O(3),\ v=0$ in the
Levi-flat hypersurface $v=0$. From this, one sees the proof of the
last statement in Corollary \ref{44.44}. Thus, the proof of
Corollary \ref{44.44} is complete.\ \ $\endpf$

\bigskip

\bigskip
\centerline{\bf REFERENCES}

\medskip

\noindent [BER]: S. Baouendi, P. Ebenfelt and L. Rothschild, Local
geometric properties of real submanifolds in complex space, {\it
Bull. Amer. Math. Soc.} (N.S.) 37, no. 3, 309--339, 2000.

\noindent [BR]: S. Baouendi and L. Rothschild,  Germs of CR maps
between real analytic hypersurfaces,  {\it Invent. Math. 93},
481-500, 1988.

\noindent [BG]: E. Bedford and B. Gaveau, Envelopes of holomorphy of
certain 2-spheres in ${\CC}^2$, {\it Amer. J. Math. 105}, 975-1009,
1983.

\noindent [BK]:  E. Bedford and W. Klingenberg,  On the envelopes of
holomorphy of a 2-sphere in ${\CC}^2$,  {\it Journal of AMS 4},
623-655, 1991.

\noindent [Bis]: E. Bishop, Differentiable manifolds in complex
Euclidean space,  {\it Duke Math. J. 32}, 1-21, 1965.

\noindent [BGr]:     T. Bloom and I. Graham,  On type  conditions
for generic real submanifolds of ${\mathbb C}^n$, {\it Invent.
Math.} Vol 40 (3),  217-243, 1977.

\noindent [Bur1]:  V. Burcea,
  A normal form for a real 2-codimensional submanifold in $M\subset {\bf C}^{N+1}$ near
 a CR singularity, {\it Adv. Math.} 243, 262-295, 2013.

\noindent [Bur2]:   V. Burcea, On a family of analytic discs
attached to a real submanifold $M$ in ${\bf C}^{N+1}$, {\it Methods
Appl. Anal.} 20, no. 1, 69-78, 2013.

%\noindent [Car]: \'E. Cartan, Sur les vari\'et\'es pseudo-conformal
%des hypersurfaces de l'espace de deux variables complexes, Ann. Mat.
%Pura Appl. (4) 11, 17-90, 1932.

\noindent [CM]: S. S. Chern and J. K. Moser, Real hypersurfaces in
complex manifolds, {\it Acta Math. 133,} 219-271, 1974.

\noindent [Co]: A. Coffman, CR singularities of real fourfolds in
$\mathbb C^3$, {\it Illinois J. Math.} 53, no. 3, 939-981, 2010.

\noindent [DTZ1]: P. Dolbeault, G. Tomassini and D. Zaitsev,
 On Levi-flat hypersurfaces with prescribed boundary, {\it Pure Appl.
Math. Q.} 6, no. 3, 725-753, 2010.

\noindent  [DTZ2]: P. Dolbeault, G. Tomassini and D. Zaitsev,
Boundary problem for Levi flat graphs, {\it Indiana Univ. Math. J.}
60, no. 1, 161-170, 2011.

\noindent [Gon1]:  X. Gong, Normal forms of real surfaces under
unimodular transformations near elliptic complex tangents, {\it Duke
Math. J. 74},  no. 1, 145--157, 1994.

\noindent [Gon2]: X.  Gong, Existence of real analytic surfaces with
hyperbolic complex tangent that are formally but not holomorphically
equivalent to quadrics, {\it Indiana Univ. Math. J. 53}, no. 1,
83--95, 2004.

\noindent [GL]: X. Gong and J. Lebl, Normal forms for CR singular
codimension two Levi-flat submanifolds, {\it Pacific J. Math.} 275, no. 1, 115-165, 2015.

\noindent [GS1]:  X. Gong and L. Stolovitch, Real submanifolds of
maximum complex tangent space at a CR singular point I,  {\it
Invent. Math.} 206, no. 1, 293-377, 2016.

\noindent [GS2]:  X. Gong and L. Stolovitch,  Real submanifolds of
maximum complex tangent space at a CR singular point II,  preprint.
(66 pages)

\noindent
 [Hu1]:  X. Huang, On  an n-manifold
 in ${\bf C}^n$ near an elliptic complex tangent,
{\it J. Amer. Math. Soc.} (11),  669--692,  1998.

\noindent [Hu2]: X. Huang, Geometric and analytic problems for a
real submanifold in  ${\mathbb C}^n$ with CR singularities, a survey
preprint dedicated to the memory of Qi-Keng Lu, to appear in {\it
Sciences in China}, 2017.

 \noindent [HK]: X. Huang and S. Krantz, On a problem of
Moser,
 {\it Duke Math. J.} Vol 78, 213-228, 1995.

\noindent [HY1]: X. Huang and  W. Yin, A Bishop surface with a
vanishing Bishop invariant,   {\it Invent. Math. } 176,  no.
3, 461--520, 2009.

\noindent [HY2]: X. Huang and  W. Yin, A codimension two CR singular
submanifold  that is formally equivalent to a symmetric quadric,
{\it Int. Math. Res. Not.  },  no. 15, 2789--2828,  2009.

\noindent [HY3]: X. Huang and W. Yin,  Flattening of CR singular
points and analyticity of the local hull of holomorphy I, {\it Math.
Ann.} 365, no. 1-2, 381-399, 2016.

\noindent [HY4]: X. Huang and W. Yin,  Flattening of  CR singular
points and analyticity of the local hull of holomorphy II, {\it Adv.
in Math.} 308,  1009-1073, 2017.

\noindent [KW1]:  C. Kenig and S. Webster, The local hull of
holomorphy of a surface in the space of two complex variables, {\it
Invent. Math. 67}, 1-21, 1982.

\noindent [KW2]: C. Kenig and S. Webster, On the hull of holomorphy
of an n-manifold in ${\bf C}^n$, {\it Annali Scoula Norm. Sup. de
Pisa IV, 11 No. 2}, 261-280, 1984.

\noindent [KZ]: I. Kossovskiy and D. Zaitsev, Convergent normal form
and canonical connection for hypersurfaces of finite type in
${\mathbb C}^2$, {\it  Adv. Math.} 281 (2015), 670-705.

\noindent [LM1]: B. Lamel and N. Mir, Parametrization of local CR
automorphisms by finite jets and applications, {\it J. Amer. Math.
Soc.} 20, no. 2, 519-572, 2007.

\noindent [LM2]: B. Lamel and N. Mir, Convergence of formal CR
mappings between  strongly pseudo-convex hypersurfaces, preprint,
2016.

 \noindent [Leb]: J. Lebl, Extension of Levi-flat
hypersurfaces past CR boundaries, {\it Indiana University
Mathematical Journal}, 57, no. 2, 699-716, 2008.

\noindent [LMSS]: J. Lebl, A. Minor, R. Shroff, Duong Son, and Yuan
Zhang, CR singular images of generic submanifolds under holomorphic
maps, {\it Arkiv Fur Matematik}, 52, no. 2, 301-327, 2014.

\noindent [LNR1]: J. Lebl, A. Noell and S. Ravisankar, Codimension
two CR singula submanifolds and extensions of CR functions, Arxiv
preprint (arXiv:1604.02073v2), 2016. (16 pages).

\noindent [LNR2]: J. Lebl, A. Noell, and S. Ravisankar, Extension of
CR functions from boundaries in ${\mathbb C}^n\times {\mathbb R}$,
arXiv:1505.05255, to appear in Indiana University Mathematics
Journal.

\noindent [Mos]: J. Moser, Analytic surfaces in ${\CC}^2$ and their
local hull of holomorphy,
 {\it Annales Aca -demi\ae Fennicae}, Series A.I. Mathematica
10, 397-410, 1985.

\noindent [MW]: J. Moser and S. Webster, Normal forms for real
surfaces in ${\CC}^2$ near complex tangents and hyperbolic surface
transformations,
 {\it Acta Math. 150}, 255-296, 1983.

\noindent [SM]: C. L. Siegel and J. Moser, Lectures on Celestial
Mechanics, {\it Classics in Mathematics}, Springer-verlag, 1995.

\noindent [Tu]: A. Tumanov, Extension of CR functions into a wedge
from a manifold of finite type, {\it Mat. Sb. 136} 128-139, 1988;
English transl. in {\it Math. USSR-Sb. 64}, 129-140, 1989.

\noindent [Zat]: D. Zaitsev, Private communication to X. Huang
during his visit to Trinity College, Dublin, 2013.

\bigskip
\noindent H. Fang and X. Huang, Department of Mathematics, Rutgers
University, New Brunswick, NJ 08903, USA.
(hf115$@$math.rutgers.edu,\ huangx$@$math.rutgers.edu)

\end{document}